%% file: Art-VEM-OptimalityAdaptivity.tex
\documentclass[11pt,a4paper]{elsarticle}
\usepackage{./Styles/ART-DFN-Adapt}
 \usepackage{amsmath}
\input{Styles/personal_packages}

\newtheorem{theorem}{Theorem}
\newtheorem{definition}{Definition}
\newtheorem{remark}{Remark}

\graphicspath{{./}}

\begin{document}
\begin{frontmatter}
  \title{A new quality preserving polygonal mesh refinement algorithm for Virtual Element Methods.\tnoteref{t1}}%
 \tnotetext[t1]{This research has been supported by INdAM-GNCS, by the MIUR projects PRIN 201744KLJL\_004 and 20204LN5N5\_003,
   and by the MIUR project ``Dipartimenti di Eccellenza 2018-2022''    (CUP E11G18000350001).
   Computational resources were partially provided by
   HPC@POLITO (\url{http://hpc.polito.it})
   {and SmartData@POLITO (\url{http://smartdata.polito.it})}.}%

\author[poli,indam]{Stefano Berrone\corref{cor1}}%
\ead{stefano.berrone@polito.it}%
\author[poli,indam]{Alessandro D'Auria}%
\ead{alessandro.dauria@polito.it}%
\address[poli]{Dipartimento di Scienze Matematiche, Politecnico di
  Torino\\Corso Duca degli Abruzzi 24, Torino, 10129, Italy}%
\address[indam]{Member of the INdAM research group GNCS}%
\cortext[cor1]{Corresponding author}%
 \begin{abstract}
Mesh adaptivity is a useful tool for efficient solution to partial differential equations in very complex geometries. In the present paper we discuss the use of polygonal mesh refinement in order to tackle two common issues: first, adaptively refine a provided good quality polygonal mesh preserving quality, second, improve the quality of a coarse poor quality polygonal mesh during the refinement process on very complex domains. For finite element methods and triangular meshes, convergence of a posteriori mesh refinement algorithms and optimality properties have been widely investigated, whereas convergence and optimality are still open problems for polygonal adaptive methods. In this article, we propose a new refinement method for convex cells with the aim of introducing some properties useful to tackle convergence and optimality for adaptive methods. The key issues in refining convex general polygons are: a refinement dependent only on the marked cells for refinement at each refinement step; a partial quality improvement, or, at least, a non degenerate quality of the mesh during the refinement iterations; a bound of the number of unknowns of the discrete problem with respect to the number of the cells in the mesh. Although these properties are quite common for refinement algorithms of triangular meshes, these issues are still open problems for polygonal meshes.
\end{abstract}

 \begin{keyword}
   Mesh adaptivity
   \sep Polygonal mesh refinement
   \sep Virtual Element Method
   \sep Convergence and Optimality

   \MSC[2010] 65N30 \sep 65N50
 \end{keyword}

\end{frontmatter}

\today

\input{Introduction} 
\input{VemApost}
\input{Definitions}
\input{Algorithms}
\input{num_res}
\input{num_res_Polygon} 
\input{conclusions}

\bibliographystyle{elsarticle-num}
\bibliography{scico2mp,dfn,articolo,VEMbibl,Apost_Meshadap}

\end{document}

%% file: Styles/personal_packages.tex
\usepackage{xcolor}
\usepackage{listings}

\usepackage{etex}
\usepackage{pictexwd}
\usepackage{etoolbox}

{\left\lbrace\begin{array}{@{}l@{}}}%
{\end{array}\right.}

\definecolor{CommentGreen}{rgb}{0.16,0.55,0.20}
\definecolor{WordBlue}{rgb}{0.08,0.28,0.6}
\definecolor{DefBlue}{rgb}{0.09,0.31,0.8}

\graphicspath{{Immagini/}}

%% file: Introduction.tex
\section{Introduction}
\label{sec:intro}
In recent years, the use of numerical methods on polygonal meshes and elements with curved edges has considerably reduced  the constraints for mesh generation in domains with high geometrical complexities. The versatility and robustness of these methods has increased their use in engineering applications in the field of mechanics, thermo and fluid dynamic, and in other engineering fields. A great interest in geophysical application has been addressed to these methods, among them underground flow problems and flows in poro-fractured media, where the geometric complexities of the domain sometimes make infeasible the creation of standard triangular meshes. Moreover, many of these applications are quite critical for safety reasons and the need of reliable simulations calls for an adaptive approach. Adaptive methods have been widely investigated both from the analytical and geometric point of view for standard finite elements \cite{Verfuerth1996,Doerfler,CasconKreuzeretal2008,Stevenson2005,Stevenson2007OptimalityOA}.
In recent years,  isotropic and anisotropic estimators have been developed and analysed  for polygonal methods \cite{BM:2008,BM:2015,Beirao-Manzini-Mascotto:2019,Cangiani_et_al-apost:2017,BBapost,ABBDVW}. The derived error estimators have many common properties with the estimators for standard finite elements, but a topic still under-investigated for polygonal methods is the refinement strategy preserving quality and with a bounded growth of unknowns.
The presence in the mesh of different type of elements that can be refined independently largely increase the difficulties of obtaining a good quality mesh \cite{Beirao2017,sorgente2021role,BBorth}.
A preliminary work for general convex polygons is recently appeared \cite{BBD}. 

Here we show an evolution of the algorithms presented in \cite{BBD} that takes into account the geometric properties of the cells such that the refinement results to preserve or improve the quality of the cell. The resulting algorithms is capable to ensure some relevant issue for an optimal refinement strategy: a refinement dependent only on the marked cells at each refinement step; a common improvement of the mesh quality during the refinement iterations; a bound of the number of unknowns of the discrete problem with the number of the cells in the mesh. Moreover, the refinement algorithm has proved to automatically converge towards good quality triangular or quadrilateral meshes very suitable for the solution in the regions subject to several refinement iterations.

The article is structured in this way. In Section~\ref{sec:VEM} a very brief introduction to the Virtual Element Method for 2-dimensional problems is provided with the geometrical assumptions on the polygonal elements. In Section~\ref{sec:Apost} the residual estimator used for the standard Poisson problem is recalled.
Virtual Element methods are introduced because the algorithm is tested in this framework, although none of the properties of these methods is assumed in the refinement algorithm that is exclusively based on the geometrical properties of the cells.
In Section~\ref{sec:Definitions} the geometric parameters used to check the quality of the mesh are provided.  In Section~\ref{sec:RefAlgo} the new refinement algorithm based on two different splitting directions is described. In the Section~\ref{sec:NumRes} we show the numerical tests solving the Poisson problem in the L-shape domain starting from a set of different initial meshes. In Section~\ref{sec:UnifPolyg} some other test problems concerning uniform refinement of regular polygons that easily yields to problematic configurations is presented.


%% file: VemApost.tex
\section{Model Problem}
\label{model}
In the present work
we consider the numerical solution with a polygonal mesh refinement method of the simple Poisson problem with piecewise constant diffusivity coefficient $\K$.
Let $\Omega\subset\mathbb{R}^2$ be a bounded open set with Lipschitz boundary $\partial\Omega$;
then, for a forcing term $f$ we look for a
function $u$ such that
\begin{equation}
  \label{eq:problem}
  \begin{cases}
    -\div\left(\K\nabla u\right) = f & \text{in $\Omega$}\,,
    \\
    u=0 & \text{on $\partial\Omega$}\,,
  \end{cases}
\end{equation}
where $\K$ is a positive function representing the \emph{diffusivity}.

Let
$\a{}{}\colon \sobho{1}{0}{\Omega}\times\sobho{1}{0}{\Omega}\to \mathbb{R}$ be
such that
\begin{equation*}
  \label{eq:defa}
  \a{w}{v} \defeq \scal{\K\nabla w}{\nabla v}\quad \forall w,v \in \sobho{1}{0}{\Omega} \,,
\end{equation*}
where $\scal{\cdot}{\cdot}$ is the $\lebl{\Omega}$ scalar product, and
given $f\in\lebl[2]{\Omega}$, the variational form of \eqref{eq:problem} is:
find $u\in\sobho{1}{0}{\Omega}$ such that
\begin{equation}
  \label{eq:exvarform}
  \a{u}{v} = \scal{f}{v}\quad \forall v\in\sobho{1}{0}{\Omega} \,.
\end{equation}

\section{The Virtual Element Method}
\label{sec:VEM}

Let $\Th$ be a discretization of $\Omega\subset \mathbb{R}^2$ with
open star-shaped polygons having an arbitrary, bounded, number of sides (even
different from one polygon to another) and let $\Eh$ be the set of
their edges.  As a mesh regularity assumption, we require that
$\forall E\in\Th$, with diameter $H_E$, there exists a constant
$\gamma>0$ such that
\begin{itemize}
\item $E$ is star-shaped with respect to a ball $B_E$ of radius larger
  than $\gamma H_E$;
\item for any two vertices $\mathbf{x}_1,\mathbf{x}_2\in E$,
  $\norm[\mathbb{R}^2]{\mathbf{x}_1-\mathbf{x}_2}\geq\gamma H_E$.
\end{itemize}
In \cite{Beirao2017} a detailed discussion on the mesh requirements and the VEM stability properties is provided.
The assumptions introduced imply the existence, on each
element $E\in\Th$, of a uniformly shape regular nested
triangulation $\Th[E]$ whose triangles $t$ are such that
\begin{equation}
  \label{eq:internal_tri_regular}
  \forall E\in\Th,\,\forall t\in\Th[E],\quad H_E\geq H_t\geq \gamma H_E \,,
\end{equation}
and each of these triangles have one edge lying on $\partial E$. This
subtriangulation can be accomplished by connecting all vertices of $E$ to
the center of the ball $B_E$, whose coordinates are
$\mathbf{x}_E=\left(x_E,y_E\right)$, although this option does not always provide the best quality subtriangulation.

Let $k\in\mathbb{N}$ be the polynomial order of the VEM
discretization, and let
$\Pi^\nabla_k\colon\sobho{1}{0}{\Omega}\to \Poly{k}{\Th}$ be the operator
such that, $\forall v\in\sobho{1}{0}{\Omega}$ and $\forall E\in\Th$,
\begin{equation*}
  \label{eq:defPinabla}
  \scal[E]{\nabla\left(v-\Pi_k^\nabla v\right)}{\nabla p}=0 \,, \forall
  p\in\Poly{k}{E} \ \text{and}\ 
  \begin{cases}
    \scal[\partial E]{\Pi^\nabla_k v}{1} = \scal[\partial E]{v}{1} &
    \text{if $k=1$}\,,
    \\
    \scal[E]{\Pi^\nabla_k v}{1} = \scal[E]{v}{1} & \text{if
      $k\geq 1$}\,,
  \end{cases}
\end{equation*}
where $\Poly{k}{\omega}$ is the space of the polynomials of degree less than or equal to $k$ on $\omega\subset\mathbb{R}^2$, \cite{Beirao2013a}.
Following \cite{Ahmad2013}, we introduce the finite
dimensional spaces
\begin{align*}
  \label{eq:defVhE}
  \begin{split}
    V_{h}^E &= \left\{ v\in \sobh{1}{E}\colon\Delta
      v\in\Poly{k}{E},\,v\in\Poly{k}{e} \,\forall e\subset\partial E,
      v_{|_{\partial E}} \in \cont{\partial E}\,, \right.
    \\
    &\quad\left. \scal[E]{v}{p}=\scal[E]{\Pi^\nabla_k v}{p}\,\forall
      p\in\Poly{k}{E}/ \Poly{k-2}{E} \right\}, \quad\forall E\in\Th,
  \end{split}
  \\
  V_\delta &= 
        \left\{
        v\in C^0(\Omega)\cap\sobh[0]{1}{\Omega} \colon v\in V_\delta^E \;\forall E\in\Th
        \right\},
\end{align*}
where $\Poly{k}{E}/ \Poly{k-2}{E}$ denotes the subspace
  of $\Poly{k}{E}$ containing polynomials that are
  $\lebl{E}$-orthogonal to $\Poly{k-2}{E}$ (see
  \cite{Beirao2015b}; other options are possible, see, for
  example, \cite{Ahmad2013}).
  
  A function $v\in V_\delta$ can be described on each polygon $E\in\Th$ by
  the following degrees of freedom:
  \begin{enumerate}
  \item the values at the vertices of the polygon;
  \item if $k\geq 2$, for each edge $e\subset\partial E$, the value of
    $v$ at $k-1$ internal points of $e$. For practical purposes, we
    choose these points to be the internal Gauss -- Lobatto quadrature
    nodes;
  \item if $k\geq 2$, the moments
    $\scal[E]{v}{m_{\boldsymbol{\alpha}}}$ for all the monomials
    $m_{\boldsymbol{\alpha}}\in\M_{k-2}(E)$
    up to
    the order $k-2$ , with
    $\boldsymbol{\alpha}=\left(\alpha_1,\alpha_2\right)$,
    $\abs{\boldsymbol{\alpha}}=\alpha_1+\alpha_2\leq k-2$, and
    \begin{equation}
      \label{eq:defmonomials}
      \forall \mathbf{x}=(x,y) \in E,\quad m_{\boldsymbol{\alpha}}(x,y) \defeq
      \frac{(x-x_E)^{\alpha_1}(y-y_E)^{\alpha_2}}
      {h_E^{\alpha_1+\alpha_2}}\,.
    \end{equation}
  \end{enumerate}
We point out that the chosen degrees of freedom uniquely identify
the polynomial expression of a function in $V_\delta$ on each edge of the
discretization, whereas inside the polygons these functions can not be
directly evaluated. These degrees of freedom ensure the computability of the projection $\Pi_k^\nabla v_\delta$, for any $v_\delta\in V_\delta$, see
\cite{Beirao2013a,Beirao2015b}, and, once it is known, to compute the
$\lebl{\Omega}$ projection of $v_\delta$ on $\Poly{k}{\Th}$, which is
indicated by $\Pi^0_k v_\delta$ in the following. Similarly,
$\Pi^0_{k-1}\nabla v_\delta$ indicates the vector containing the
$\lebl{\Omega}$ projection on $\Poly{k-1}{\Th}$ of the partial
derivatives of $v_\delta$, computable by the degrees of
freedom.

To introduce the VEM discretization of the Poisson problem we suppose
to know, for each $E\in\Th$, a symmetric bilinear form
$\vemstab[E]{}{}\colon V_\delta\times V_\delta\to \mathbb{R}$ that scales like
$\aE{}{}$ on the kernel of $\Pi_k^\nabla$, \cite{Beirao2017}, i.e.
$\exists c_\ast,c^\ast >0$ such that
\begin{equation}
  \label{eq:SEscales}
  \forall v_\delta\in V_\delta\text{ with }\Pi^\nabla_k v_\delta=0,\quad c_\ast\aE{v_\delta}{v_\delta} \leq \vemstab[E]{v_\delta}{v_\delta}\leq c^\ast\aE{v_\delta}{v_\delta} \,,
\end{equation}
where $\aE{v}{w}\defeq \scal[E]{\K\nabla v}{\nabla w}$. Once
$\vemstab[E]{}{}$ is given, we can define the following local and
global discrete bilinear forms:
\begin{gather*}
  \label{eq:defahE}
  \begin{split}
    \forall E\in\Th,\forall u_\delta,v_\delta\in
    V_\delta,\quad\ahE{u_\delta}{v_\delta}&\defeq\scal[E]{\K\Pi_{k-1}^0\nabla
      u_\delta}{\Pi^0_{k-1} \nabla v_\delta}
    \\
    &\quad +
    \vemstab[E]{\left(I-\Pi^\nabla_k\right)u_\delta}{\left(I-\Pi^\nabla_k\right)v_\delta}
    \,,
  \end{split}
  \\
  \label{eq:defah}
  \forall u_\delta,v_\delta\in
  V_\delta,\quad\a[h]{u_\delta}{v_\delta}\defeq\sum_{E\in\Th}\ahE{u_\delta}{v_\delta} \,.
\end{gather*}
With the above definitions, we can formulate the Virtual Element
method as the solution to the following discrete problem: find
$u_\delta\in V_\delta$ such that
\begin{equation}
  \label{eq:vemvarform}
  \a[h]{u_\delta}{v_\delta} = \scal{f_\delta}{v_\delta} \quad \forall v_\delta\in V_\delta \,,
\end{equation}
where $f_\delta\defeq\Pi^0_k f$, that is the best approximation of $f$ that
allows the computability of the scalar product with a VEM function,
since $\scal{\Pi^0_k f}{v_\delta} = \scal{f}{\Pi^0_k v_\delta}$ and we can
compute $\Pi^0_k v_\delta$ using the degrees of freedom. The well-posedness
of this problem simply follows by noticing that, thanks to
\eqref{eq:SEscales}, $\a[h]{}{}$ is coercive on $V_\delta$; optimal orders
of convergence are proved in \cite{Beirao2015b}.

One possible choice for $\vemstab[E]{}{}$ is the scalar product between the two
  vectors containing the degrees of freedom of the two functions
  involved, i.e., if we indicate by $\chi_r$ the operator which
  associates to each function in $V_\delta$ its $r$-th degree of freedom,
  \begin{equation}
    \label{original_stab}
    \vemstab[E]{u_\delta}{v_\delta} \defeq \K_E \sum_{r=1}^{N_E^{\rm dof}}\chi_r\left(u_\delta\right)\chi_r\left(v_\delta\right) \quad\forall E\in\Th,\forall u_\delta,v_\delta\in V_\delta\,,
  \end{equation}
  where $N_E^{\rm dof}$ indicates the number of degrees of freedom on element
  $E$.

The stabilization term can be avoided resorting to a slightly different VEM formulation \cite{berrone2021lowest}.

\section{A posteriori error estimate}
\label{sec:Apost}
The proposed refinement strategies are tested on a simple Poisson problem with constant diffusivity parameter. This simple case fits the assumptions made in \cite{BBapost} to prove an equivalent relation between the error $e_\delta^\pi=u-u_\delta^\pi$, where $u_\delta^\pi=\Pi^\nabla_k u_\delta$, measured in the broken norm
\begin{equation}
  \label{eq:defennorm}
  \ennorm{v} \defeq \sup_{w\in\sobho{1}{0}{\Omega}}\frac{\sum_{E\in\Th}   \scal[E]{\K\nabla v}{\nabla w}}{\norm{\sqrt{\K}\nabla w}}
  \,
\end{equation}
and the computable, residual-based error estimator introduced in the following. The results provided in \cite{BM:2015,Beirao-Manzini-Mascotto:2019,Cangiani_et_al-apost:2017} can be applied as well.

In order to define the error estimator mesh-size parameters are needed: let $H_E$ be the longest edge of the element $E\in\Th$, and $h_e$ be  the length of the edge $e\in\Eh$.

For each $e\in\Eh$, let $R,L\in\Th\,$ be the two polygons sharing $e$
and $\mathcal{N}_{e}\defeq\{R,L\}$ the set of the elements (right and left) sharing the edge $e$. 

\begin{definition}
  For any internal edge $e\in\Eh$ let us define a unit normal vector
  $\mathbf{n}_e$ as the outward unit normal vector for the element on
  the right of $e$ ($\mathbf{n}_e=\mathbf{n}_R$) and the jump of
  the co-normal derivative of $u_\delta^\pi$
  \begin{displaymath}
    \jmpnormder{u_\delta^\pi}{e}{\K} = \K\nabla{u_\delta^\pi}_{|_{R}}\cdot\mathbf{n}_{R}
    -\K\nabla{u_\delta^\pi}_{|_{L}}\cdot\mathbf{n}_{L}=\K\nabla{u_\delta^\pi}_{|_{R}}\cdot\mathbf{n}_{e}
    +\K\nabla{u_\delta^\pi}_{|_{L}}\cdot\mathbf{n}_{e}.
  \end{displaymath}
Moreover, let $\eta_{R,E}$ the local residual estimator $\forall E\in\Th$
\begin{displaymath}
  \eta_{R,E}^2
  \defeq
  \frac{H^2_E}{\K}\norm[E]{f_\delta+\div\left(\K \nabla u^\pi_\delta\right)}^2
  +
  \frac12\sum_{e\in\Eh\cap\partial E} \frac{h_e}{\K}
  \norm[e]{\jmpnormder{u_\delta^\pi}{e}{\K}}^2 \,,
\end{displaymath}
and finally, we define the error estimator for the computable projection $u_\delta^\pi$ of the VEM solution $u_\delta$: 
\begin{equation}
  \label{eq:defetaR}
  \eta_R \defeq \left\{\sum_{E\in\Th} \eta_{R,E}^2 \right\}^{\frac12} \,.
\end{equation}
\end{definition}

In the following, the relation $a\lesssim b$ means that there exists a constant $C$ independent of $a, b$ such that  $a\leq C\, b$.

\begin{theorem}\cite{BBapost}
  \label{teo:upper_bound}
  Let $u$ be the solution to \eqref{eq:exvarform}, $u_\delta^\pi=\Pi^\nabla_k u_\delta$, and $f_\delta = \Pi^0_k f$. Then,
  \begin{equation}
    \ennorm{u-u^\pi_\delta}\lesssim \left\{ \eta_R^2
      + \sum_{E\in\Th}\frac{H_E^2}{\K}\norm[E]{f-f_\delta}^2\right\}^{\frac12}
    \,,
  \end{equation}
  being the constant independent of the meshsize and dependent on the mesh quality and on the VEM order $k$.
\end{theorem}

\begin{theorem}\cite{BBapost}
  \label{teo:lower_bound}
  Let $u$ be the solution to \eqref{eq:exvarform}, $u_\delta^\pi=\Pi^\nabla_k u_\delta$, and  $f_\delta = \Pi^0_k f$. Then, 
  \begin{equation}
    \eta_{R} \lesssim 
    \left\{ \ennorm[]{u-u^\pi_\delta}^2
      +
      \sum_{E\in\Th}\left(1+H^2_{E}\right)
      \frac{H^2_{E}}{\K}
      \norm[E]{f-f_\delta}^2 \right\}^{\frac12} \,.
  \end{equation}
\end{theorem}

Results contained in Theorems~\ref{teo:upper_bound} and \ref{teo:lower_bound} provide an equivalent relation between the error estimator $\eta_R$ and the error $\ennorm[]{u-u^\pi_\delta}$ up to the approximation $\norm[E]{f-f_\delta}$ of the right hand side $f$, they are proved in \cite{BBapost} provided some assumption on the polygonal mesh.

In order to evaluate the quality, reliability and efficiency of the estimator the effectivity index $e.i.$ is commonly used:
\begin{equation}
  \label{eq:efeffectivity}
  e.i. \defeq \frac{err}{\eta_R} 
  \quad \text{where}\quad 
  err\defeq \left\{
    \sum_{E\in\Th}\norm[E]{\sqrt{\K}\nabla(u-u^\pi_\delta)}^2
  \right\}^{\frac12}\,.
\end{equation}

Being the effectivity index strongly influenced by several aspects of mesh quality, in the following we will investigate its behaviour during the refinement process and we will analyse how the different refinement strategies considered can impact on it.


%% file: Definitions.tex


\section{Definition of mesh quality parameters}
\label{sec:Definitions}
In order to evaluate and check the evolution of the mesh quality during the refinement process we introduce a set of geometric parameters based on the following geometric quantities, $\forall E\in\Th$ and $e\in\Eh$ let us define: 
\begin{itemize}
\item  $H_E$ the longest edge of the element $E$;
\item  $h_E$ the shortest edge of the element $E$;
\item  $h_e$ the length of the edge $e$;
\item $C_E$ the centroid of the element $E$;
\item $R_E$ the largest distance between the centroid and the vertices of  $E$;  
\item $r_E$ the smallest distance between the centroid and the edges of  $E$.
\end{itemize}

The quality of the mesh elements is analysed through  the following parameters:

\begin{itemize}
\item $AR^{Rr}_E=R_E/r_E$;
\item $AR^{Hr}_E=H_E/r_E$; 
\item $AR^{Hh}_E=H_E/h_E$;  
\item $AR^{edge}_E={\rm max}_{e=1}^{N_E}\frac{{\rm max}\{h_{e},h_{(e+1)\%N_E}\}}{{\rm min}\{h_{e},h_{(e+1)\%N_E}\}}$;
\item $K^{\Pi^\nabla}_E=\kappa_2(\Pi^\nabla_E)$, condition number of the matrix $\Pi^\nabla_E\in\mathbb{R}^{{\rm dim}(\Poly{k}{E})\times\#{\rm dofs}_E}$.
\end{itemize}

The matrix $\Pi^\nabla_E$ is computed using the scaled monomial basis for $\Poly{k}{E}$ introduce in Section~\ref{sec:VEM}, this quantity is relevant only for VEM discretization.

Some of these parameters provide similar information, and we will discuss which of them can be considered significant in order to discuss the quality of the elements of a polygonal partition for a VEM discretization. These parameters are selected because their computation does not require the solution of complex or expensive problems. For example, instead of using $R_E$ and $r_E$ the computation of the radius of the smallest circumscribed circle and of the largest inscribed circle (the ball $B_E$ for which the convex element $E$ is star shaped) should be more significant, but their computation for a generic convex polygon can be quite expensive. The quantities $R_E$ and $r_E$ can possibly replace these quantities in a cheaper way. The quality of the cells depends on the best sub-triangular mesh, that is non-trivially constructed as well.

Moreover, we consider some parameter to characterize the mesh in term of number of points, cells, and number of vertices of the polygons.
\begin{itemize}
\item $\#{\Th}$: number of cells of the mesh ${\Th}$;
\item $\#P_{\Th}$: number of points of the mesh ${\Th}$;
\item $\#\triangle_{\Th}$: number of cells of ${\Th}$ with three vertices;
\item $\#\diamond_{\Th}$: number of cells of ${\Th}$ with four vertices;
\item $E^{\#P}_{\#{\Th}}=\#P_{\Th} / \#{\Th}$;
\item $R^\triangle_{\Th}=\#\triangle_{\Th} / \#{\Th}$;
\item $R^\diamond_{\Th}=\#\diamond_{\Th} / \#{\Th}$.
\end{itemize}

In the following sometimes we will use the previous symbols $\#\triangle_{\Th}$, $\#\diamond_{\Th}$, $R^\triangle_{\Th}$, $R^\diamond_{\Th}$ dropping the mesh symbol $\Th$.

For the Virtual Element Methods, the parameter $E^{\#P}_{\#{\Th}}$ can be considered as a sort of cost parameter in the sense that each cell is providing an approximation to the solution with a fixed polynomial degree, the degree of the VEM elements, at a cost that is increasing with the number of the points of the polygon. Considering this aspect, the most efficient mesh for approximating the solution with a piecewise polynomial function (the function $u_\delta^\pi$ in the VEM a posteriori error estimator) is the triangular mesh displaying the minimum number of ${\rm dofs}$ per cell required by the polynomial degree $k$ and for which $E^{\#P}_{\#{\Th}}\simeq 0.5$.


%% file: Algorithms.tex
\section{Refinement algorithms}
\label{sec:RefAlgo}
In this section we introduce the refinement algorithms used for refining marked convex cells generating two new convex sub-cells. 
In the algorithms we assume that each edge of the original mesh is provided with a marker uniquely identifying the edge, it is inherited by the children edges generated by the refining algorithm. This marker is used in order to identify aligned edges of a cell produced by the refinement of some neighbouring cell without any additional computational cost. We do not consider aligned edges of the initial mesh present in the same element.

Triangular cells are always refined with the longest edge refinement criterion. This refining criterion is also applied when a cell is recognised as a triangle in spite of the presence of more than three vertices when aligned edges are present. For other polygons the first refinement algorithm presented is based on the idea introduced in the article \cite{BBD} of using the maximum inertia moment as cutting direction of the element (Algorithm~\ref{alg:MaximumMomentum}). The second one, is a natural extension of the longest edge refinement applied to polygons (Algorithm~\ref{alg:LongestDiagonal}) considering the longest diagonal of the cell. This second option does not require the computation of integrals on the element but the (squared) distance between all the couple of vertices.

In the Algorithm \ref{alg:SmoothingDirection} we describe how we generate an optimal cutting direction. The algorithm consists of a double loop over the edges of the cell
to prevent the collapsing of the cutting line $Cl$ on an edge of the element when a change of $Cl$ is suitable (collapsing on a vertex).

In the Algorithm \ref{alg:Refinement} we describe the procedure to refine the polygonal cells with a cut along the inertia axis with the maximum moment improved with respect to  \cite{BBD}.
In order to prevent the creation of small edges in the refined cells a small change of the refining direction was proposed in \cite{BBD} such that new edges shorter than a fraction of the cut edge were collapsed.
In this new algorithm we define an adaptive collapsing tolerance that depends on the shape and the edges of the refining polygon $M$ and its neighbouring cells $\mathcal{N}_M$ .
For each edge $e\in\Eh$ we define the set of the neighbouring cells $\mathcal{N}_{e}$.
In order to preserve or improve the quality of the produced cells we define
$\rho_E = \min \lbrace h_E, r_E \rbrace$ for $E\in\Th$,
and then we define for each edge $e\in\Eh$
\begin{equation}
  \label{rho}
  \rho_{e} = \max \lbrace \rho_{E}\ :\ E\in\mathcal{N}_{e}\rbrace.
\end{equation}
Let $\Mh\subset\Th$ the set of the marked elements for refinement, let us consider $M\in\Mh$ and let $e\in\Eh[M]$ an edge of the element $M$ that should be refined for refining $M$ (the edge intersected by the cutting line $Cl$). 
The parameter $\rho_e$ is designed in order to prevent the formation of small edges in the refinement of the edge $e$ that could degenerate the quality of the mesh either reducing the quality of the children of $M$ or of the neighbouring element of $M$ faced to the edge $e$.
The refinement algorithm is based on four main rules:
\begin{enumerate}
\item New edges obtained refining the edge $e$ are created only if they are larger then $c_{\rho}\rho_e $ for some arbitrary parameter $c_\rho\ge 0$, otherwise the intersection is collapsed to the closest vertex.
\item Regardless of the intersection point between $Cl$ and $e$, the edges is cut in its mid-point if the intersection does not collapse in one of its vertices.
\item If the intersected edge is part of a set of aligned edges (descendant of an edge of the initial mesh) of the element $M$, the set of aligned edges of $M$ is cut in the mid-point of the set. We remark that due to the chosen procedure this point is always a point of $M$.
\item A cell with three edges or three sets of aligned edges (a geometric triangle, regardless of the number of vertices) is cut with the Longest Diagonal rule.
\end{enumerate}

  \begin{figure}[ht]
    \centering   
    \includegraphics[width=0.3\linewidth,bb= 1 0 348 160]{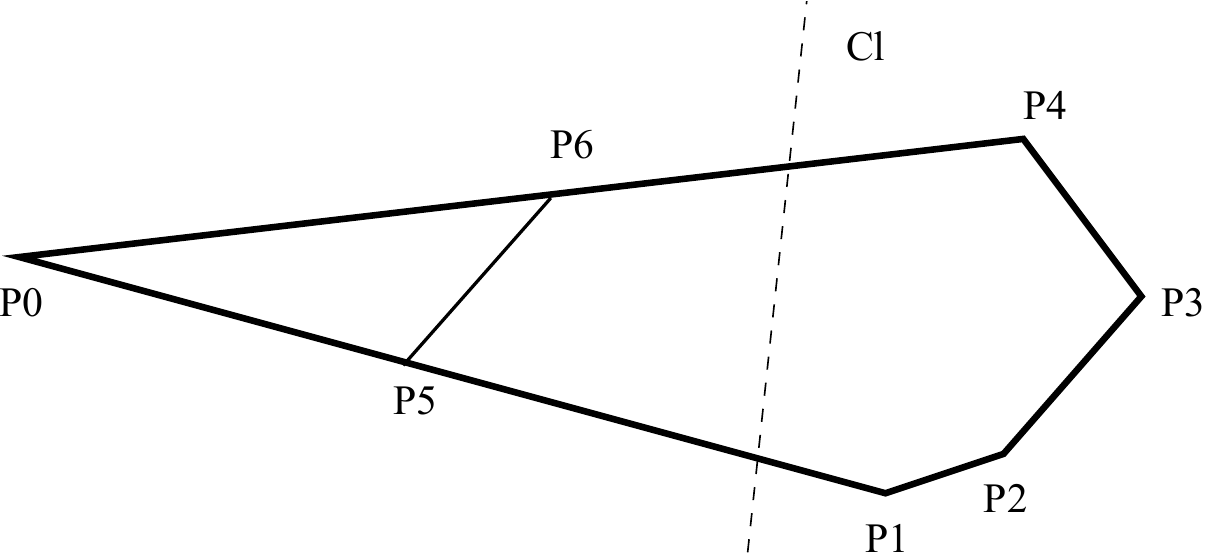}
        \caption{Example of polygon producing a child with more vertices, $M=P_0\!-\!P_6$, children $P_0P_5P_5$ and $P_6P_5P_1\!-\!P_4$.}
    \label{fig:MorePoints}
  \end{figure}

    \begin{figure}
    \centering
 \includegraphics[width=0.3\linewidth,bb= 65 20 353 296]{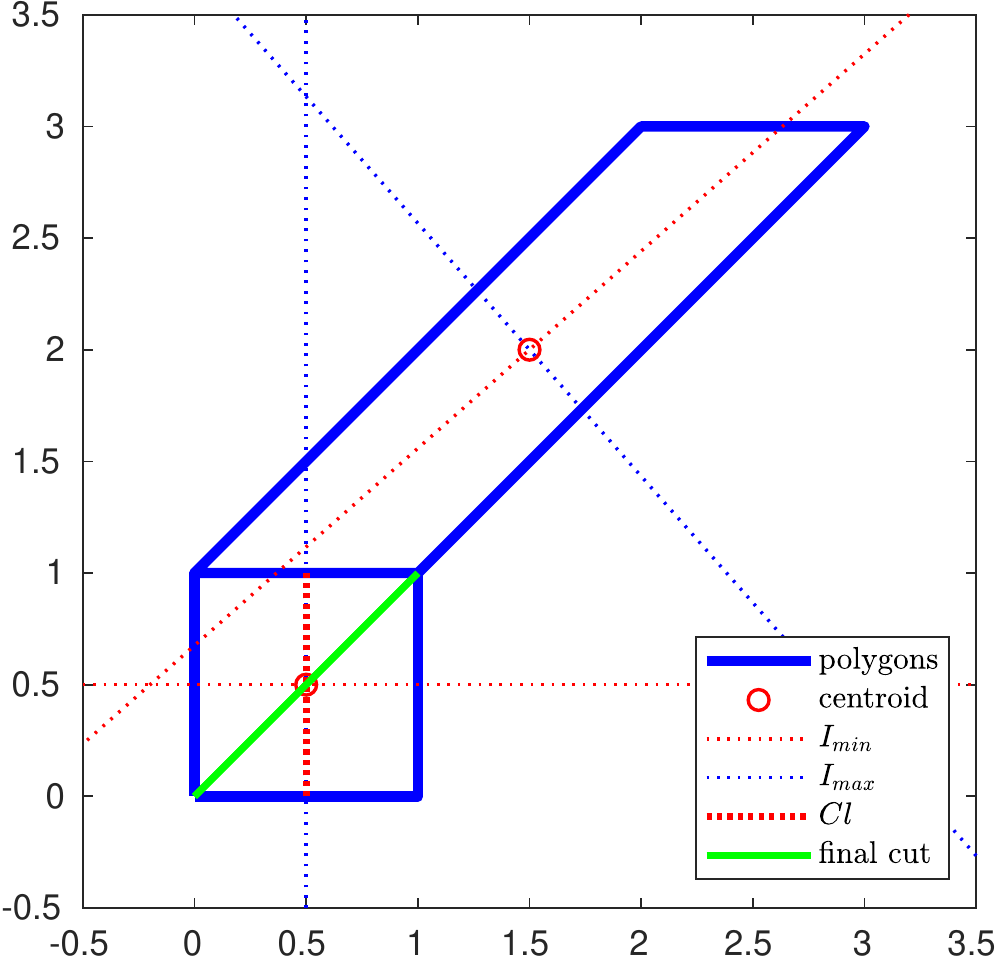}
    \caption{Example of refining conditioned by the neighbouring cell: bottom square $\rho_E=0.5$, upper trapezoid $\rho_E=0.5224$, the square is not cut by the vertical dotted red segment because this cut would produce an edge $|\hat e|/2 <0.5224$ not viable for the upper cell, so the cut is collapsed on the first vertex of the edge and a new $Cl$ is set passing from this vertex and the centroid. The new $Cl$ is passing through a vertex and the new cut results in the green continuous segment. This example is valid with $c_\rho \geq \frac{0.5}{0.5224}$ }
    \label{fig:test2}
  \end{figure}
  
\begin{algorithm}
    \caption{Compute Maximum Moment Direction}
    Given a convex polygon $M\in\Mh$
    \begin{algorithmic}[1]
		\STATE Compute the centroid $C_M$
		\STATE Compute the inertia tensor respect to $C_M$
		\STATE Compute eigenvalues and eigenvectors
		\STATE Set the cut-line $Cl$ to the line parallel to eigenvector corresponding to the maximum eigenvalue trough $C_M$
    \end{algorithmic}
    \label{alg:MaximumMomentum}
  \end{algorithm}

  \begin{remark}
    The main differences between the proposed algorithm and the one presented in \cite{BBD} is that in deciding if an edge can be cut we involve the effects of the cut on the neighbouring element sharing the edge. If the cut of an edge is suitable for the quality of the element that is performing the cut, but negative for the quality of the other facing element the cut is not performed and the intersection between the cutting line and the edge is collapsed to the closest vertex. Moreover, the cut, if performed, always occurs in the middle of the edge in this way an edge is cut always in the same way independently of the element that first processes it if the two sharing elements may cut the same edge. This simple improvement prevent some bad situations for the quality of the produced cells.
  \end{remark}
  
  
  \begin{algorithm}
    \caption{Compute Longest Diagonal Direction}
    Given a convex polygon $M\in\Mh$
    \begin{algorithmic}[1]
      \STATE Compute the longest diagonal $D_M$
      \STATE Set the cut-line $Cl$ to the line orthogonal to the longest diagonal trough $C_M$
    \end{algorithmic}
    \label{alg:LongestDiagonal}
  \end{algorithm}

  
    \begin{algorithm}
    \caption{Smoothing Direction}
    Given a convex polygon $M\in\Mh$ and the cut-line $Cl$
    \begin{algorithmic}[1]
      \STATE Define the set $\mathcal{N}_M$ of the elements sharing an edge with $M$
      \STATE Compute $\rho_E = \min \lbrace h_E, r_E \rbrace \qquad \forall E \in \mathcal{N}_M$
      \FOR{Each edges $e$ of the cell}
      \STATE Compute the intersection between $e$ and $Cl$
      \IF{There is an intersection}  \label{rho-eA}
      \STATE Set $\rho_e = \max \rho_E \qquad E \in \mathcal{N}_e$ 
      \IF{$e$ has some neighbour aligned edges of $M$ generated by a previous refinement}
      \STATE Find the mid-point of the ancestor largest edge $\hat e\supset e$ of $M$
      \ELSE
      \STATE Find the mid-point of $\hat e=e$ and set it as standard cutting point.
      \ENDIF
      \IF{$|\hat e|/2 <= c_{\rho} \rho_e$} \label{rho-eB}
      \STATE Set the closest vertex as candidate collapsed cutting point
      \ENDIF
      \ENDIF \label{collaps}
      \ENDFOR
      \smallskip
      \IF{We have two standard cutting points}
      \STATE Split the cell with the segment through the two standard cutting points
      \ELSIF{We have two non-consecutive candidate collapsing points}
      \STATE Split the cell with the segment through the two collapsed cutting points
      \ELSIF{We have a standard cutting point and a candidate collapsing point}
      \STATE Set the cutting line $Cl$ as the line passing through the collapsed point and the centroid.
      Compute the intersection between $Cl$ and the edges nonadjacent to the collapsed point. Repeat points \ref{rho-eA}-\ref{collaps}.
      \ELSIF{We have two consecutive candidate collapsing points}
      \STATE Set the cutting line $Cl$ as the line passing through the collapsed point closest to the original intersection between the two intersected edges and the centroid.
      Compute the intersection between $Cl$ and the edges nonadjacent to the collapsed point. Repeat points \ref{rho-eA}-\ref{collaps}.
      \ENDIF
    \end{algorithmic}
    \label{alg:SmoothingDirection}
  \end{algorithm}
  
\begin{algorithm}[t]
    \caption{Refinement algorithm}
    Given a convex polygon $M\in\Mh$
    \begin{algorithmic}[1]
      \IF{$M$ is a triangle}
      \STATE Refine the cell with the longest edge refinement
      \ELSE
      \STATE Compute the centroid $\mathbf{C}_M$
      \STATE Compute the cut-line $Cl$
	  \STATE Use the Smoothing Direction Algorithm \ref{alg:SmoothingDirection}
      \STATE Cut the cells in the two children cells
      \STATE Update the neighbourhood elements
      \ENDIF
    \end{algorithmic}
    \label{alg:Refinement}
  \end{algorithm}

  \begin{remark}
    The presented refinement algorithm can produce new children cells with at most the number of vertices of the father cell increased by one.
A child cell have more vertices than the father cell  when the cut refines to consecutive edges of the father cell, see Figure~\ref{fig:MorePoints}. In all the other cases the vertices of the children cells are less or equal to the number of vertices of the refined cell.
  \end{remark}

    \begin{remark}
    The absence of a conformity recovering step needed, for example, to get a conforming triangular mesh in the refinement of a triangular grid implies that the number of new cells generated by the refinement algorithm that splits each refined cell in two children cells is exactly the number of marked cells. Neighbouring cells of refined cells may change due to the refined edges shared with refined cells. In order to prevent the growing of number of vertices for non refined cells, after the refinement of all the marked cells a refinement of the un-refined cells modified with an increment of the number vertices larger than a given number or a refinement of the cells with many aligned edges can be considered. We have not implemented this facility and by our experiments we did not notice any problem concerning the growing of vertices per cell in the mesh.
    \end{remark}
    
    \begin{remark}
      Being the longest diagonal the maximum distance between two points of the cells, the Longest Diagonal criterion applied to a cell with a triangular shape independently of the number of vertices is equivalent to the longest edge refinement applied to the same cell but with glued aligned edges.
    \end{remark}



%% file: num_res.tex

\begin{figure}
  \centering
  \begin{subfigure}[b]{0.32\linewidth}
    \includegraphics[width=\linewidth]{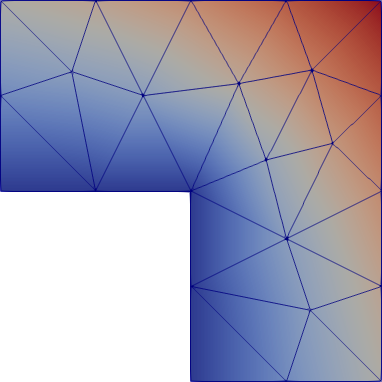}
  \end{subfigure}
  \hfill
  \begin{subfigure}[b]{0.32\linewidth}
    \includegraphics[width=\linewidth]{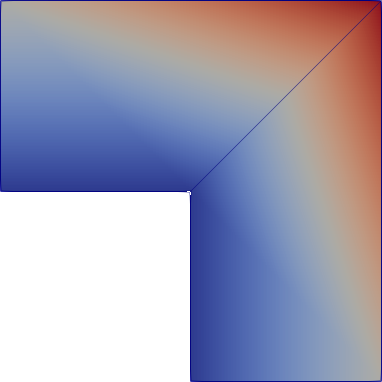}
  \end{subfigure}
  \hfill
  \centering
  \begin{subfigure}[b]{0.32\linewidth}
    \includegraphics[width=\linewidth]{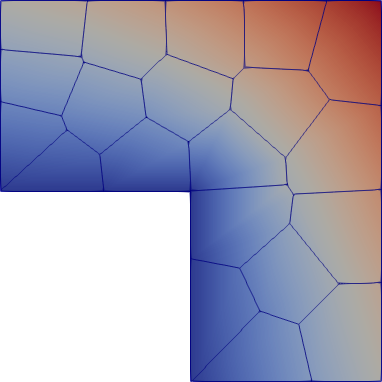}
  \end{subfigure}
  \caption{Starting Meshes}
  \label{fig:StartingMeshes}
\end{figure}

\begin{figure}
  \centering
  \begin{subfigure}[b]{0.42\linewidth}
    \includegraphics[width=\linewidth]{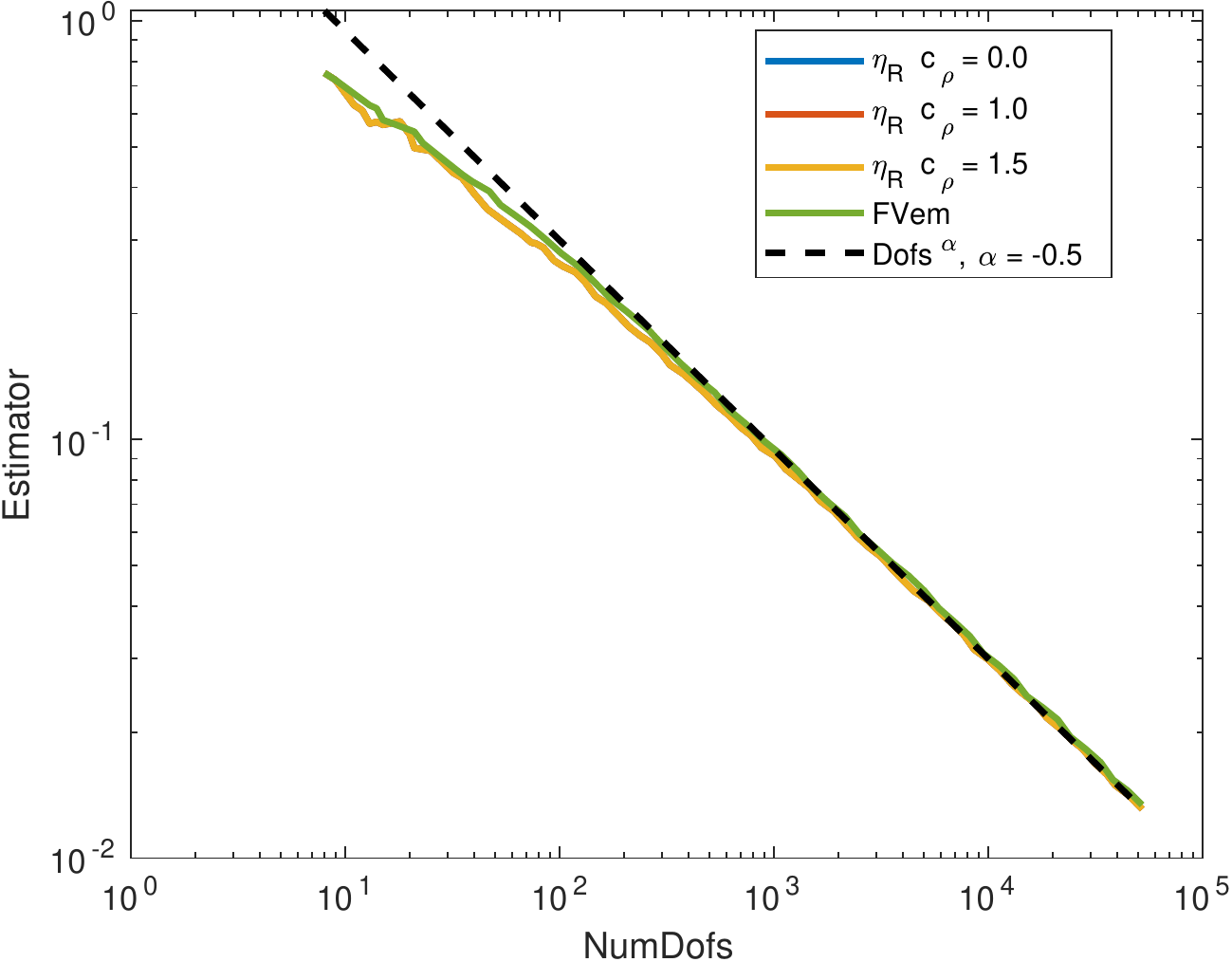}
    \caption{TRI\_MM-LE}
    \label{fig:TMVLF_EstDof}
  \end{subfigure}
  \hfill
  %
    \begin{subfigure}[b]{0.42\linewidth}
    \includegraphics[width=\linewidth]{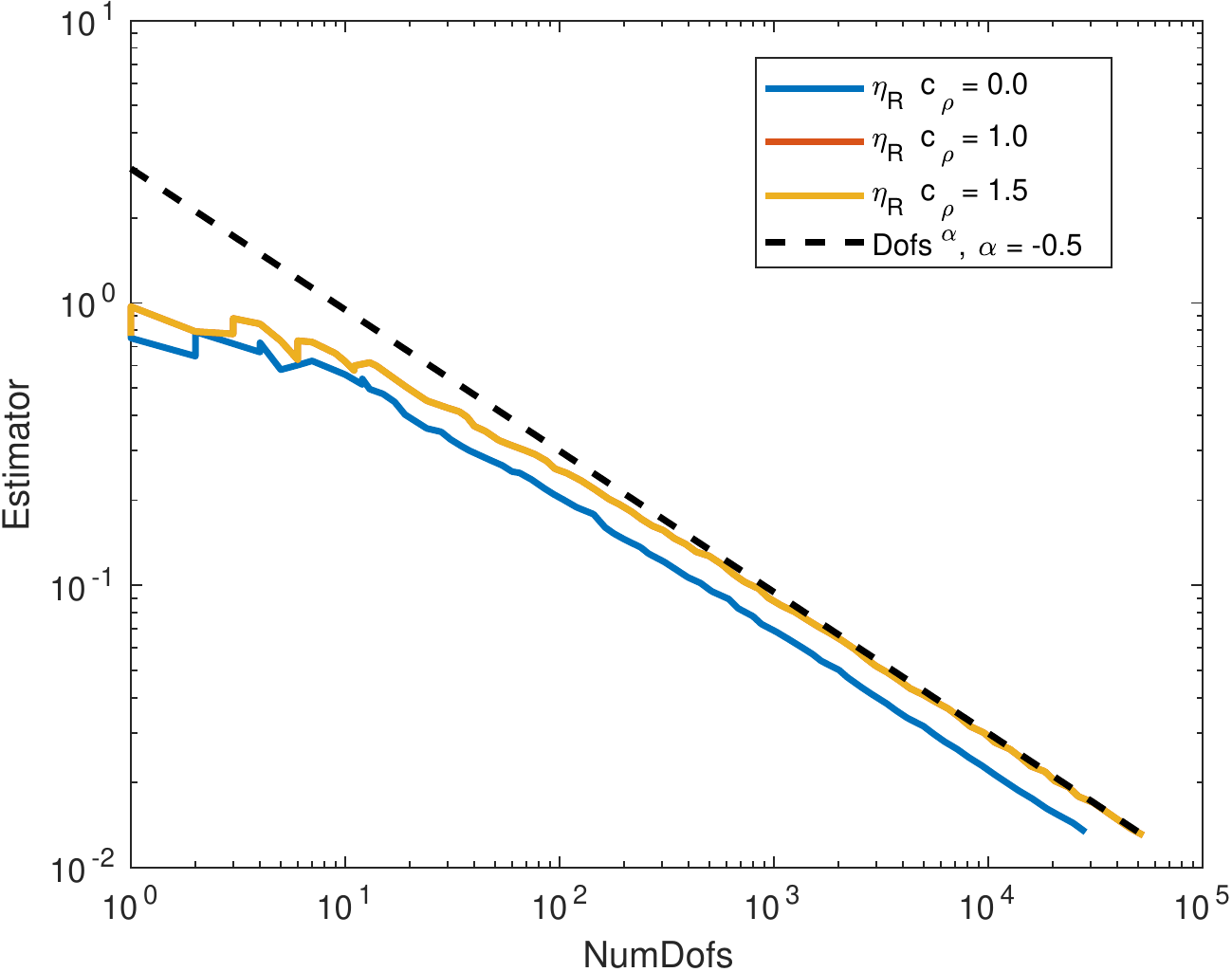}
    \caption{TRAP\_MM}
    \label{fig:TMV_EstDof}
  \end{subfigure}
  \hfill
  %
  \begin{subfigure}[b]{0.42\linewidth}
    \includegraphics[width=\linewidth]{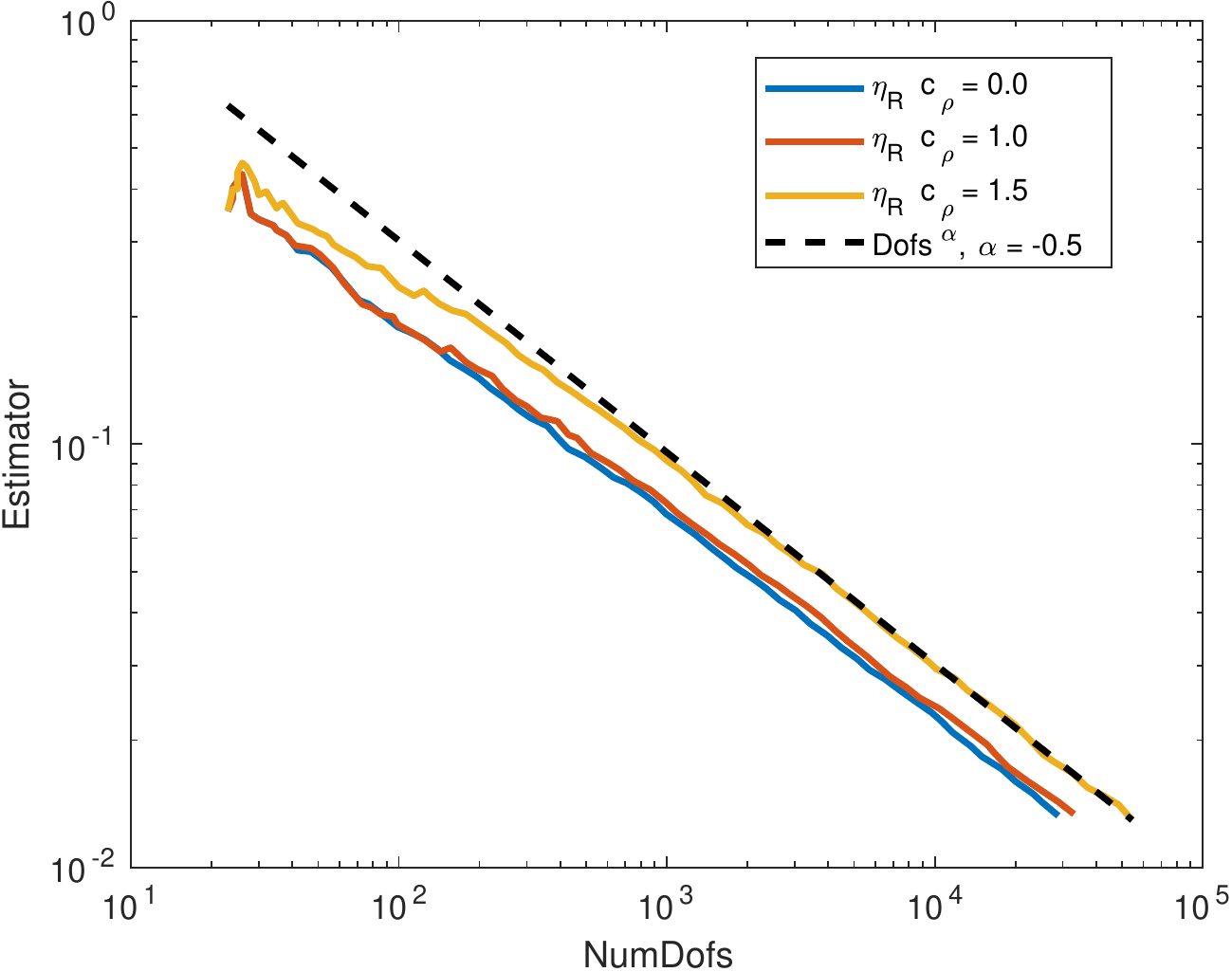}
    \caption{POLY\_MM}
    \label{fig:PMV_EstDof}
  \end{subfigure}
  \hfill
  \begin{subfigure}[b]{0.42\linewidth}
    \includegraphics[width=\linewidth]{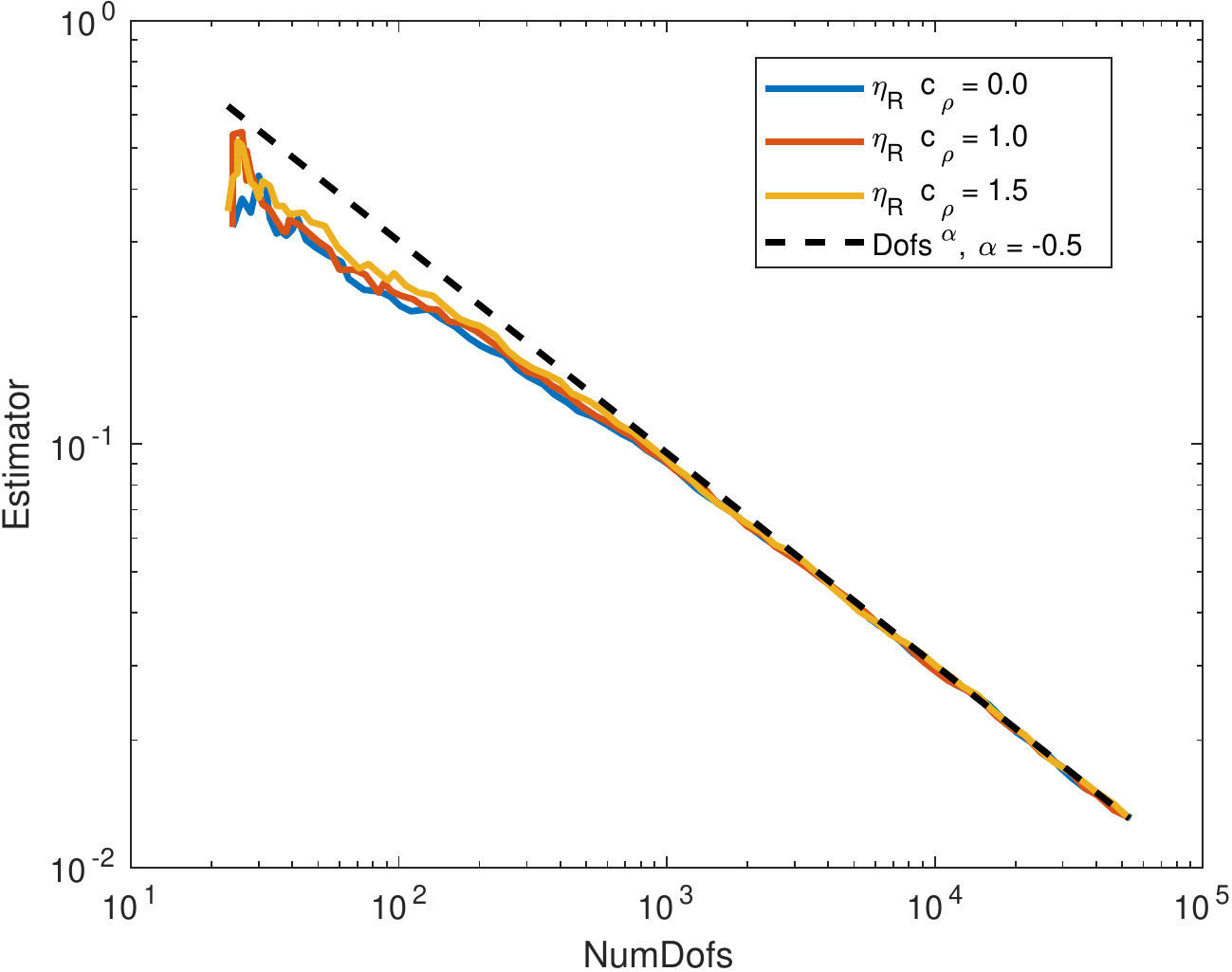}
    \caption{POLY\_LD}
    \label{fig:PLV_EstDof}
  \end{subfigure}
  \caption{Error Estimator $\eta_R$ versus \#dofs}
  \label{fig:EstDofs}
  \centering
  \begin{subfigure}[b]{0.42\linewidth}
    \includegraphics[width=\linewidth]{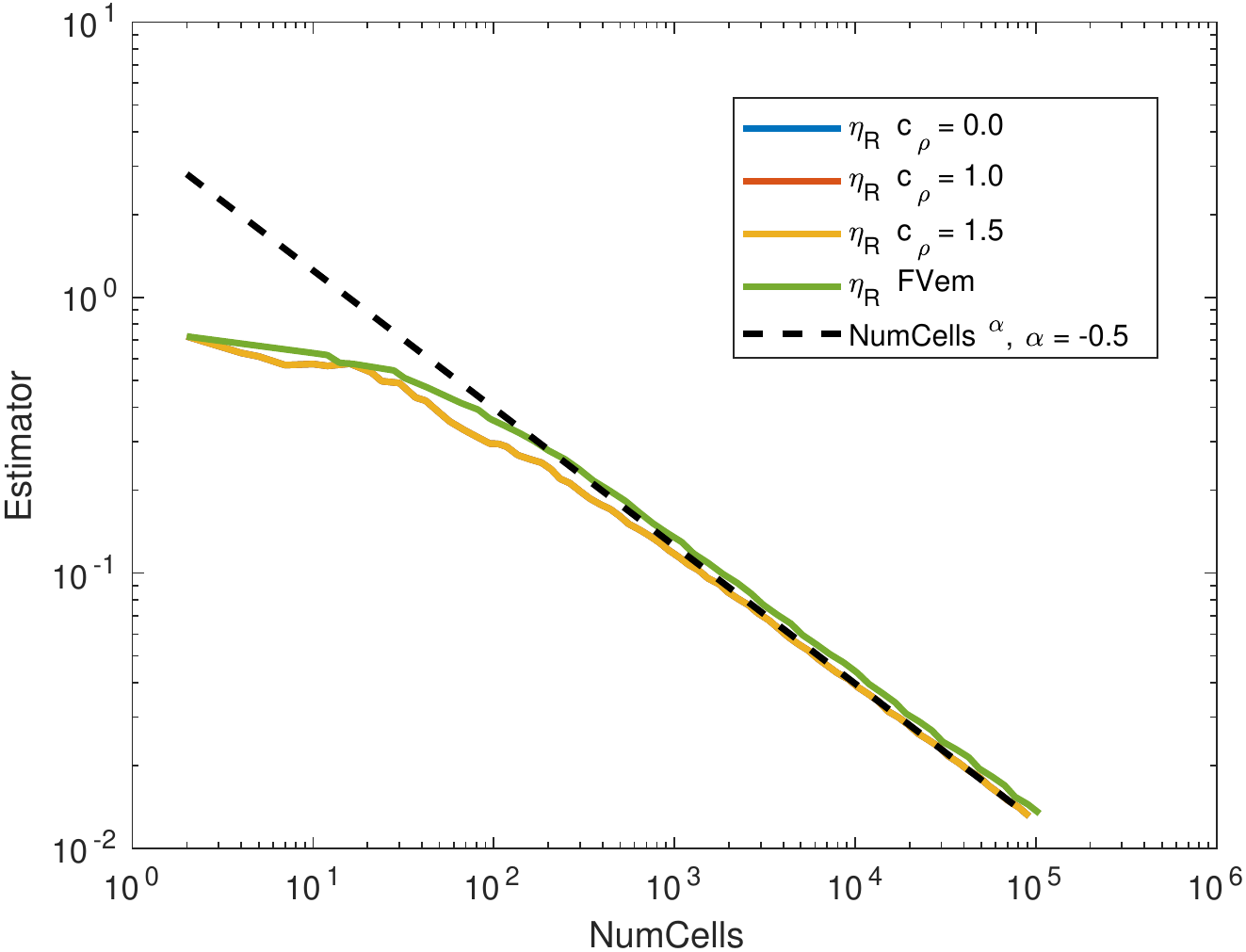}
    \caption{TRI\_MM-LE}
    \label{fig:TMVLF_EstCells}
  \end{subfigure}
  \hfill
  %
    \begin{subfigure}[b]{0.42\linewidth}
    \includegraphics[width=\linewidth]{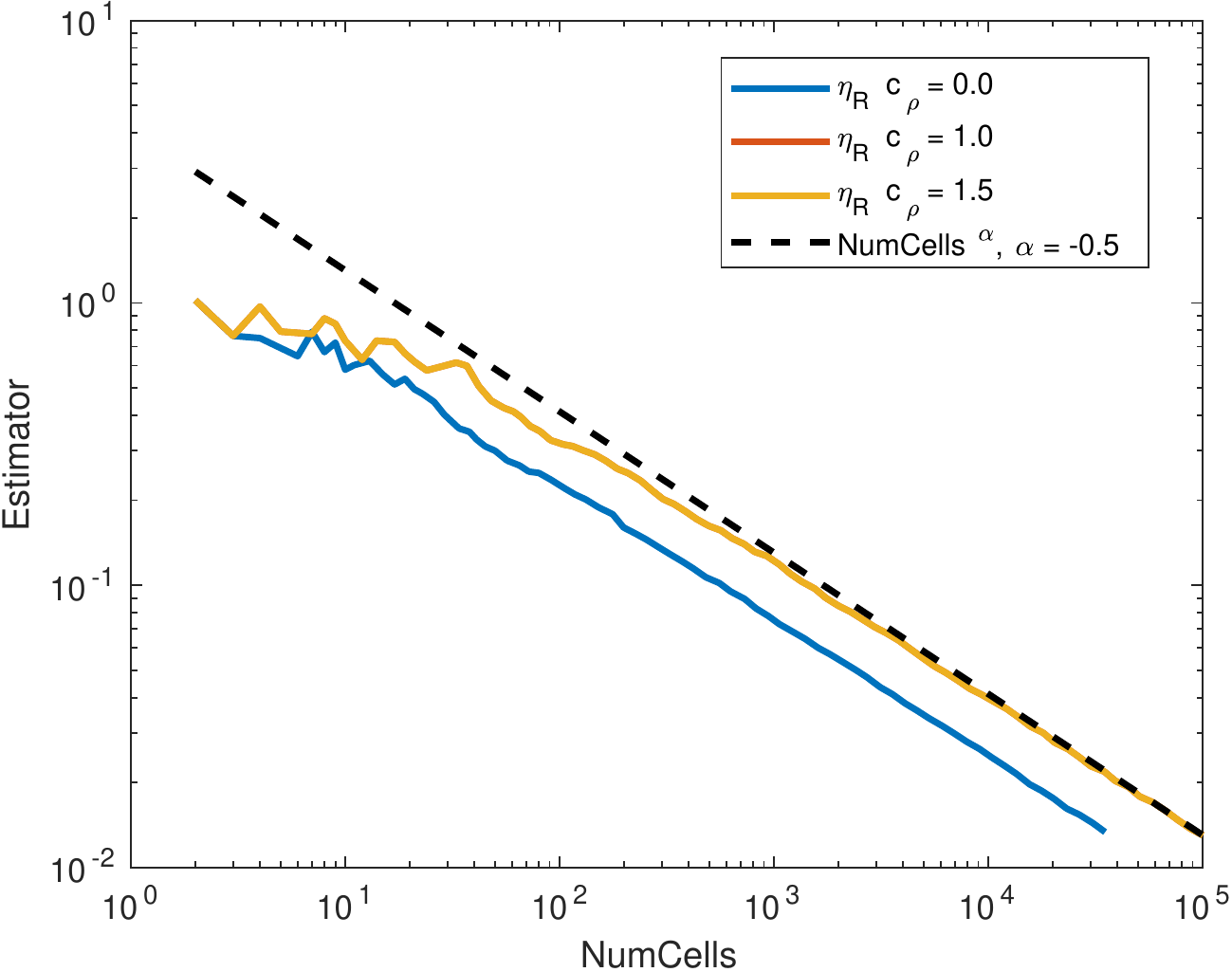}
    \caption{TRAP\_MM}
    \label{fig:TMV_EstCells}
  \end{subfigure}
  \hfill
  %
  \begin{subfigure}[b]{0.42\linewidth}
    \includegraphics[width=\linewidth]{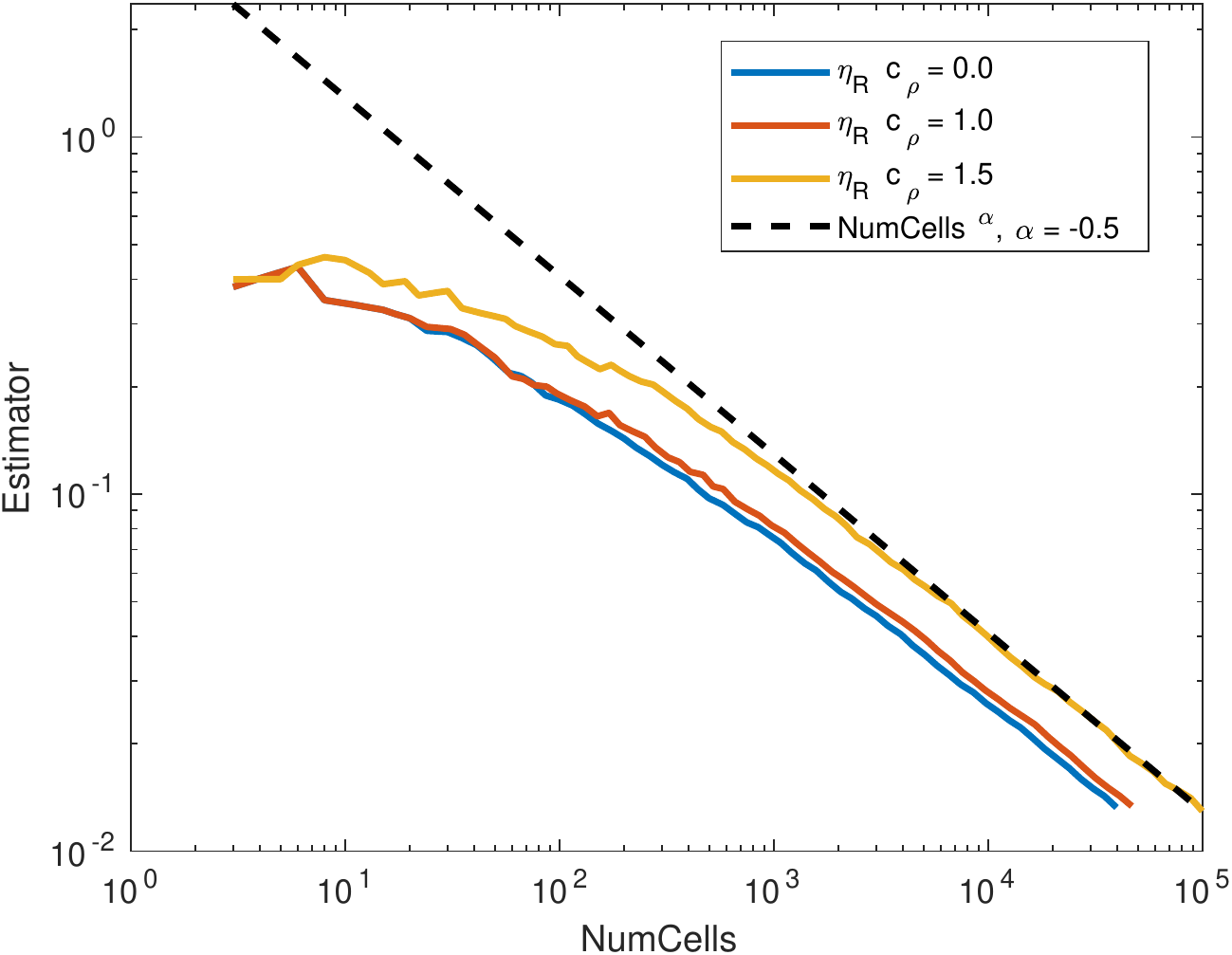}
    \caption{POLY\_MM}
    \label{fig:PMV_EstCells}
  \end{subfigure}
  \hfill
  \begin{subfigure}[b]{0.42\linewidth}
    \includegraphics[width=\linewidth]{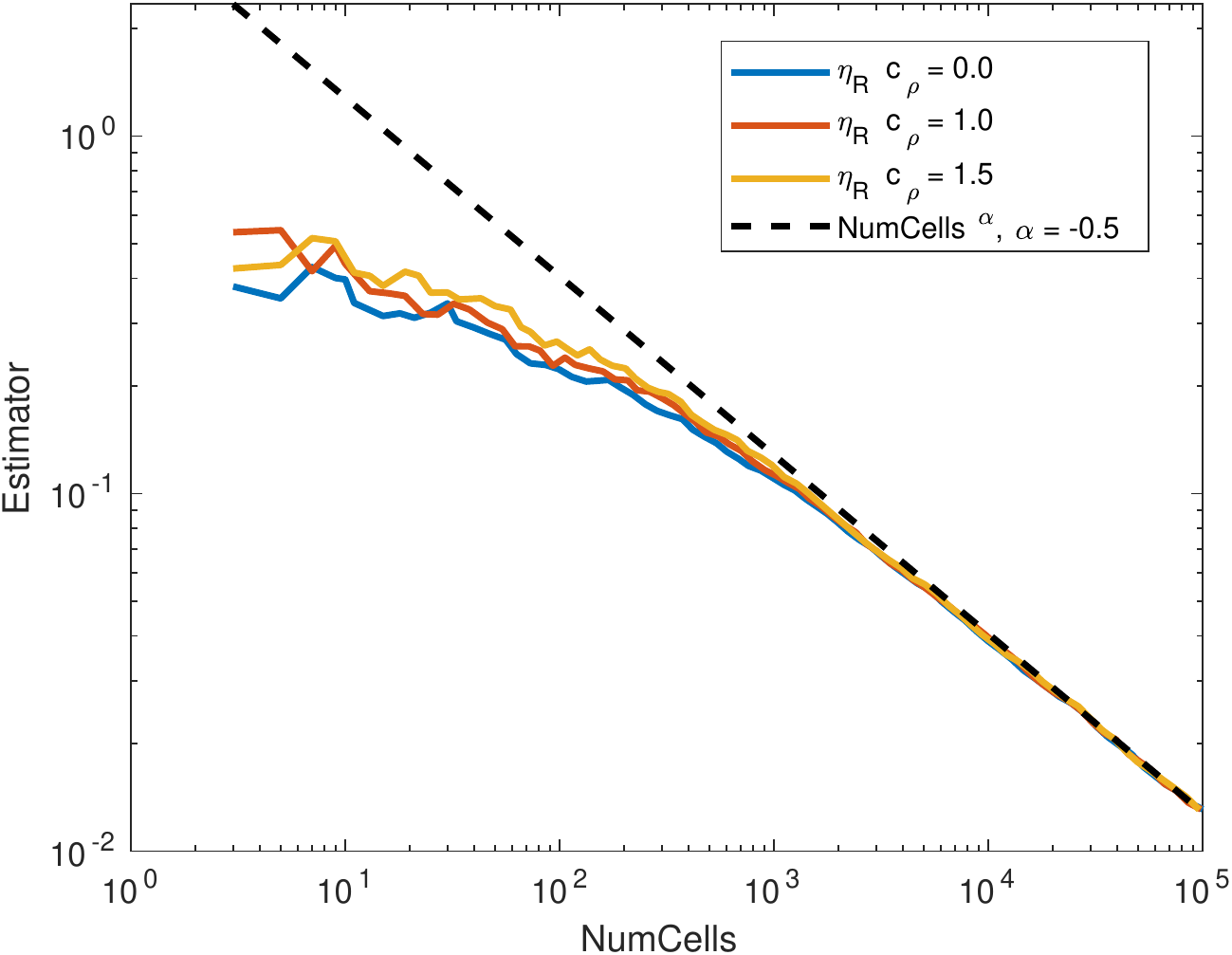}
    \caption{POLY\_LD}
    \label{fig:PLV_EstCells}
  \end{subfigure}
  \caption{Error Estimator $\eta_R$ versus $\#\Th^k-\#\Th^0$}
  \label{fig:EstCells}
\end{figure}

\begin{figure}
  \centering
  \begin{subfigure}[b]{0.42\linewidth}
    \includegraphics[width=\linewidth]{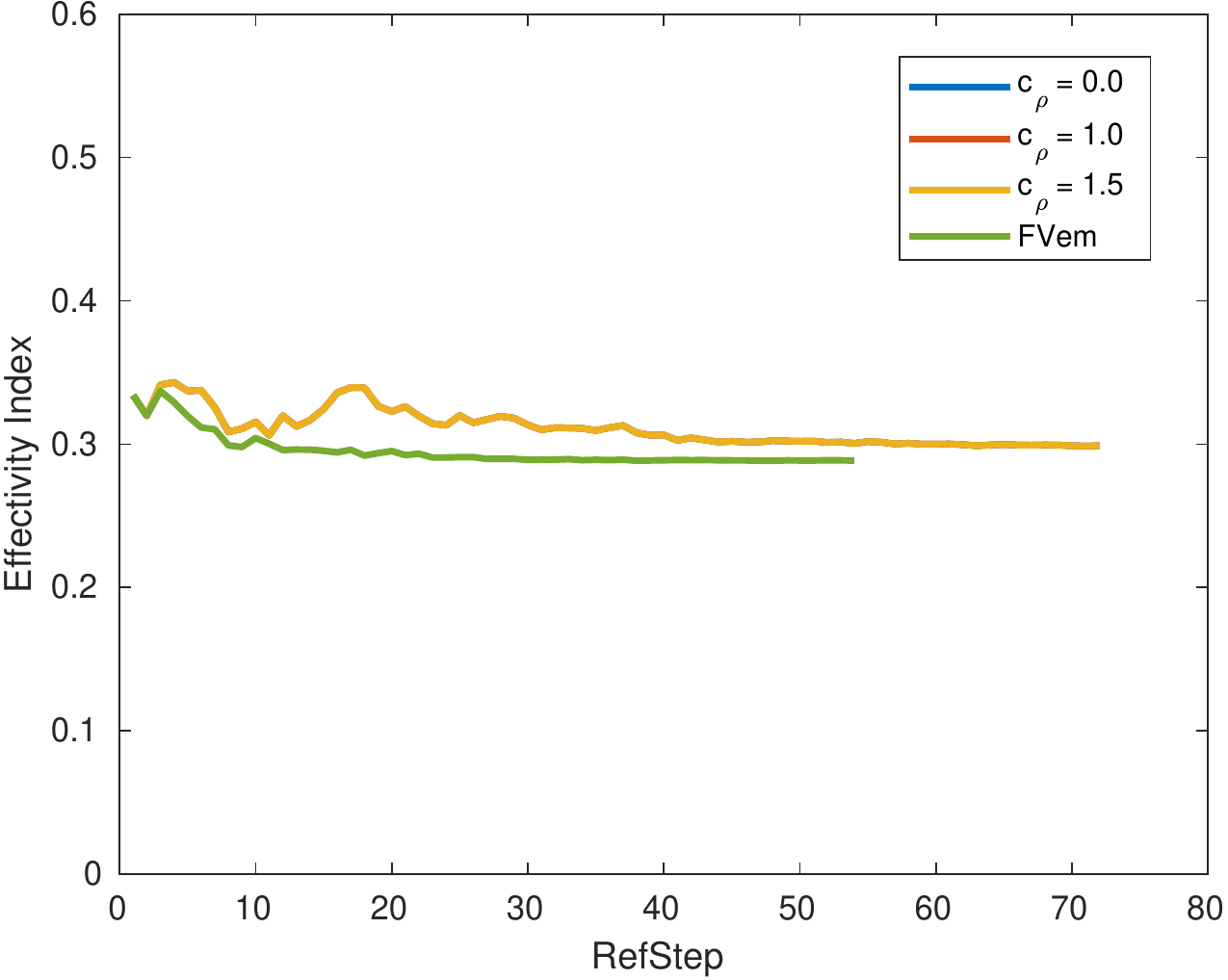}
    \caption{TRI\_MM-LE}
    \label{fig:TMVLF_EffectivityIndex}
  \end{subfigure}
  \hfill
  %
    \begin{subfigure}[b]{0.42\linewidth}
    \includegraphics[width=\linewidth]{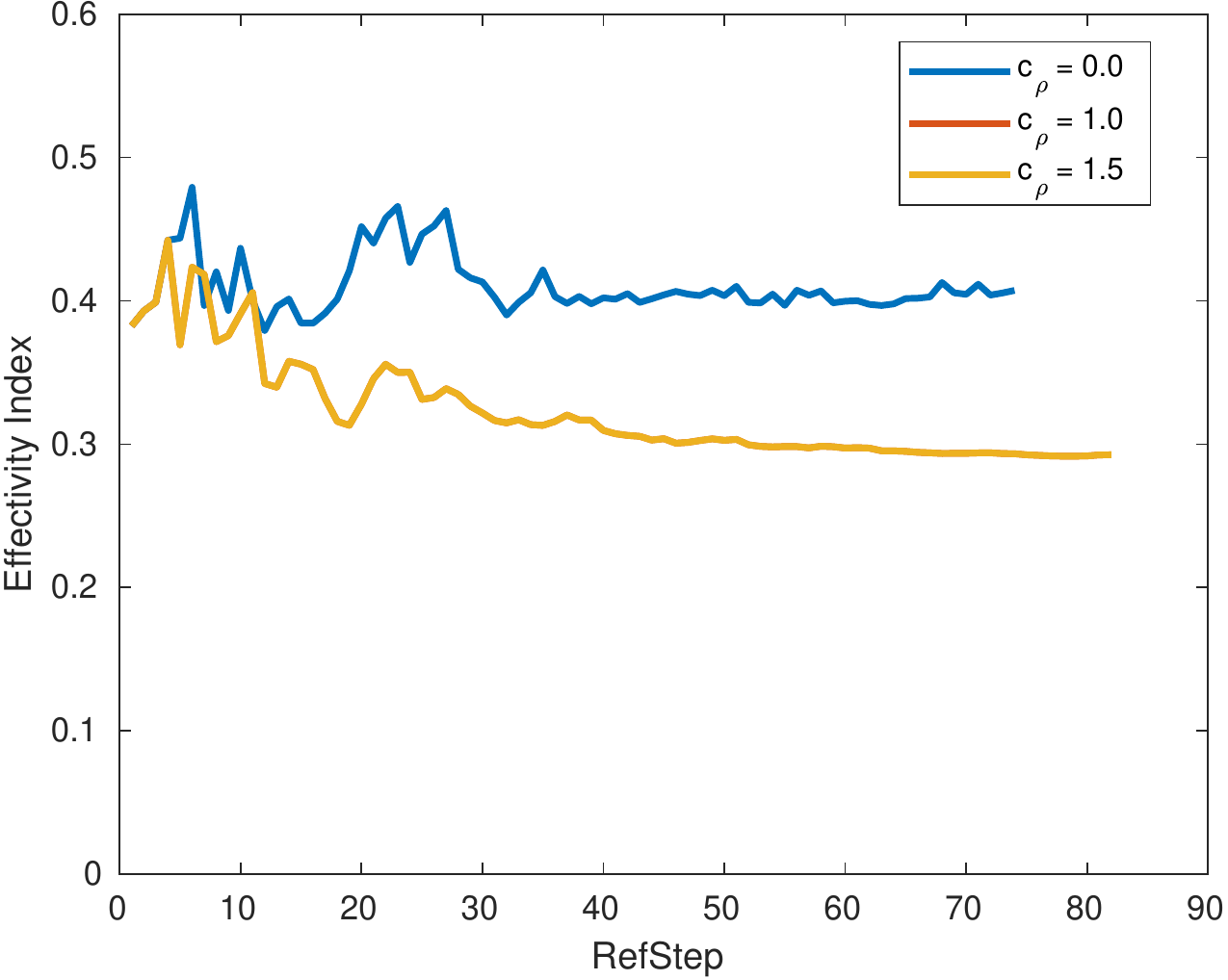}
    \caption{TRAP\_MM}
    \label{fig:TMV_EffectivityIndex}
  \end{subfigure}
  \hfill
  %
  \begin{subfigure}[b]{0.42\linewidth}
    \includegraphics[width=\linewidth]{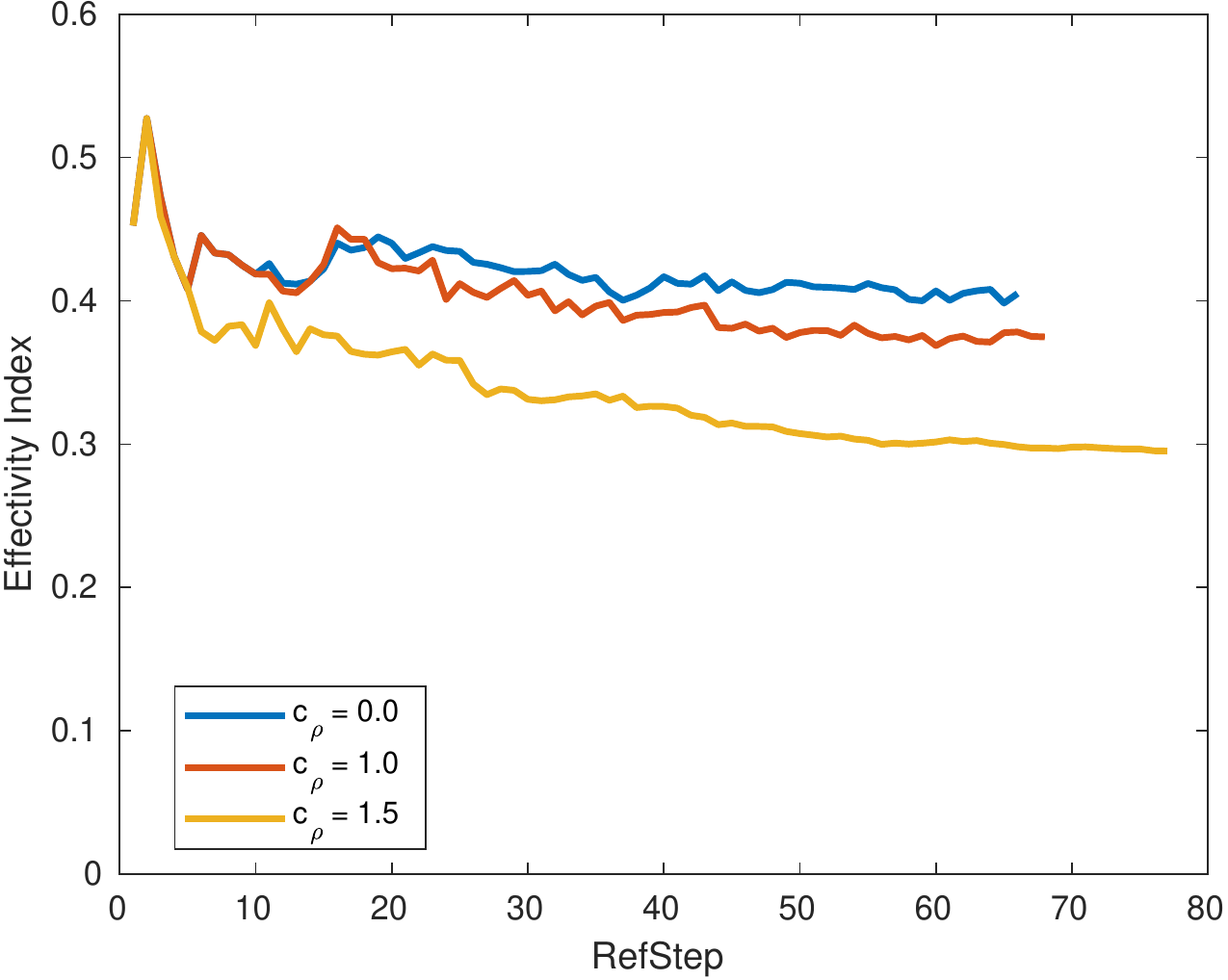}
    \caption{POLY\_MM}
    \label{fig:PMV_EffectivityIndex}
  \end{subfigure}
  \hfill
  \begin{subfigure}[b]{0.42\linewidth}
    \includegraphics[width=\linewidth]{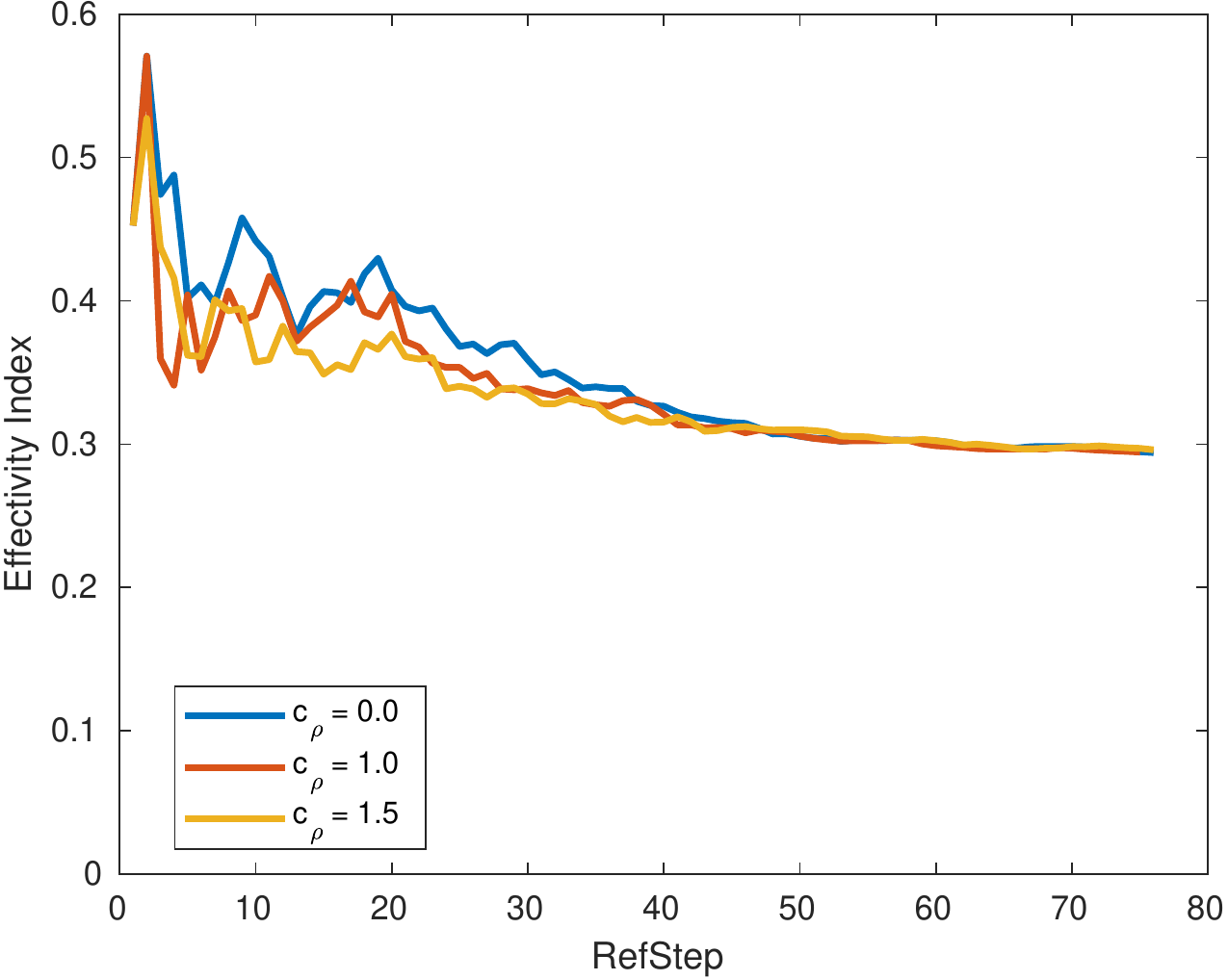}
    \caption{POLY\_LD}
    \label{fig:PLV_EffectivityIndex}
  \end{subfigure}
  \caption{Effectivity Index $e.i.$ versus refinement step}
  \label{fig:EffectivityDofs}
  \centering
  \begin{subfigure}[b]{0.42\linewidth}
    \includegraphics[width=\linewidth]{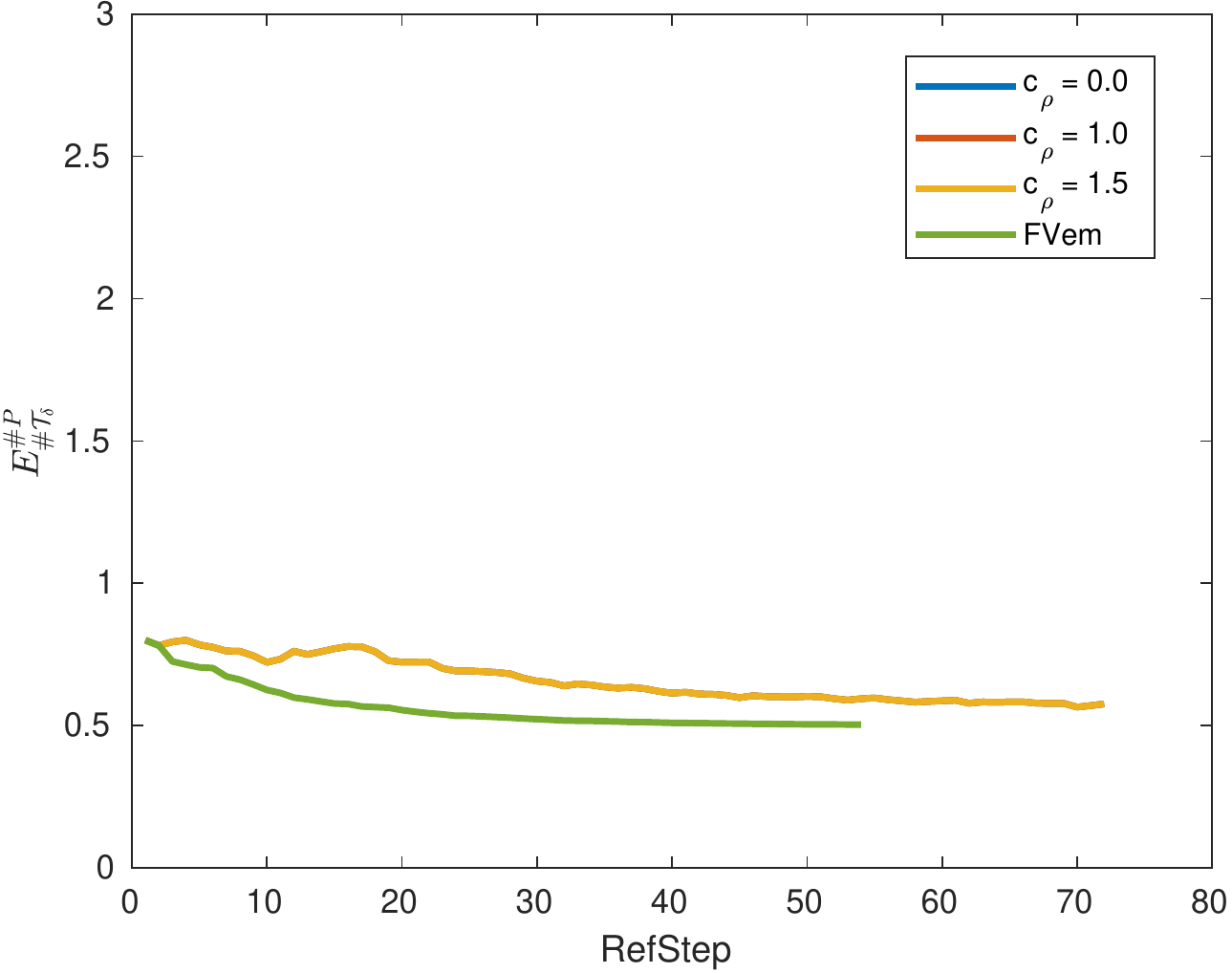}
    \caption{TRI\_MM-LE}
    \label{fig:TMVLF_NumDofsNumCells}
  \end{subfigure}
  \hfill
    \begin{subfigure}[b]{0.42\linewidth}
    \includegraphics[width=\linewidth]{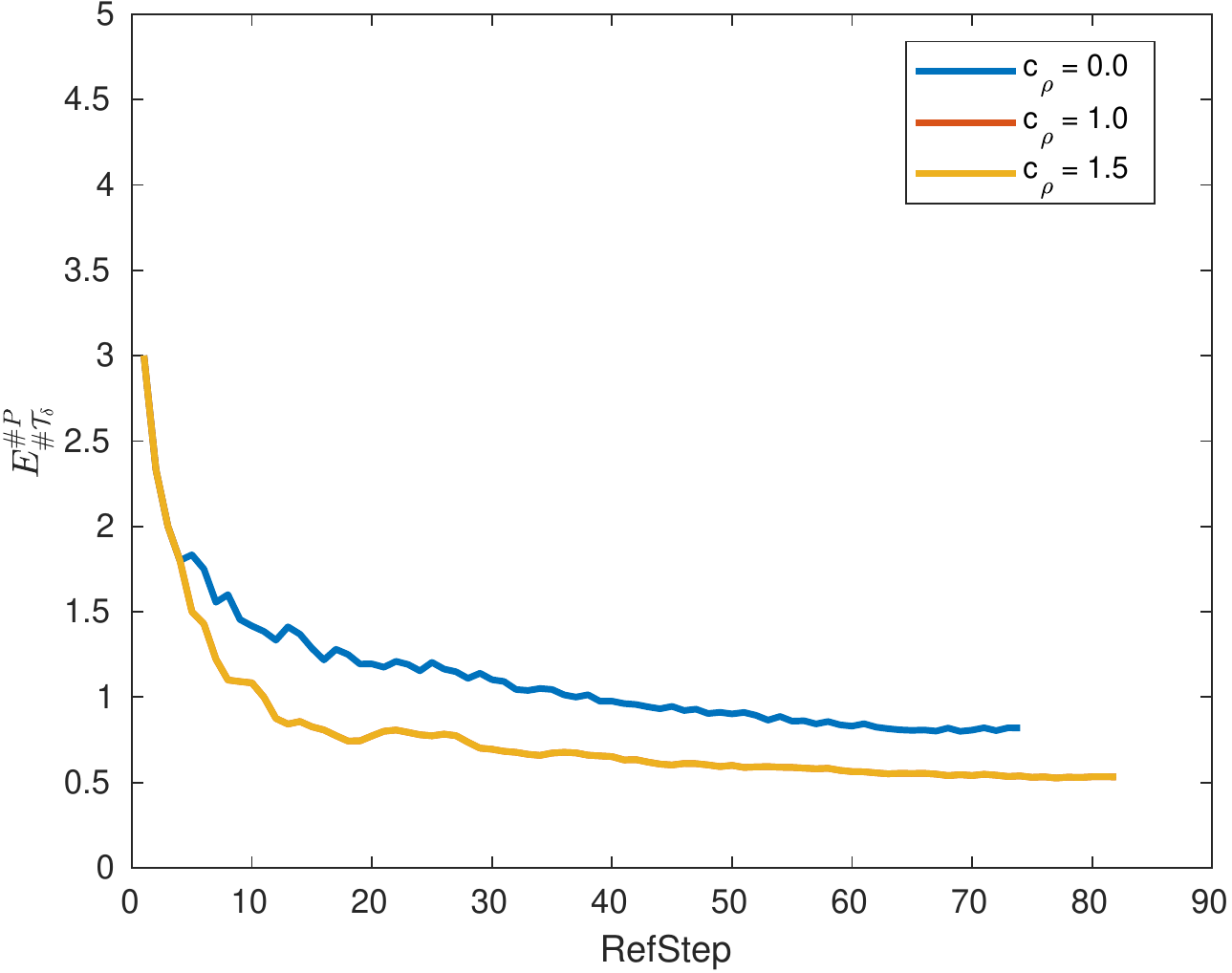}
    \caption{TRAP\_MM}
    \label{fig:TMV_NumDofsNumCells}
  \end{subfigure}
  \hfill
  \begin{subfigure}[b]{0.42\linewidth}
    \includegraphics[width=\linewidth]{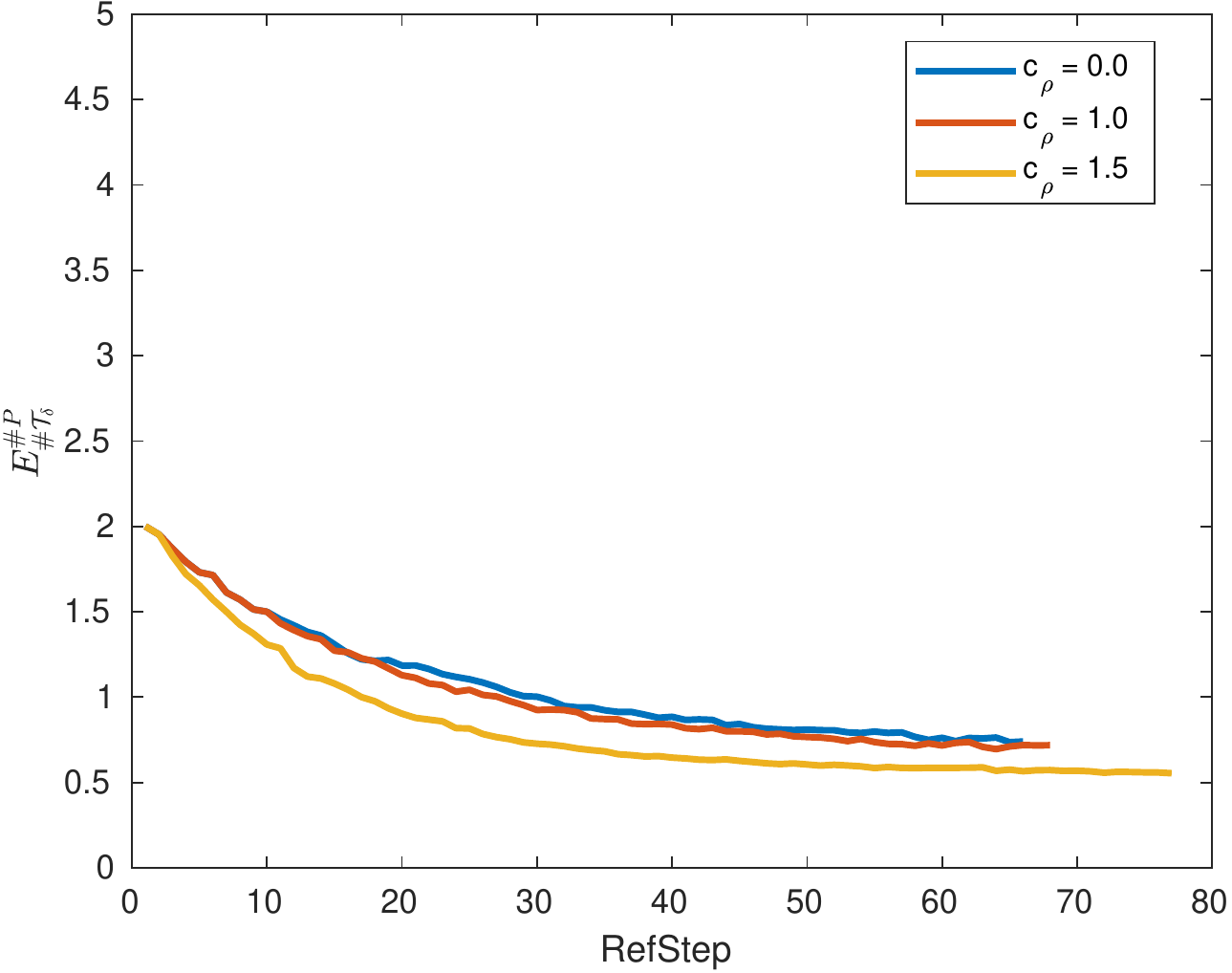}
    \caption{POLY\_MM}
    \label{fig:PMV_NumDofsNumCells}
  \end{subfigure}
  \hfill
  \begin{subfigure}[b]{0.42\linewidth}
    \includegraphics[width=\linewidth]{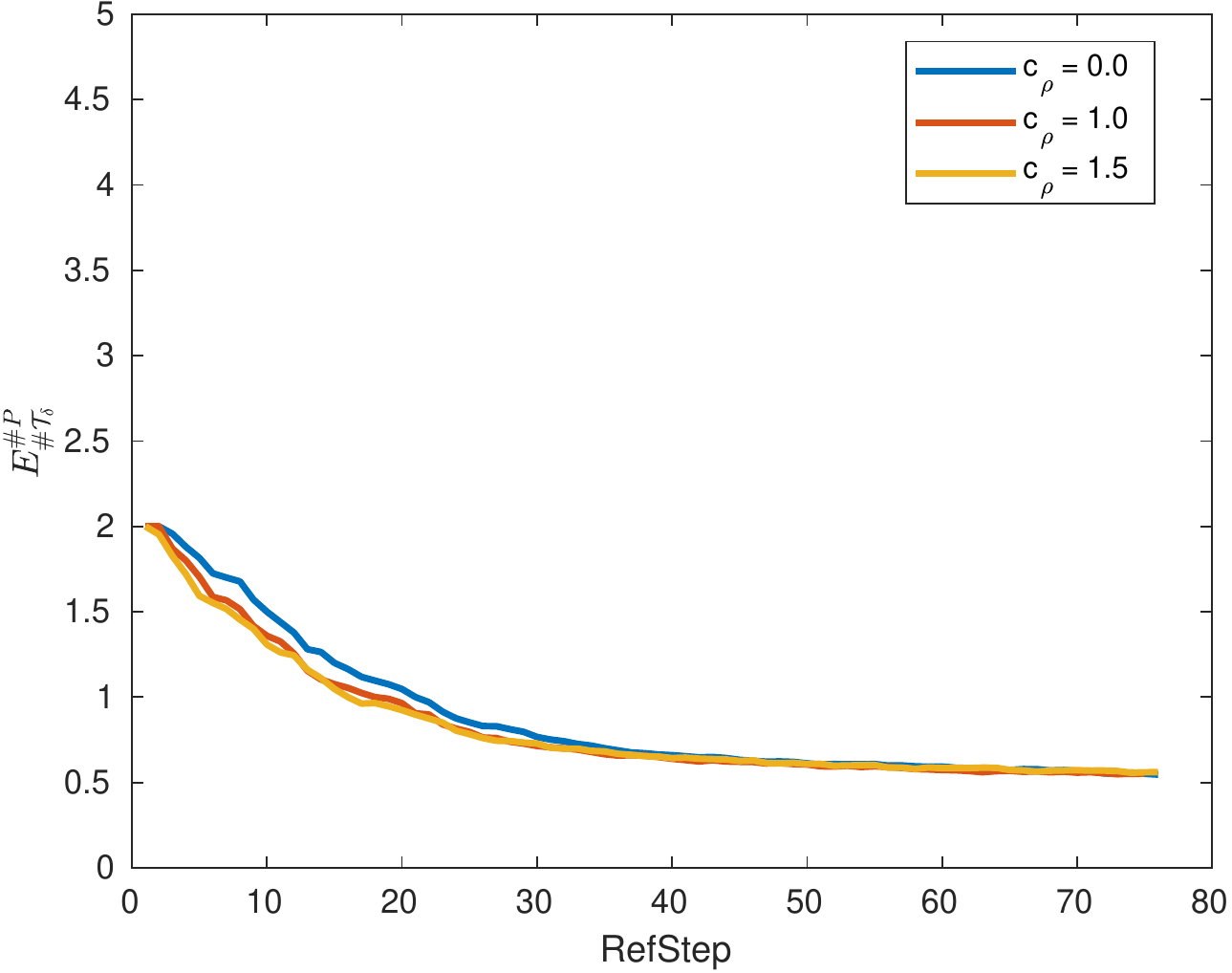}
    \caption{POLY\_LD}
    \label{fig:PLV_NumDofsNumCells}
  \end{subfigure}
  \caption{$E^{\#P}_{\#{\Th}}$}
  \label{fig:NumDofsNumCells}
\end{figure}


\section{Numerical Results: convergence and mesh quality}
\label{sec:NumRes}
In this section we analyse the performances and the overall behaviour of the refinement Algorithm \ref{alg:Refinement} applied to the solution of the Laplace problem in the $L$-shape domain $\Omega=(-1,1)^2\backslash (-1,0)^2$ (see Figure~\ref{fig:StartingMeshes}) with exact solution
\begin{align*}
 u(r,\theta) &= r^{\frac{2}{3}} \sin(\frac{2}{3}(\theta+\frac{\pi}{2})),
\end{align*}
where $r, \theta$ are the polar coordinate.
The forcing function is $f = 0$, homogeneous Dirichlet boundary conditions are imposed on the edges which generate the re-entrant corner, and  suitable non homogeneous Dirichlet boundary conditions obtained by the evaluation of the exact solution are set on the other edges. This is a quite common to test problem for a posteriori error estimates and mesh adaptive methods due to the bounded regularity of the solution \cite{doi:10.1137/1.9781611972030}.

 The adaptive algorithms is based on the standard loop:
\begin{equation}
SOLVE \rightarrow ESTIMATE \rightarrow MARK \rightarrow REFINE.
\end{equation}
The marking strategy follows from \cite{Stevenson2007OptimalityOA}, Section 5, with the parameter $\theta=0.5$ using the error estimator introduced in Section~\ref{sec:Apost} and the stopping criterion
\begin{equation}
\frac{\eta_R^2}{\norm[]{\sqrt{\K}\nabla u^\pi_\delta}^2  } \leq 1.0E^{-4},
\end{equation}
with a maximum number of iterations set to $100$.

We consider different starting meshes and analyse the behaviour of the algorithm with VEM of order one.
 The notation used to define the meshes and the refinement options is:
\begin{itemize}
\item TRI\_MM-LE
  \begin{itemize}
  \item Starting mesh: triangular mesh, Figure~\ref{fig:StartingMeshes} (left).
  \item Refinement strategy: Maximum Moment (Algorithm~\ref{alg:MaximumMomentum}) with quality check ($c_\rho$=0, 1, 1.5).
  \item Refinement strategy (FVem): Longest Edge with conformity recovery for comparison with FEM P1.
  \end{itemize}
\item TRAP\_MM
  \begin{itemize}
  \item Starting mesh: trapezoidal mesh, Figure~\ref{fig:StartingMeshes} (center).
  \item Refinement strategy: Maximum Moment (Algorithm~\ref{alg:MaximumMomentum}) with quality check ($c_\rho$=0, 1, 1.5).
  \end{itemize}
\item POLY\_MM
  \begin{itemize}
  \item Starting mesh: polygonal mesh, Figure~\ref{fig:StartingMeshes} (right).
  \item Refinement strategy: Maximum Moment (Algorithm~\ref{alg:MaximumMomentum}) with quality check ($c_\rho$=0, 1, 1.5).
  \end{itemize}
\item POLY\_LD
  \begin{itemize}
  \item Starting mesh: polygonal, Figure~\ref{fig:StartingMeshes} (right).
  \item Refinement strategy: Longest Diagonal (Algorithm~\ref{alg:LongestDiagonal}) with quality check ($c_\rho$=0, 1, 1.5).
  \end{itemize}
\end{itemize}

The triangular mesh is generated by the Triangle library \cite{shewchuk96b,SHEWCHUK200221}, the polygonal one is generated by PolyMesher \cite{polymesher}.
For the Maximum Moment strategy (Algorithm~\ref{alg:MaximumMomentum}) and for the Longest Diagonal strategy (Algorithm~\ref{alg:LongestDiagonal}) we consider three different values of the coefficient $c_\rho = {0.0,1.0,1.5}$ (Line~\ref{rho-eB} in Algorithm~\ref{alg:SmoothingDirection}).
The case with $c_\rho=0$ is the case in which we always cut the intersected edge in the mid-point, or, in case of the intersection with aligned edges, we cut the set of aligned edges in the midpoint of the set.

In all the cases we analyse the rate of convergence of the {\em a posteriori} error estimator $\eta_R$, the effectivity index $e.i.$ and the ratio between the number of mesh nodes and the number of cells in the mesh $E^{\#P}_{\#{\Th}}$. We also analyse the quality of the mesh by observing the parameters $AR_E^{Rr}$,  $AR^{edge}_E$, $AR_E^{Hr}$, $AR_E^{Hh}$,  $K_E^{\Pi^\nabla}$. Finally, we report the fraction of cells that are triangle (with exactly three vertices, i.e. without aligned edges) $R^\triangle_{\Th}$, and the fraction of cells that are quadrilateral polygons (with exactly four vertices) $R^\diamond_{\Th}$ .

With the test TRI\_MM-LE we compare adaptive results obtained starting from a triangular mesh using the maximum moment refinement with the results obtained by the longest edge refinement with triangular conformity recovery (case denoted by FVem in the figures), i.e. with the results obtained on a suitable mesh by conforming FEM elements P1.

In Figure~\ref{fig:EstDofs} we report the convergence behaviour of the estimator $\eta_R$ with respect to the \#dofs. We note that the behaviour is essentially the same on the different meshes and using Maximum Moment and Longest Diagonal strategies. We remark that the reference rate of convergence $\alpha=-0.5$ is the largest possible with elements of order one. 
We also remark that all the collapsing strategies ($c_\rho=0.0, 1.0, 1.5$) behaves similarly and in the case TRI\_MM-LE  behaves essentially like the FEM of order one on all conforming triangular meshes.
We remark that in Figure \ref{fig:TMV_EstDof} the case $c_\rho=1.0$ is overlapped to the case $c_\rho=1.5$ (the same will happen also in the next figures concerning POLY\_MM), whereas in Figure \ref{fig:PMV_EstDof} it is very close to $c_\rho=0.0$. 

In Figure~\ref{fig:EstCells} we report the behaviour of the estimator with respect to the added cells in the refinement process, that is also the number of cells marked for refinement.
Figure~\ref{fig:EffectivityDofs} displays the effectivity index with respect to the refinement steps. We can observe that its behaviour is quite stable,
in the case TRI\_MM-LE all the collapsing strategies produce overlapped curves, very close to the FVem case, in the other plots the different collapsing strategies may slightly impact on the values. 

\begin{figure}
  \centering
    \begin{subfigure}[b]{0.32\linewidth}
    \includegraphics[width=\linewidth]{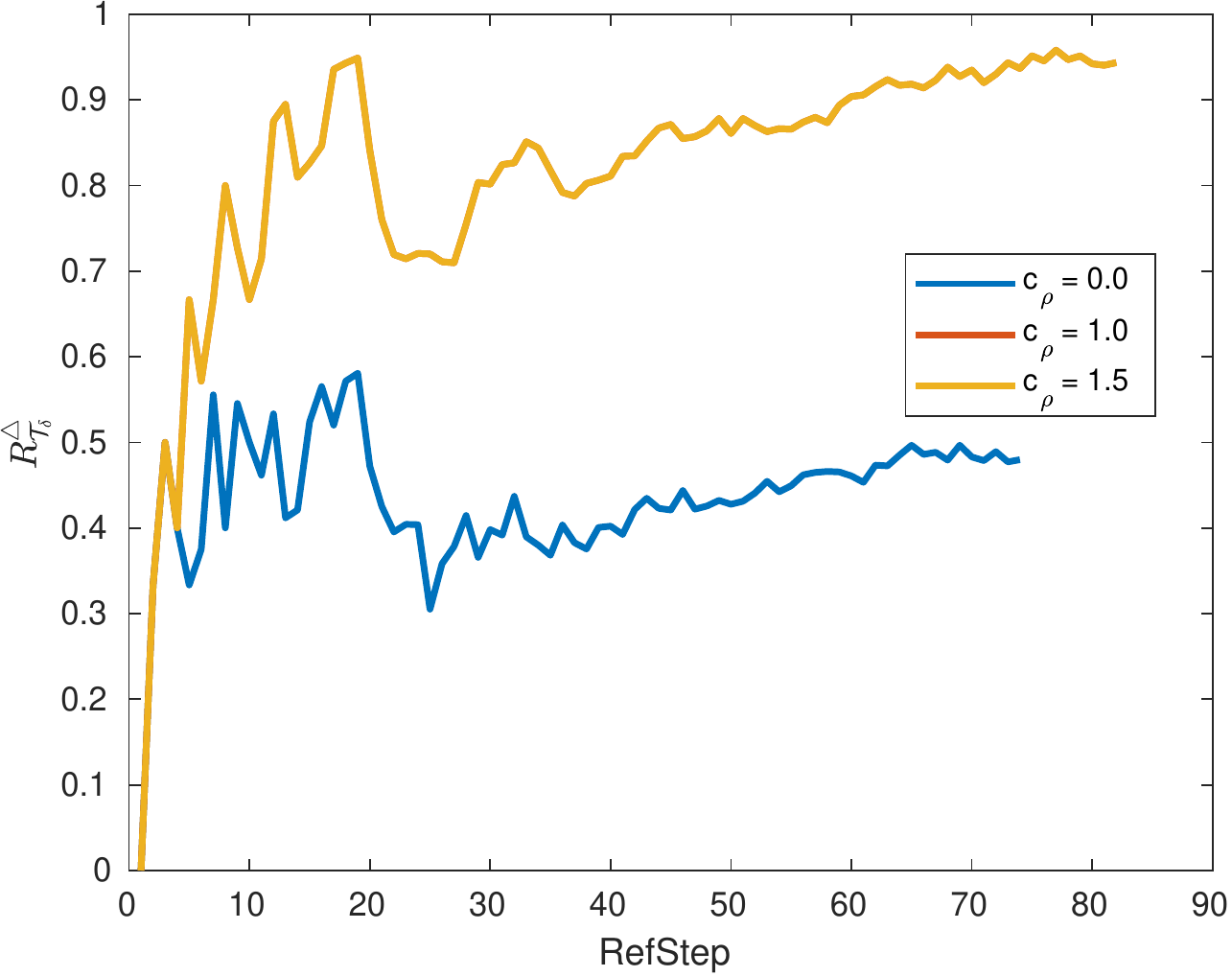}
    \caption{TRAP\_MM}
    \label{fig:TMV_NumTriangles}
  \end{subfigure}
  \hfill
  \begin{subfigure}[b]{0.32\linewidth}
    \includegraphics[width=\linewidth]{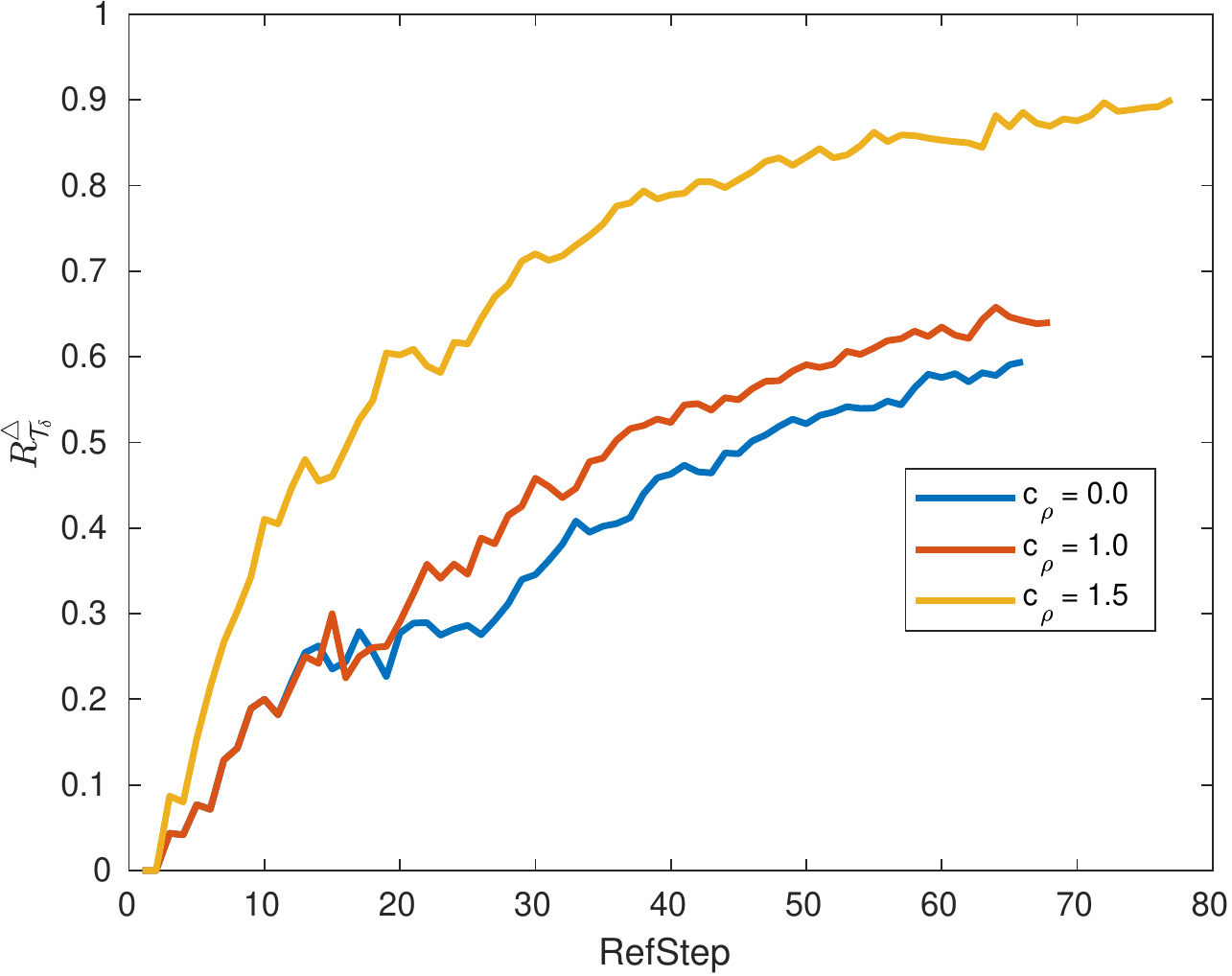}
    \caption{POLY\_MM}
    \label{fig:PMV_NumTriangles}
  \end{subfigure}
  \hfill
  \begin{subfigure}[b]{0.32\linewidth}
    \includegraphics[width=\linewidth]{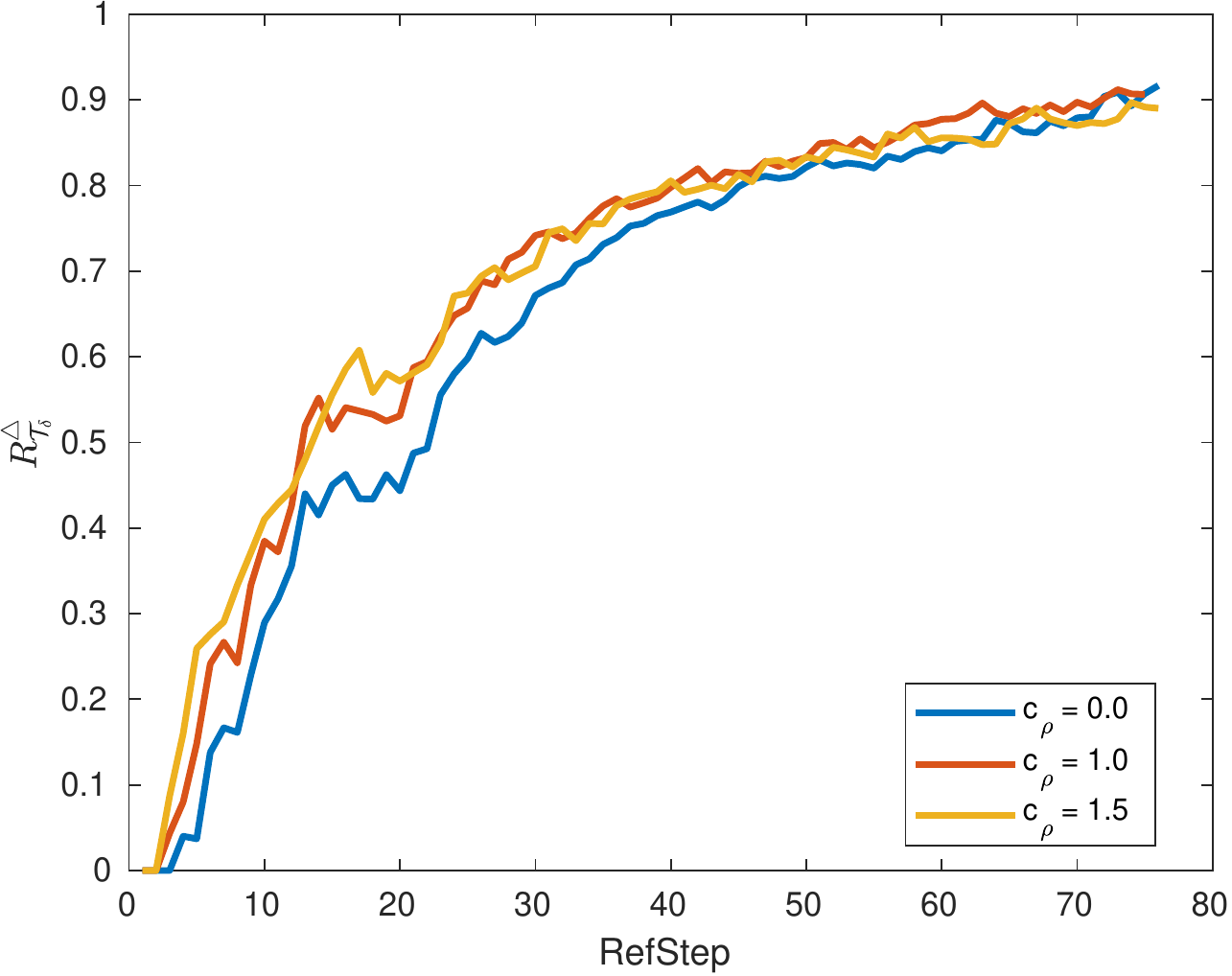}
    \caption{POLY\_LD}
    \label{fig:PLV_NumTriangles}
  \end{subfigure}
  \caption{$R^\triangle_{\Th}$}
  \label{fig:NumTri}
  \centering
    \begin{subfigure}[b]{0.32\linewidth}
    \includegraphics[width=\linewidth]{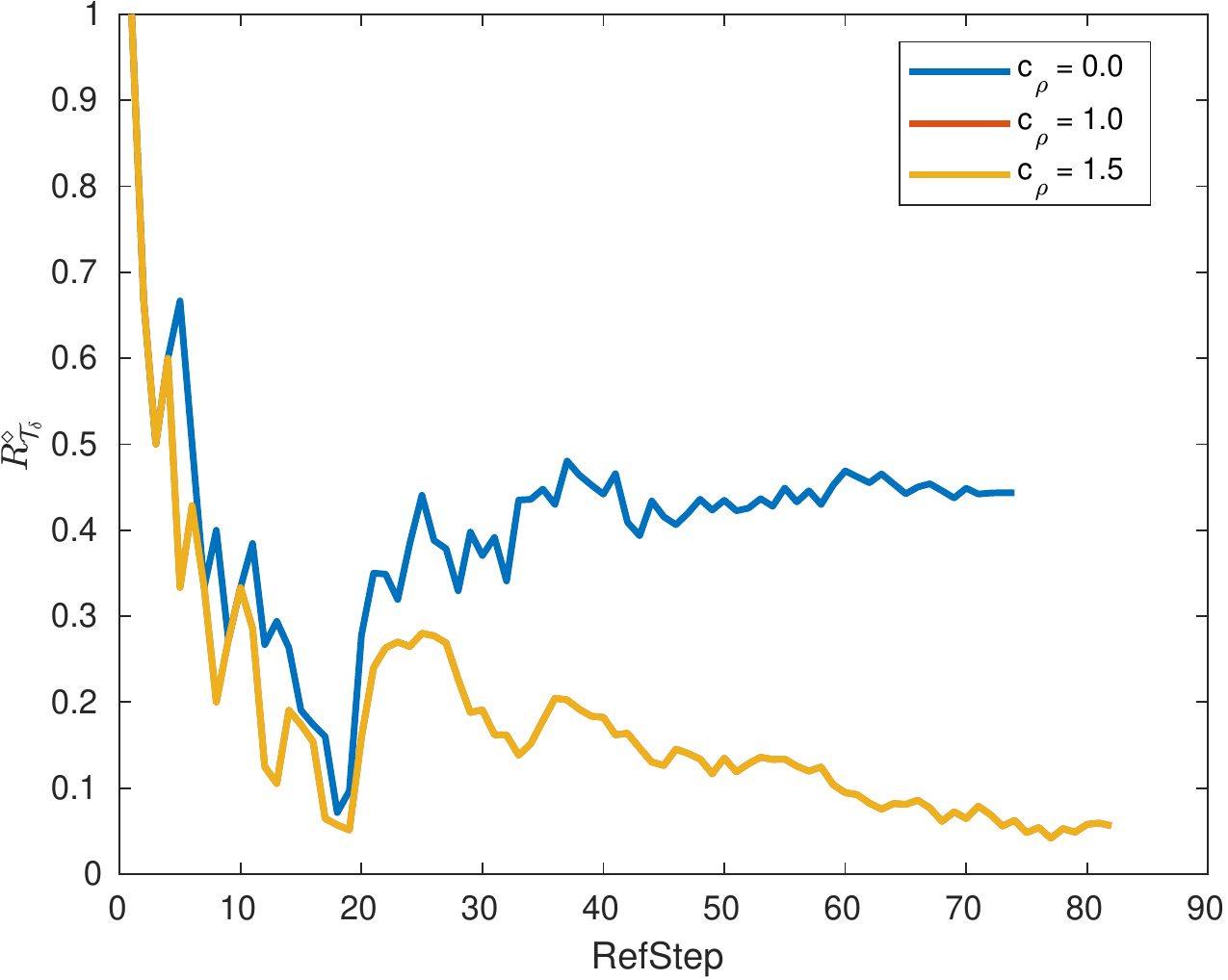}
    \caption{TRAP\_MM}
    \label{fig:TMV_NumSquares}
  \end{subfigure}
  \hfill
  \begin{subfigure}[b]{0.32\linewidth}
    \includegraphics[width=\linewidth]{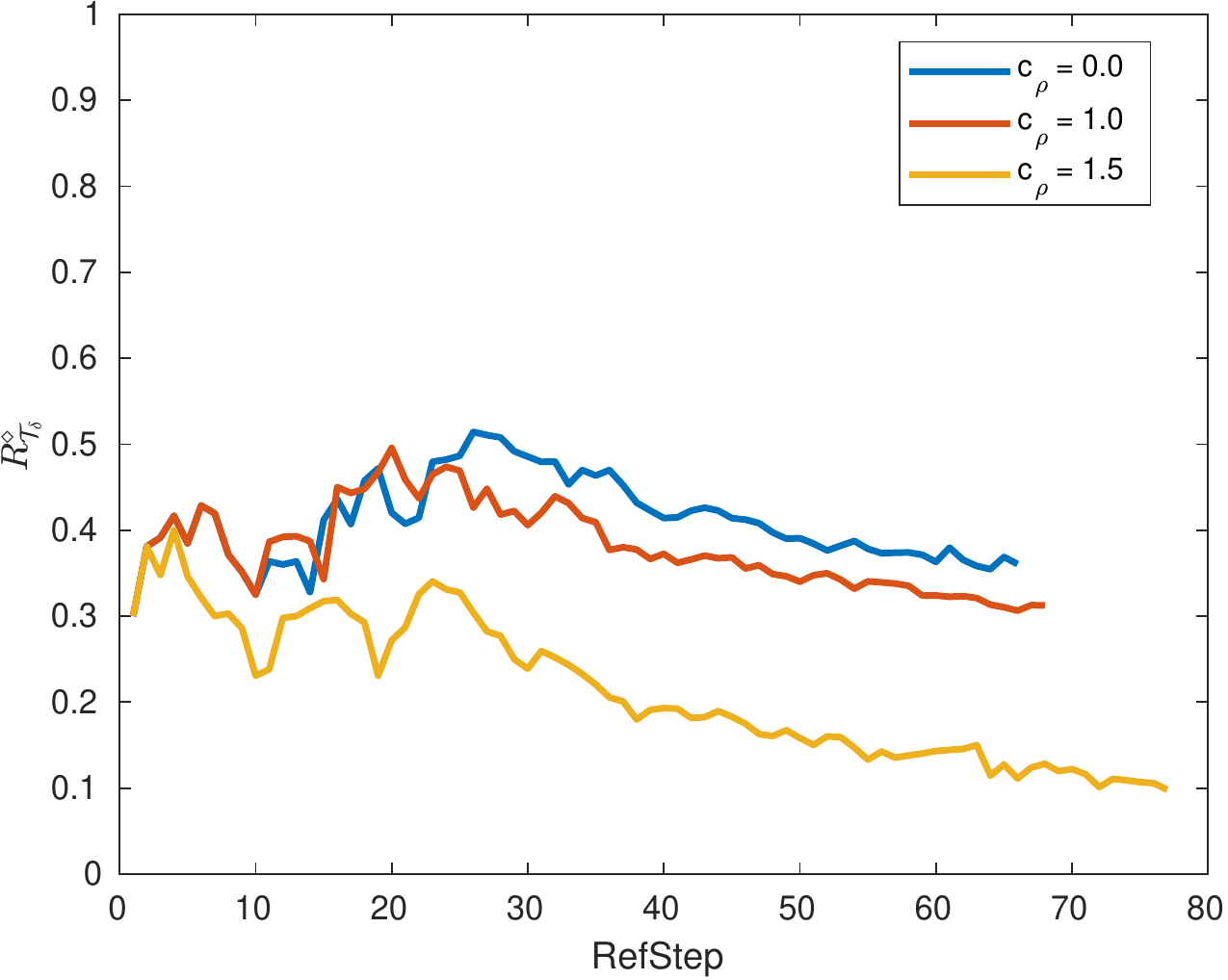}
    \caption{POLY\_MM}
    \label{fig:PMV_NumSquares}
  \end{subfigure}
  \hfill
  \begin{subfigure}[b]{0.32\linewidth}
    \includegraphics[width=\linewidth]{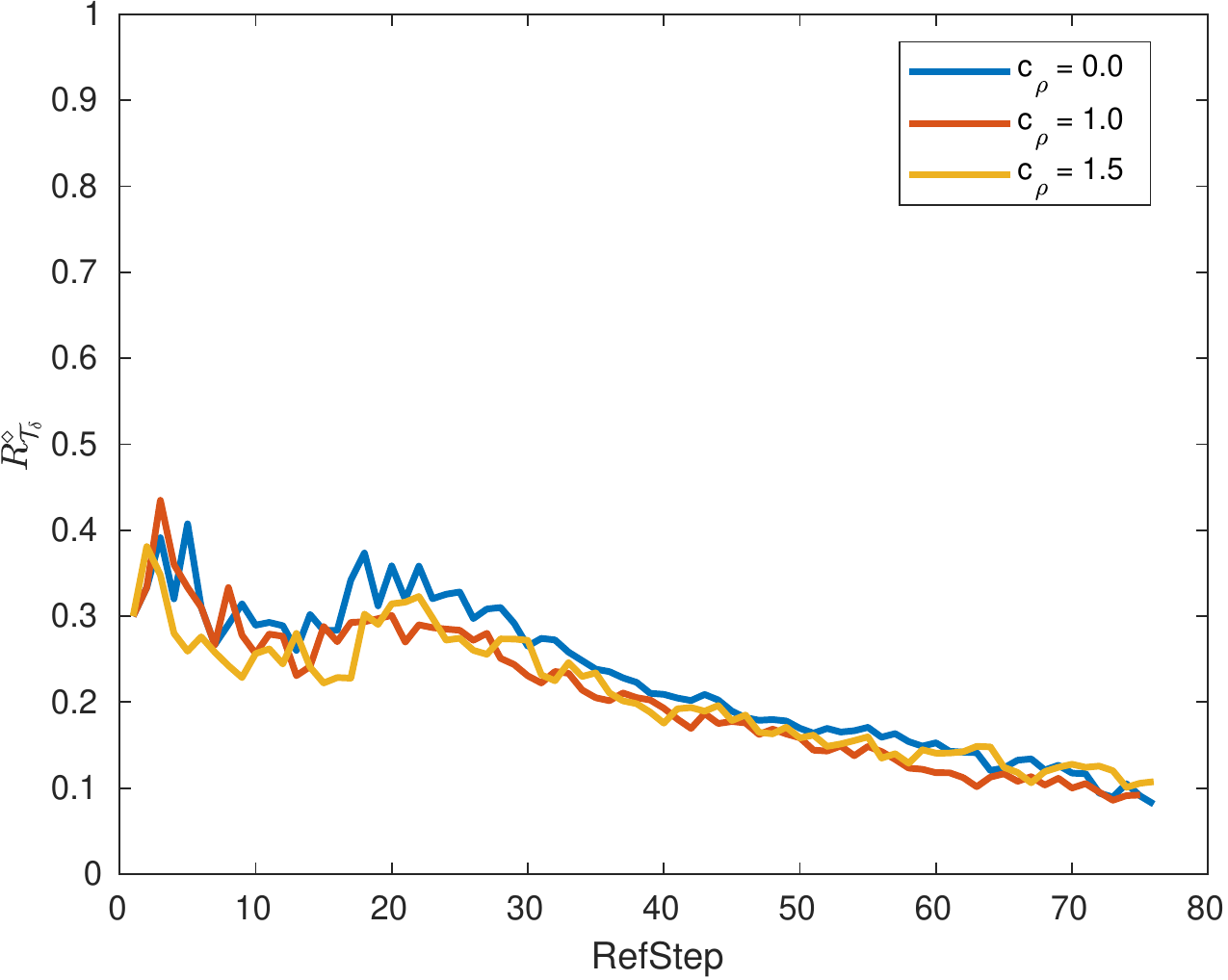}
    \caption{POLY\_LD}
    \label{fig:PLV_NumSquares}
  \end{subfigure}
  \caption{$R^\diamond_{\Th}$}
  \label{fig:NumQuad}
  \end{figure}

In Figure~\ref{fig:NumDofsNumCells} the ratio $E^{\#P}_{\#{\Th}}$ between the number of mesh points and cells is reported, we can observe that this ratio is bounded and decreasing. The meshes that during refinement progressively tend to triangular meshes (see Figure~\ref{fig:NumTri}) have this ratio approaching the value 0.5, whereas the TRAP-mesh with $c_\rho=0$ that converge to a mesh of quadrilateral elements (see Figure~\ref{fig:NumQuad}) has a limit for this ratio that approaches 1. The boundedness of this ratio is relevant for efficiency reasons being the number of cells related to the approximation capacity of the VEM solution, whereas the number of points being related to the dimension of the linear systems to be solved. The degrowth of this ratio implies that refinement does improve approximation efficiency of the new cells.


\begin{figure}
  \centering
  \begin{subfigure}[b]{0.49\linewidth}
    \includegraphics[width=\linewidth]{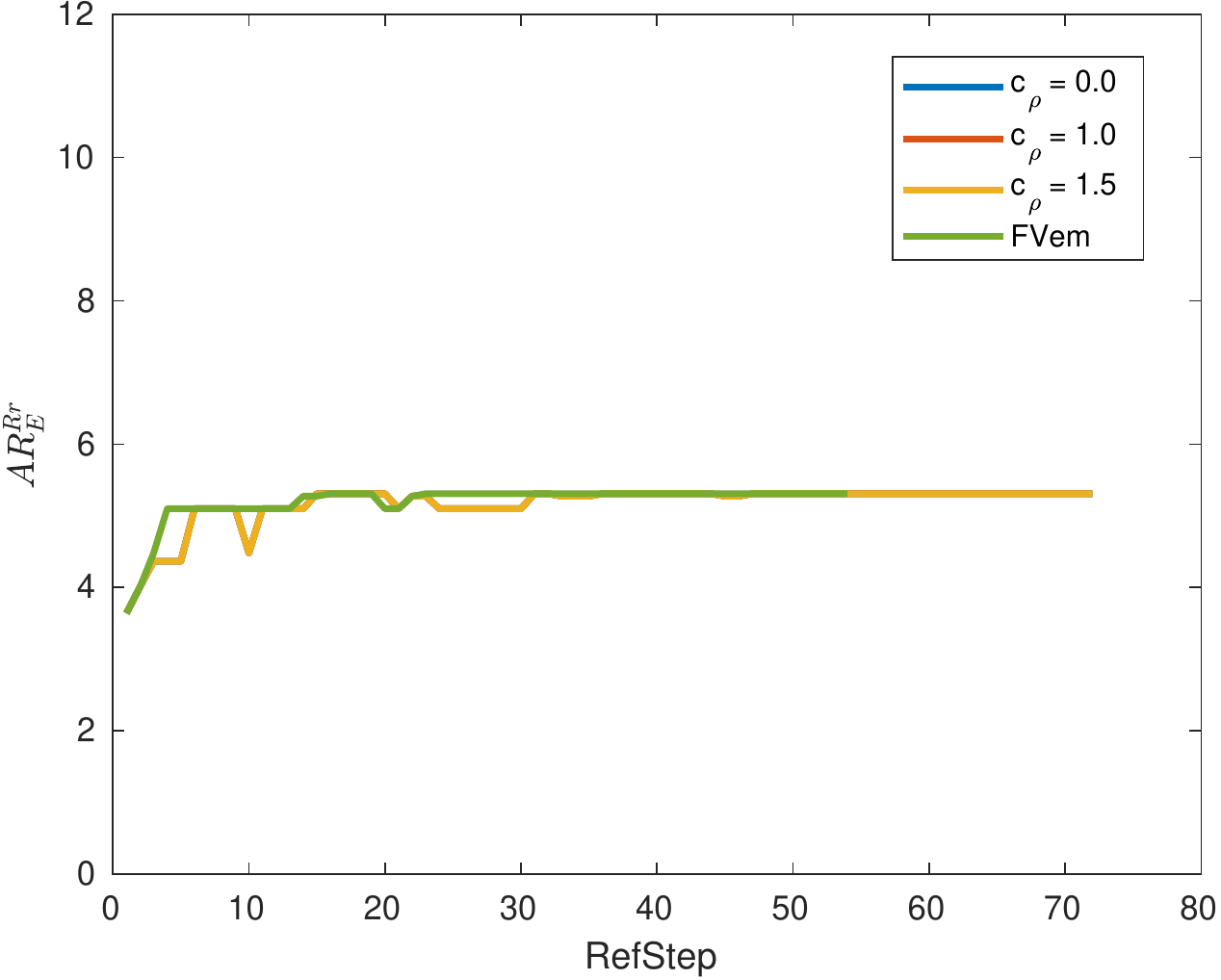}
    \caption{TRI\_MM-LE}
    \label{fig:TMVLF_MaxARR}
  \end{subfigure}
  \hfill
    \begin{subfigure}[b]{0.49\linewidth}
    \includegraphics[width=\linewidth]{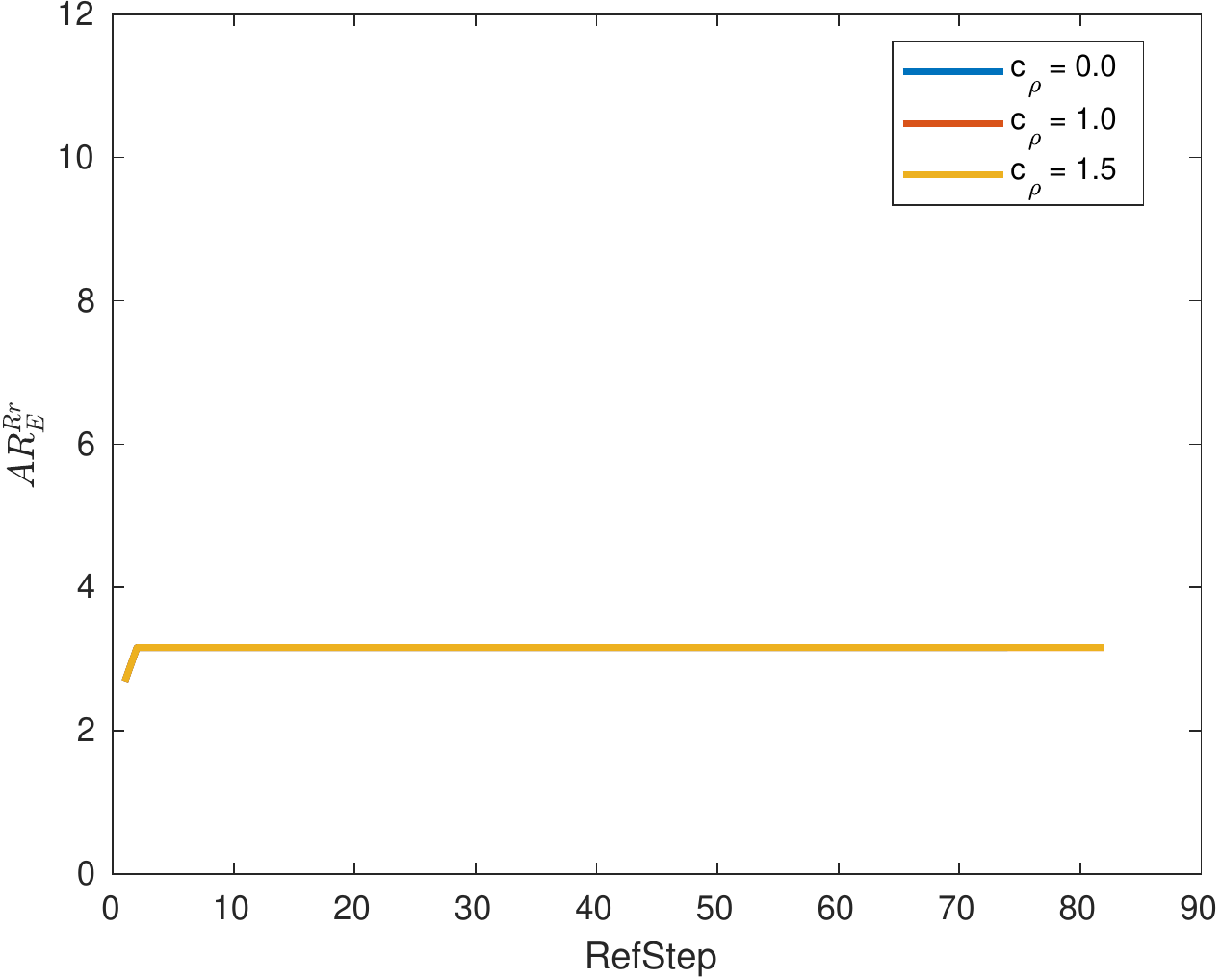}
    \caption{TRAP\_MM}
    \label{fig:TMV_MaxARR}
  \end{subfigure}
  \hfill
  \begin{subfigure}[b]{0.49\linewidth}
    \includegraphics[width=\linewidth]{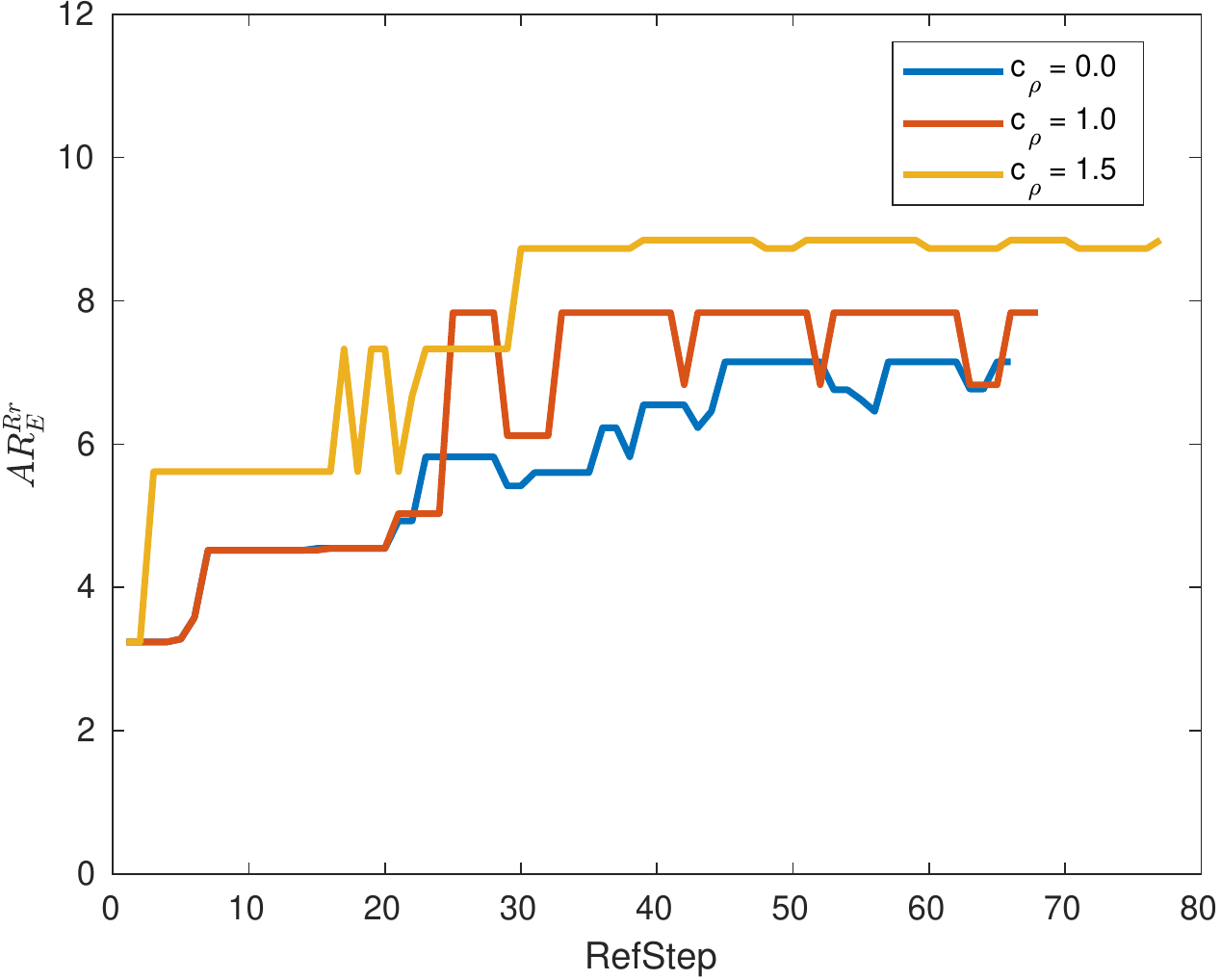}
    \caption{POLY\_MM}
    \label{fig:PMV_MaxARR}
  \end{subfigure}
  \hfill
  \begin{subfigure}[b]{0.49\linewidth}
    \includegraphics[width=\linewidth]{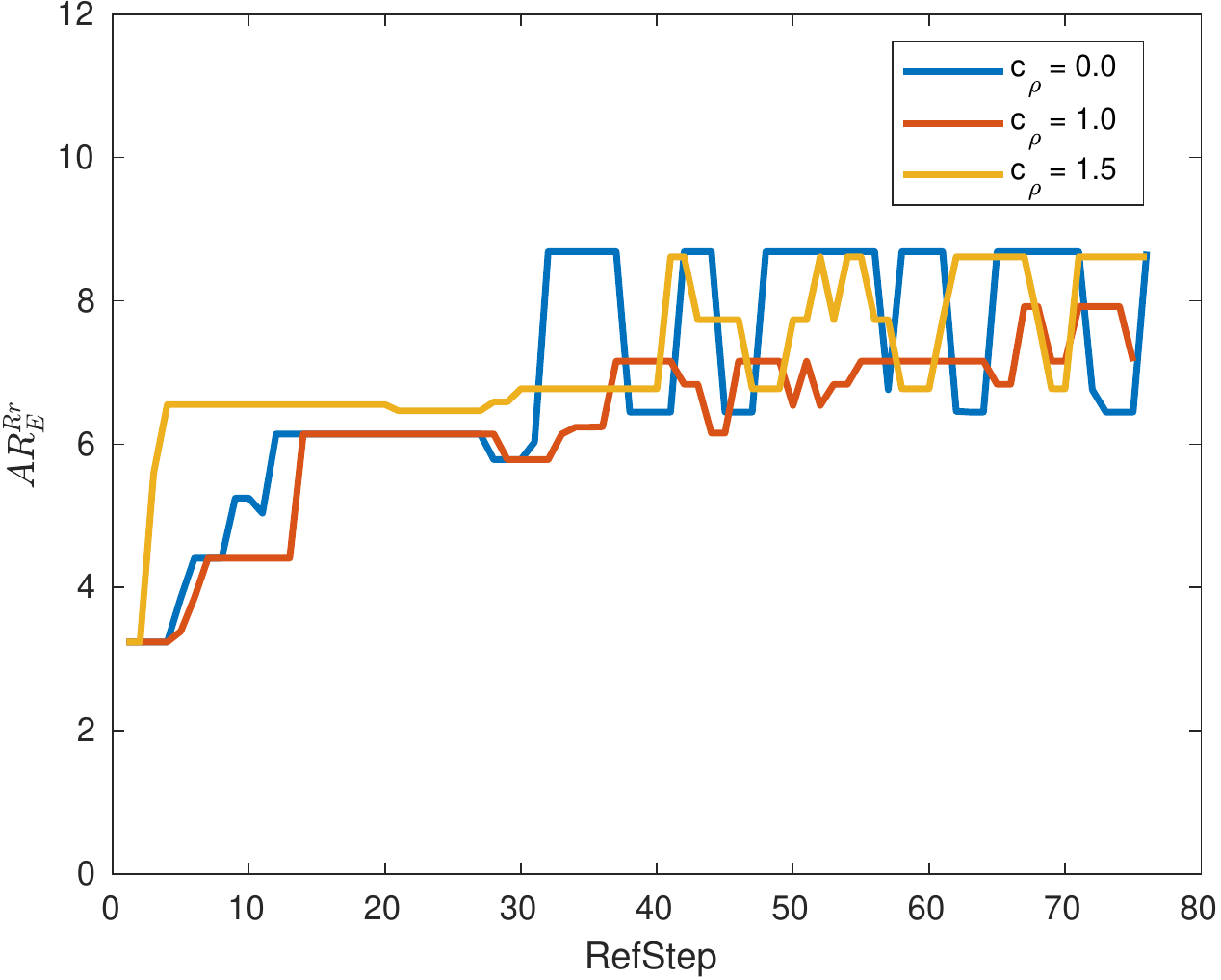}
    \caption{POLY\_LD}
    \label{fig:PLV_MaxARR}
  \end{subfigure}
  \hfill
  \begin{subfigure}[b]{0.49\linewidth}
    \includegraphics[width=\linewidth]{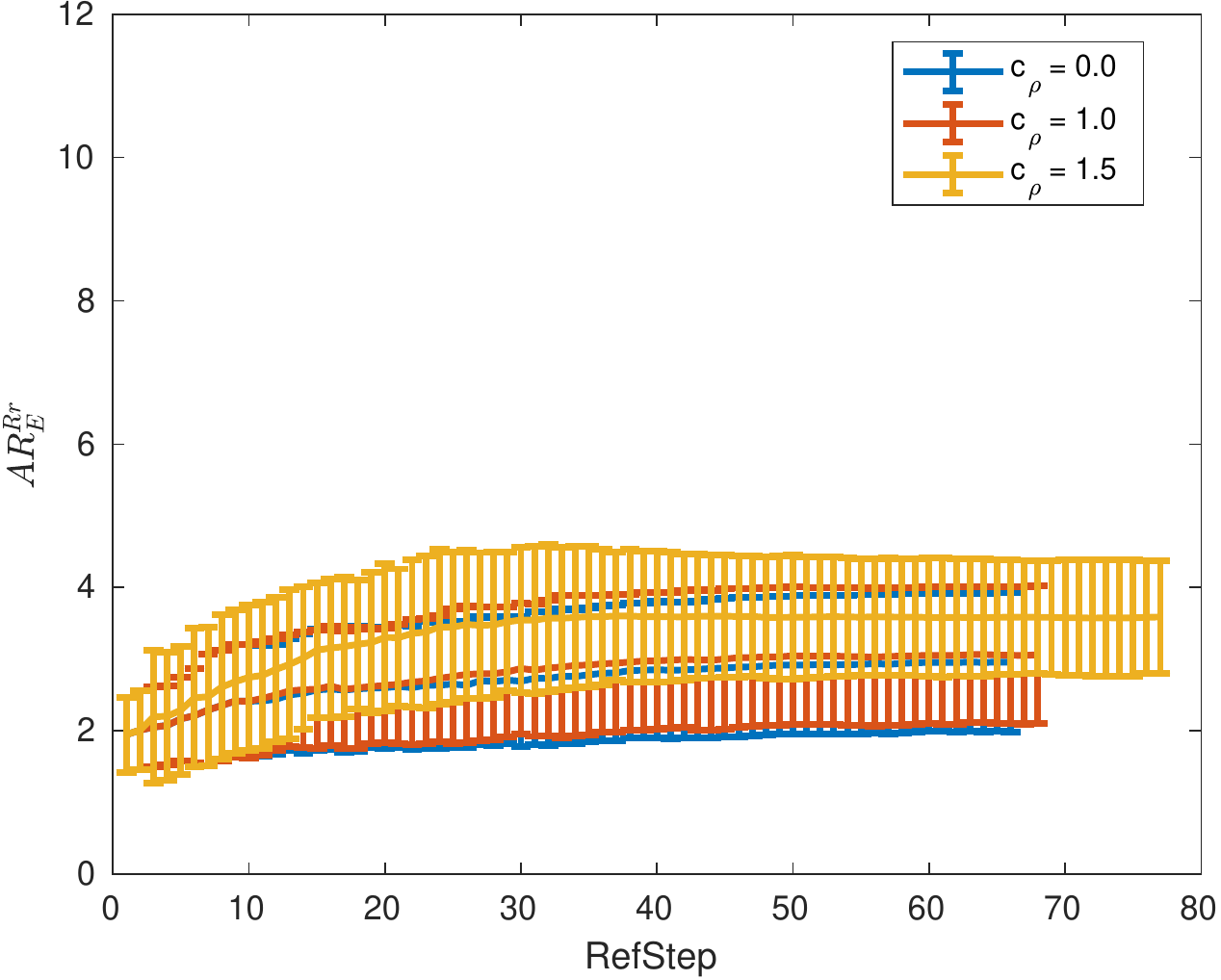}
    \caption{POLY\_MM}
    \label{fig:PMV_MeanARR}
  \end{subfigure}
  \hfill
  \begin{subfigure}[b]{0.49\linewidth}
    \includegraphics[width=\linewidth]{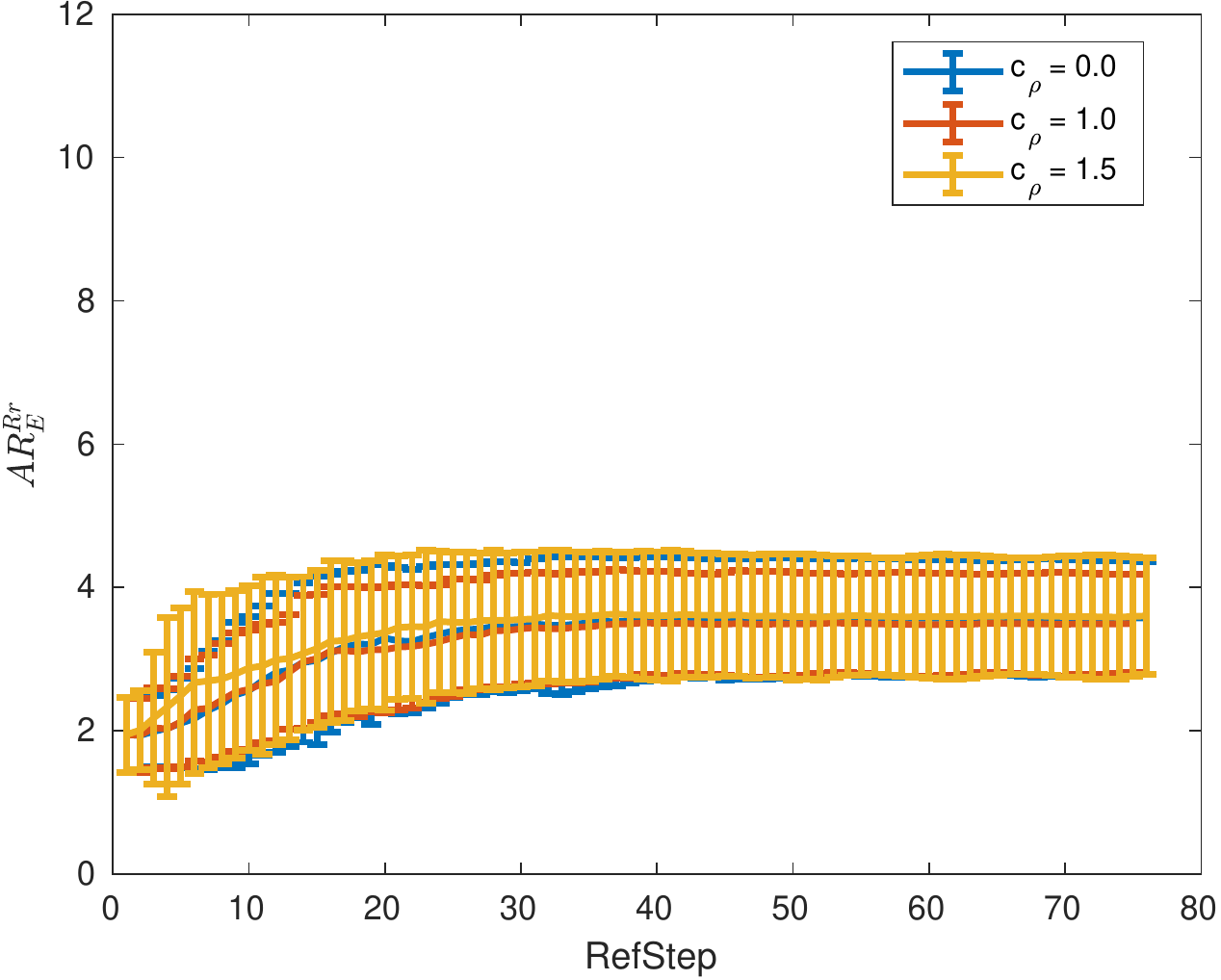}
    \caption{POLY\_LD}
    \label{fig:PLV_MeanARR}
  \end{subfigure}
  \caption{Max $AR^{Rr}_E$: (a), (b), (c), (d), Mean $AR^{Rr}_E \pm$ standard deviation: (e), (f)}
  \label{fig:MaxARR}
\end{figure}

\begin{figure}
  \centering
  \begin{subfigure}[b]{0.49\linewidth}
    \includegraphics[width=\linewidth]{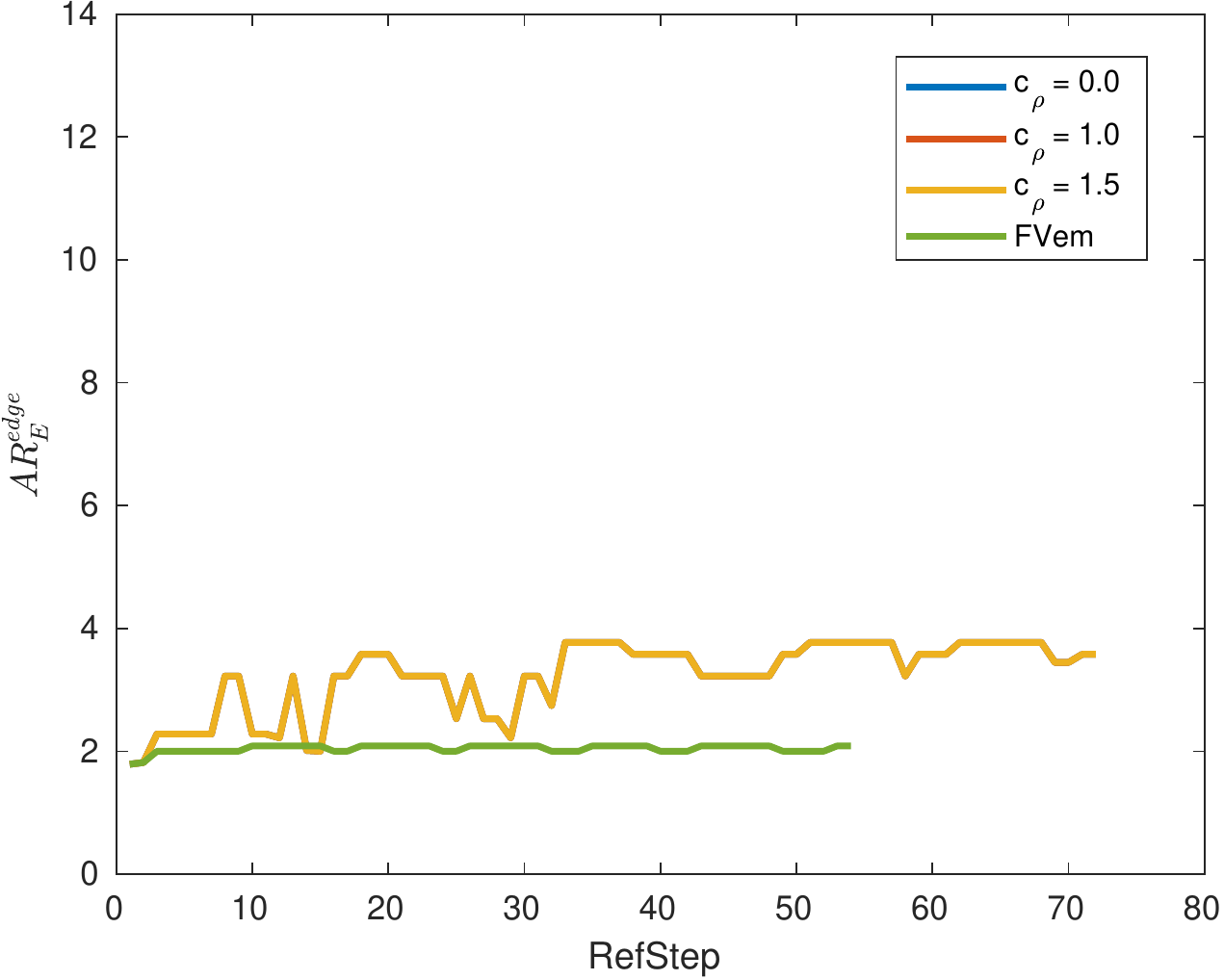}
    \caption{TRI\_MM-LE}
    \label{fig:TMVLF_MaxARH}
  \end{subfigure}
  \hfill
    \begin{subfigure}[b]{0.49\linewidth}
    \includegraphics[width=\linewidth]{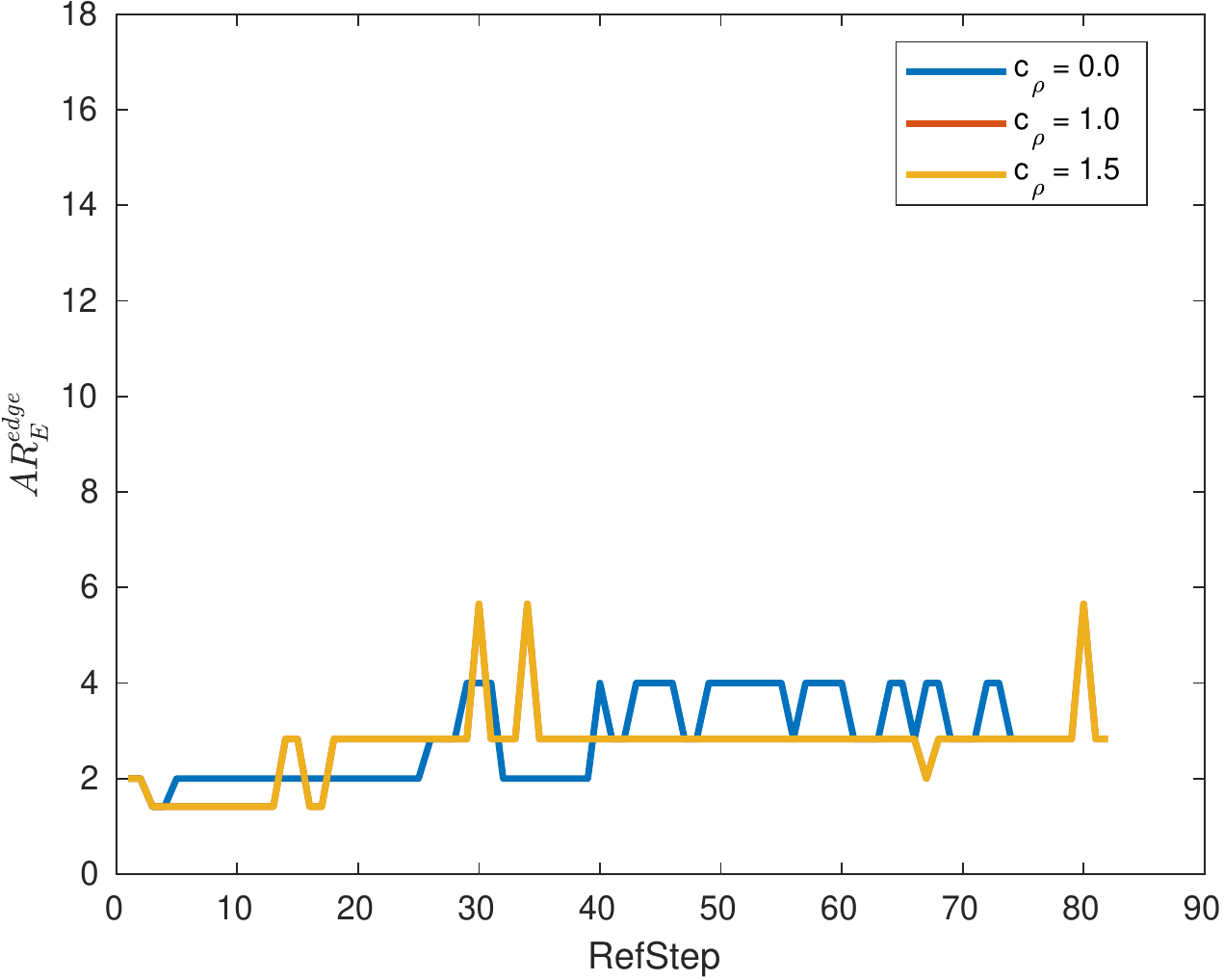}
    \caption{TRAP\_MM}
    \label{fig:TMV_MaxARH}
  \end{subfigure}
  \hfill
  \begin{subfigure}[b]{0.49\linewidth}
    \includegraphics[width=\linewidth]{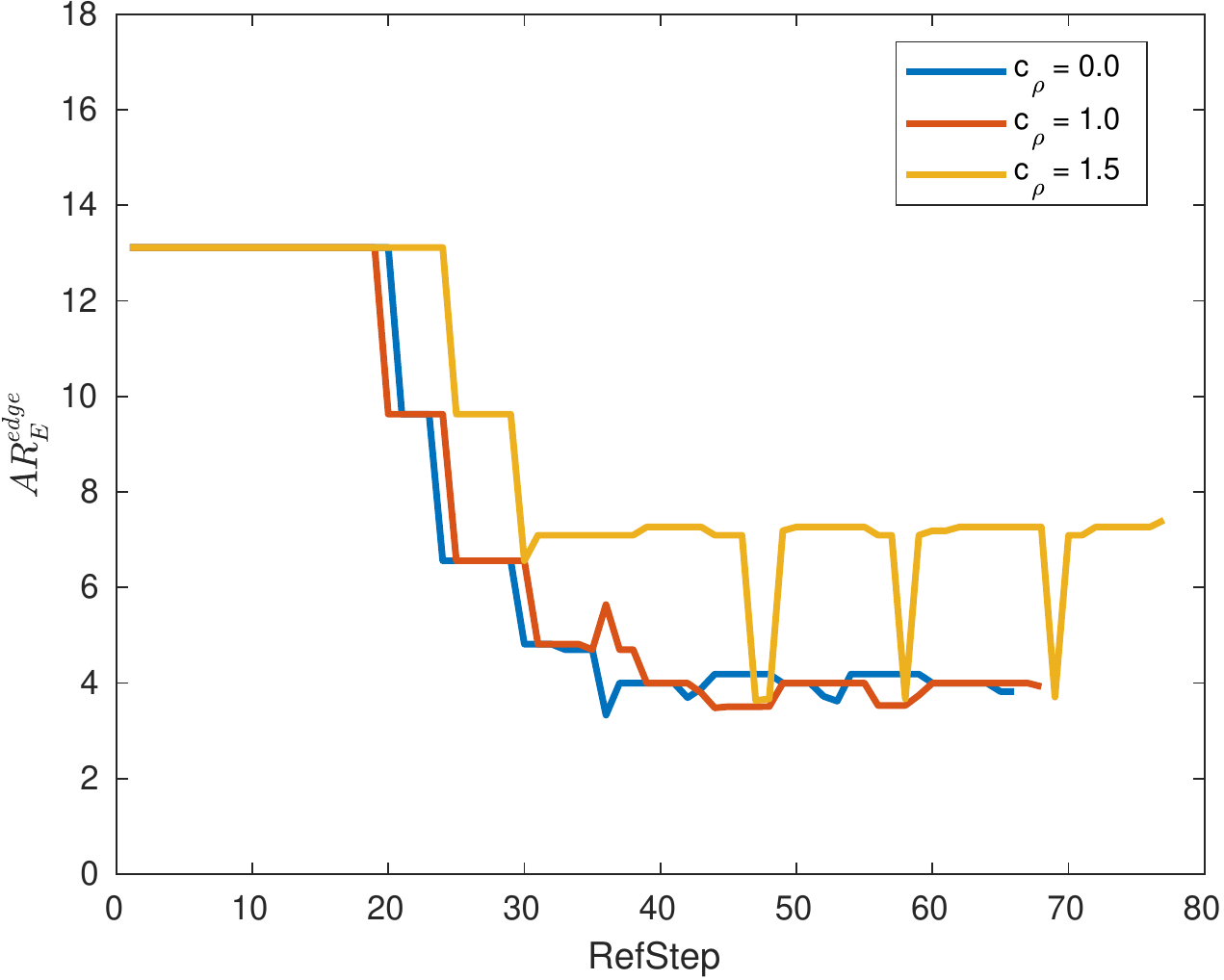}
    \caption{POLY\_MM}
    \label{fig:PMV_MaxARH}
  \end{subfigure}
  \hfill
  \begin{subfigure}[b]{0.49\linewidth}
    \includegraphics[width=\linewidth]{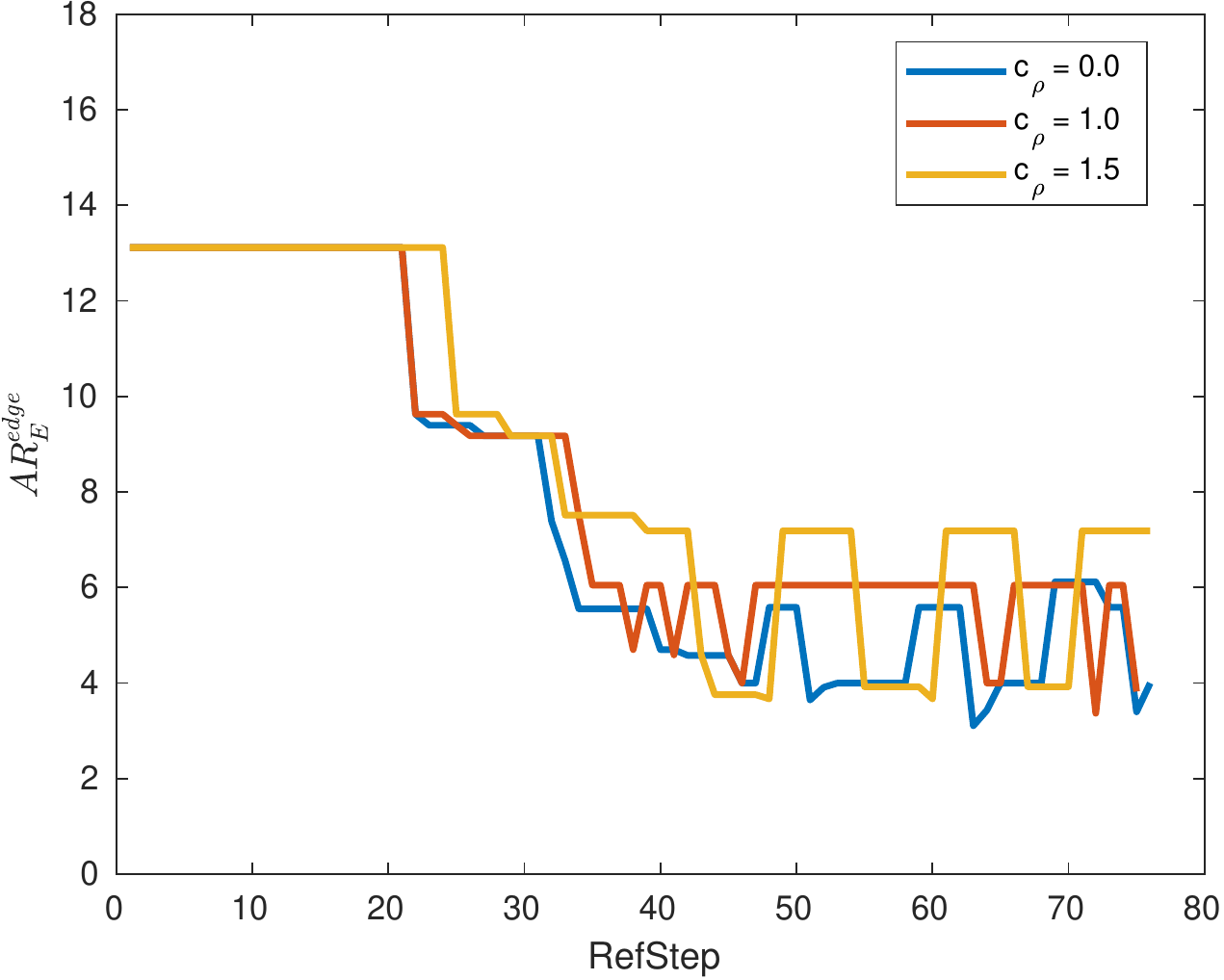}
    \caption{POLY\_LD}
    \label{fig:PLV_MaxARH}
  \end{subfigure}
  \hfill
  \begin{subfigure}[b]{0.49\linewidth}
    \includegraphics[width=\linewidth]{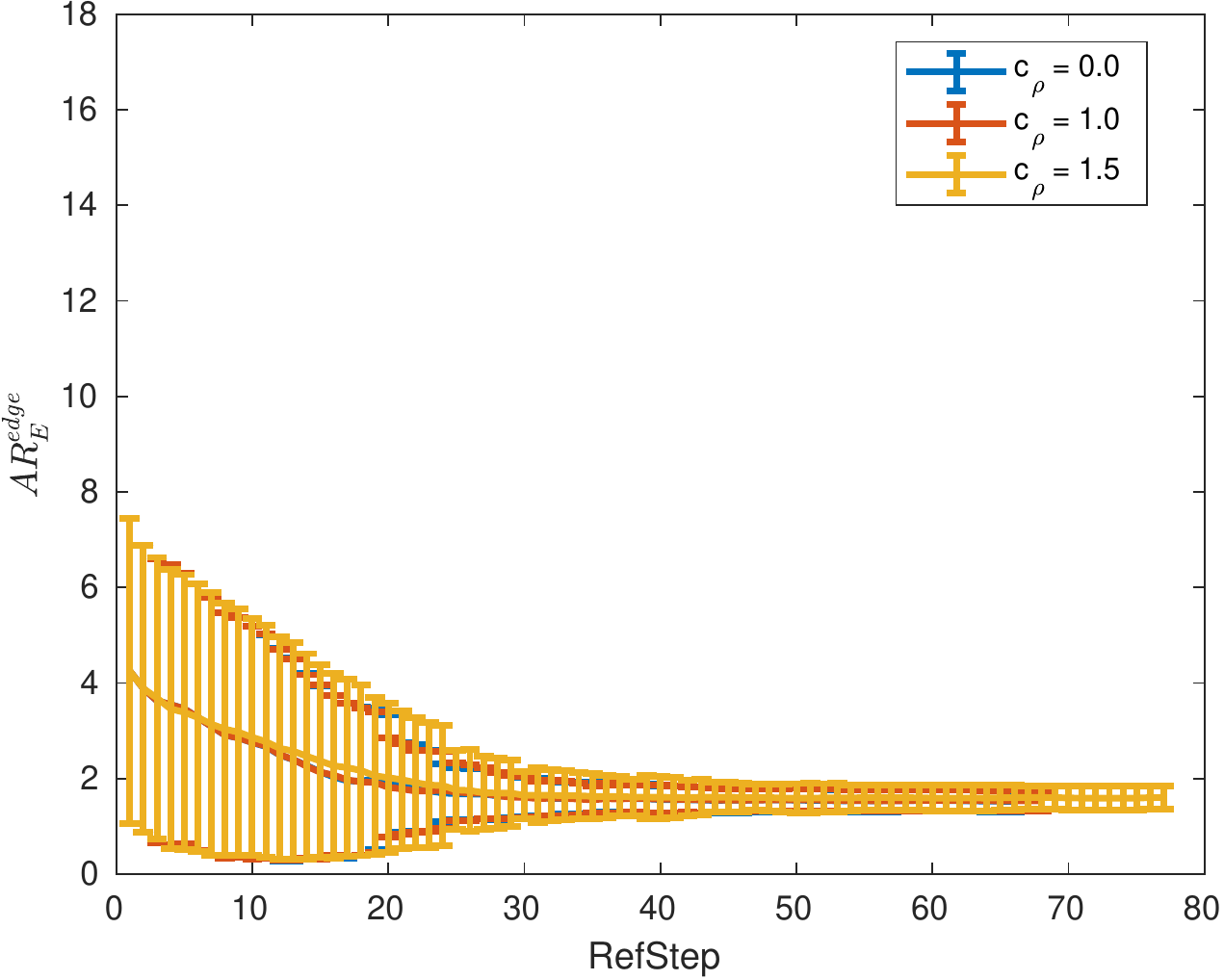}
    \caption{POLY\_MM}
    \label{fig:PMV_MeanARH}
  \end{subfigure}
  \hfill
  \begin{subfigure}[b]{0.49\linewidth}
    \includegraphics[width=\linewidth]{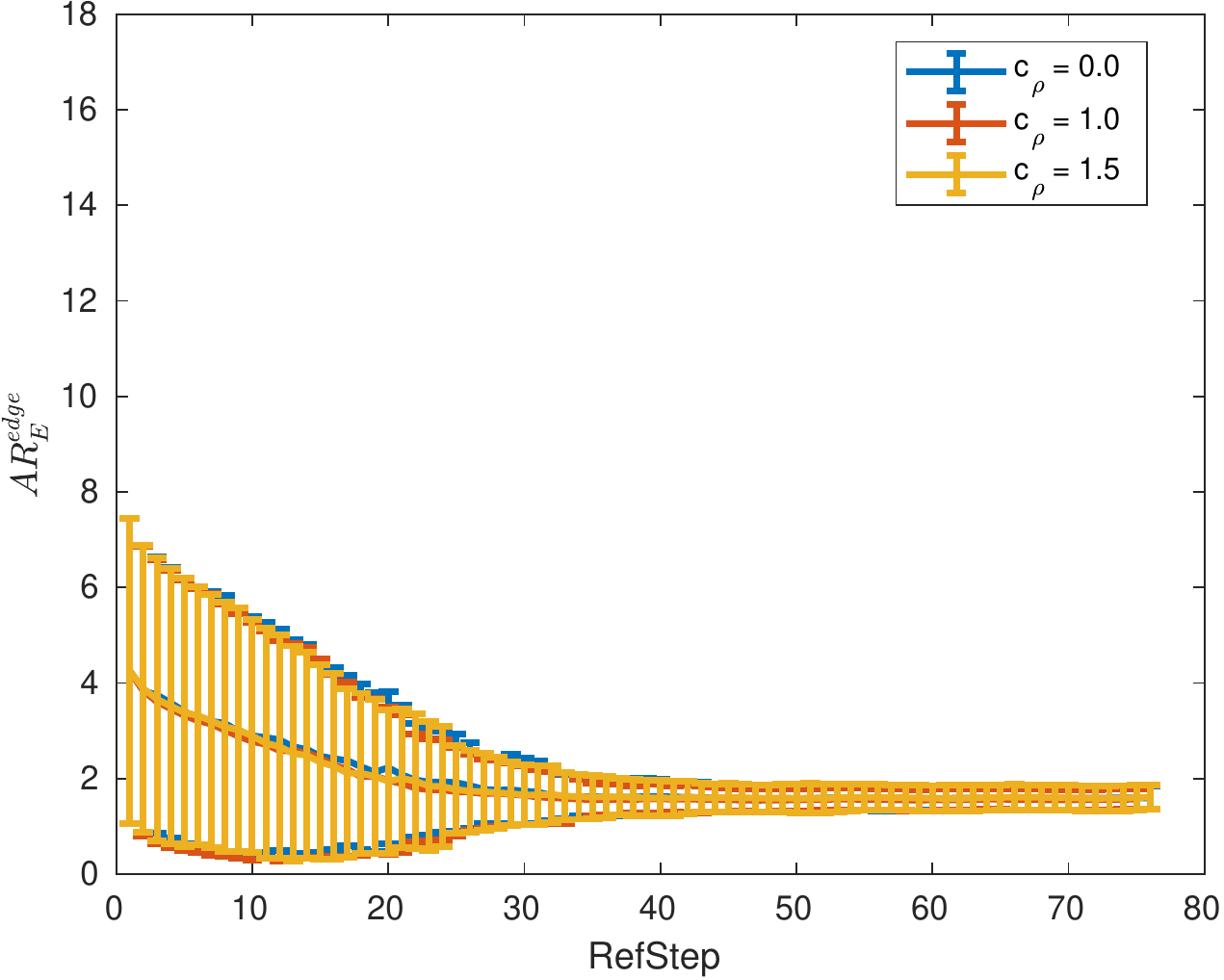}
    \caption{POLY\_LD}
    \label{fig:PLV_MeanARH}
  \end{subfigure}
  \caption{Max $AR^{edge}_E$: (a), (b), (c), (d), Mean $AR^{edge}_E \pm$ standard deviation: (e), (f)}
  \label{fig:MaxARH}
\end{figure}

\begin{figure}
  \centering
    \begin{subfigure}[b]{0.32\linewidth}
    \includegraphics[width=\linewidth]{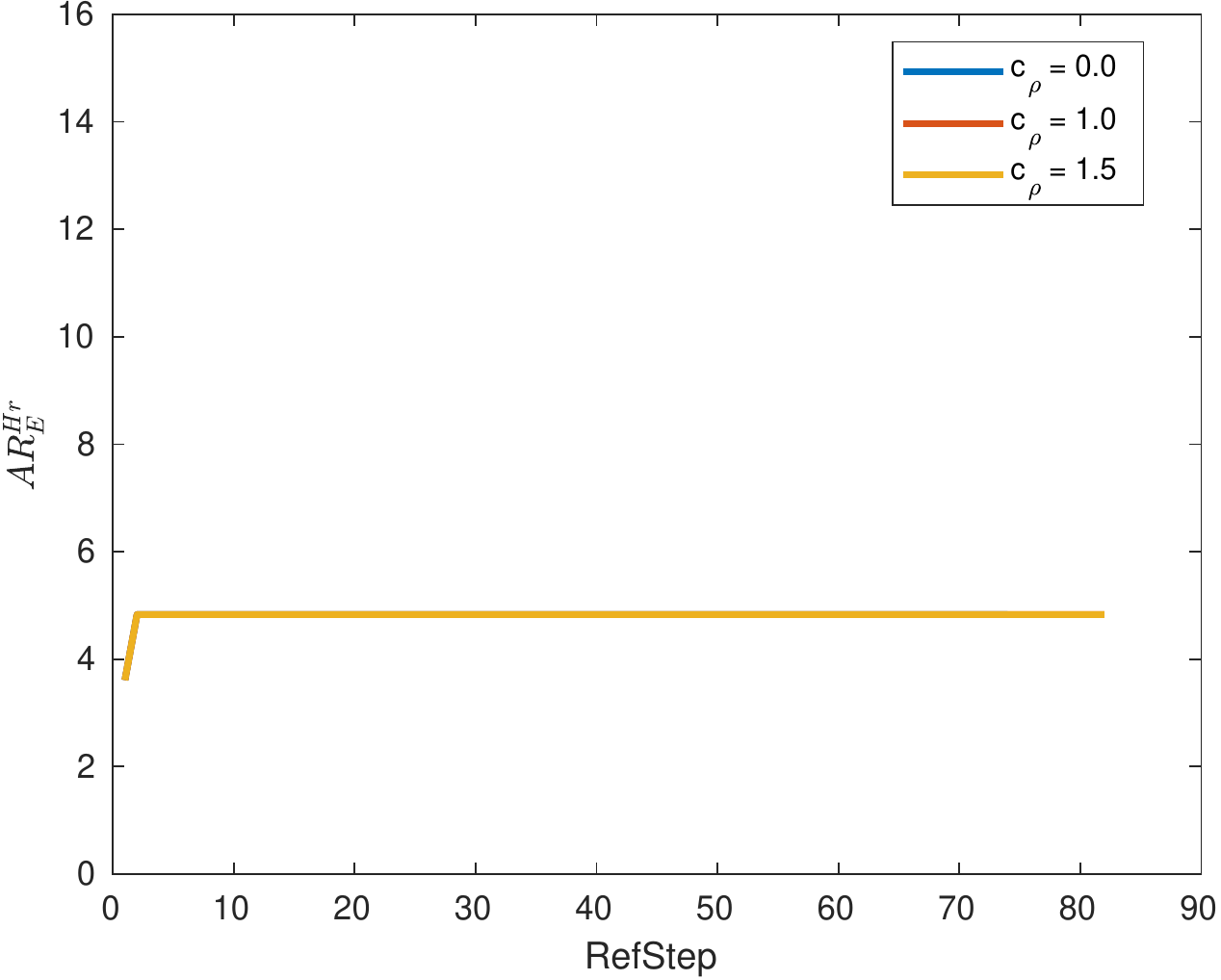}
    \caption{TRAP\_MM}
    \label{fig:TMV_MaxGamma}
  \end{subfigure}
  \hfill
  \begin{subfigure}[b]{0.32\linewidth}
    \includegraphics[width=\linewidth]{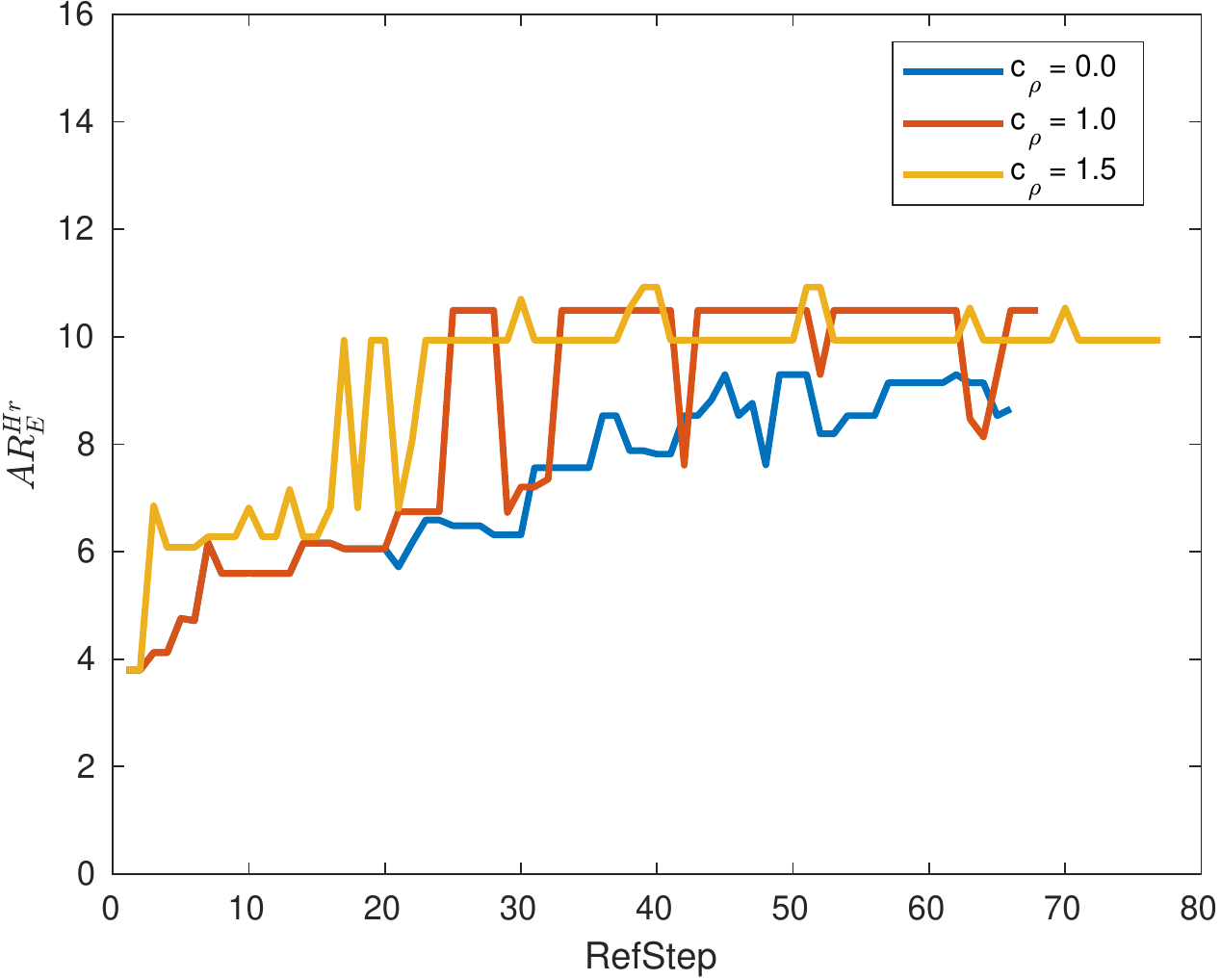}
    \caption{POLY\_MM}
    \label{fig:PMV_MaxGamma}
  \end{subfigure}
  \hfill
  \begin{subfigure}[b]{0.32\linewidth}
    \includegraphics[width=\linewidth]{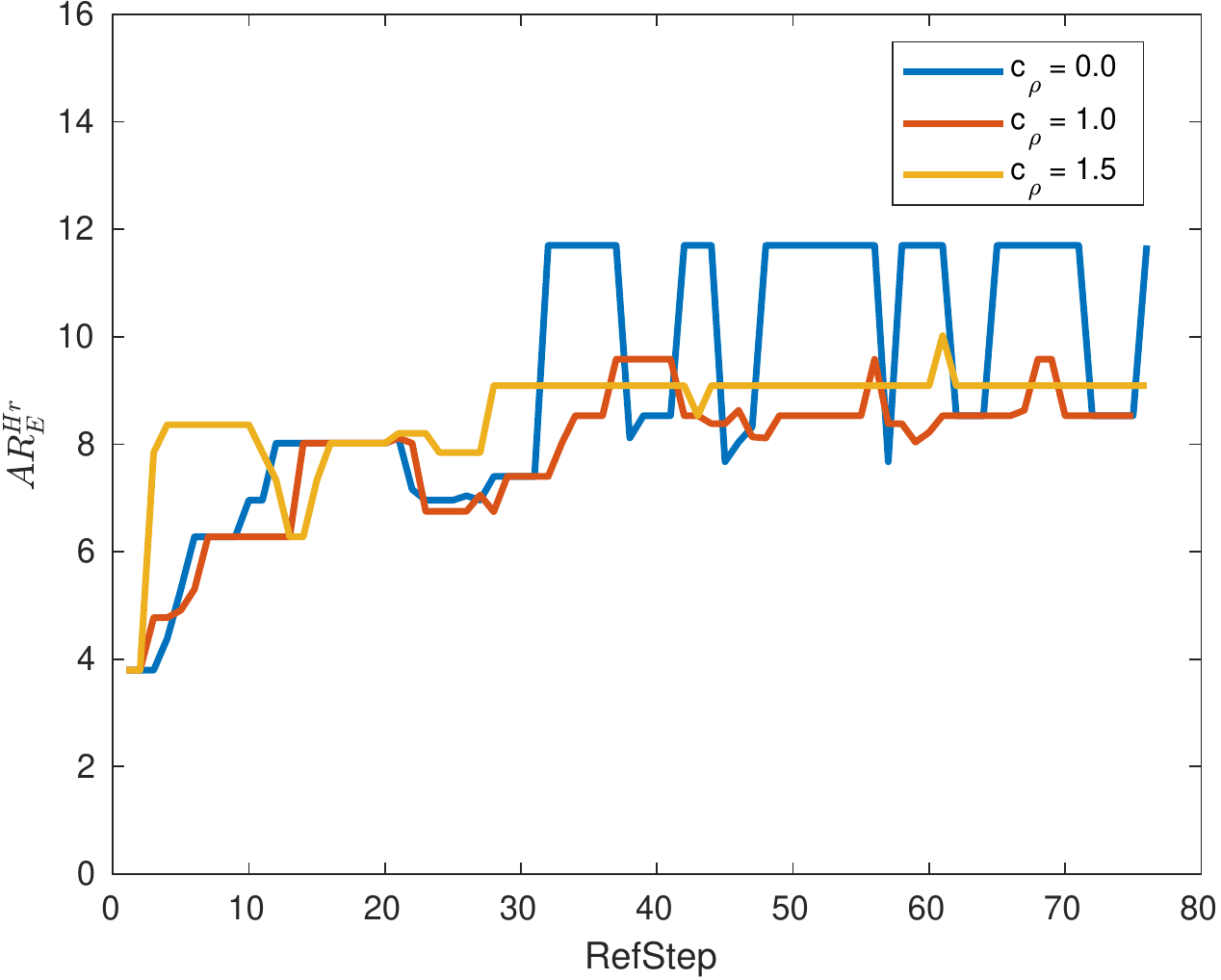}
    \caption{POLY\_LD}
    \label{fig:PLV_MaxGamma}
  \end{subfigure}
  \caption{Max $AR^{Hr}_E$}
  \label{fig:MaxGamma}
  \centering
    \begin{subfigure}[b]{0.32\linewidth}
    \includegraphics[width=\linewidth]{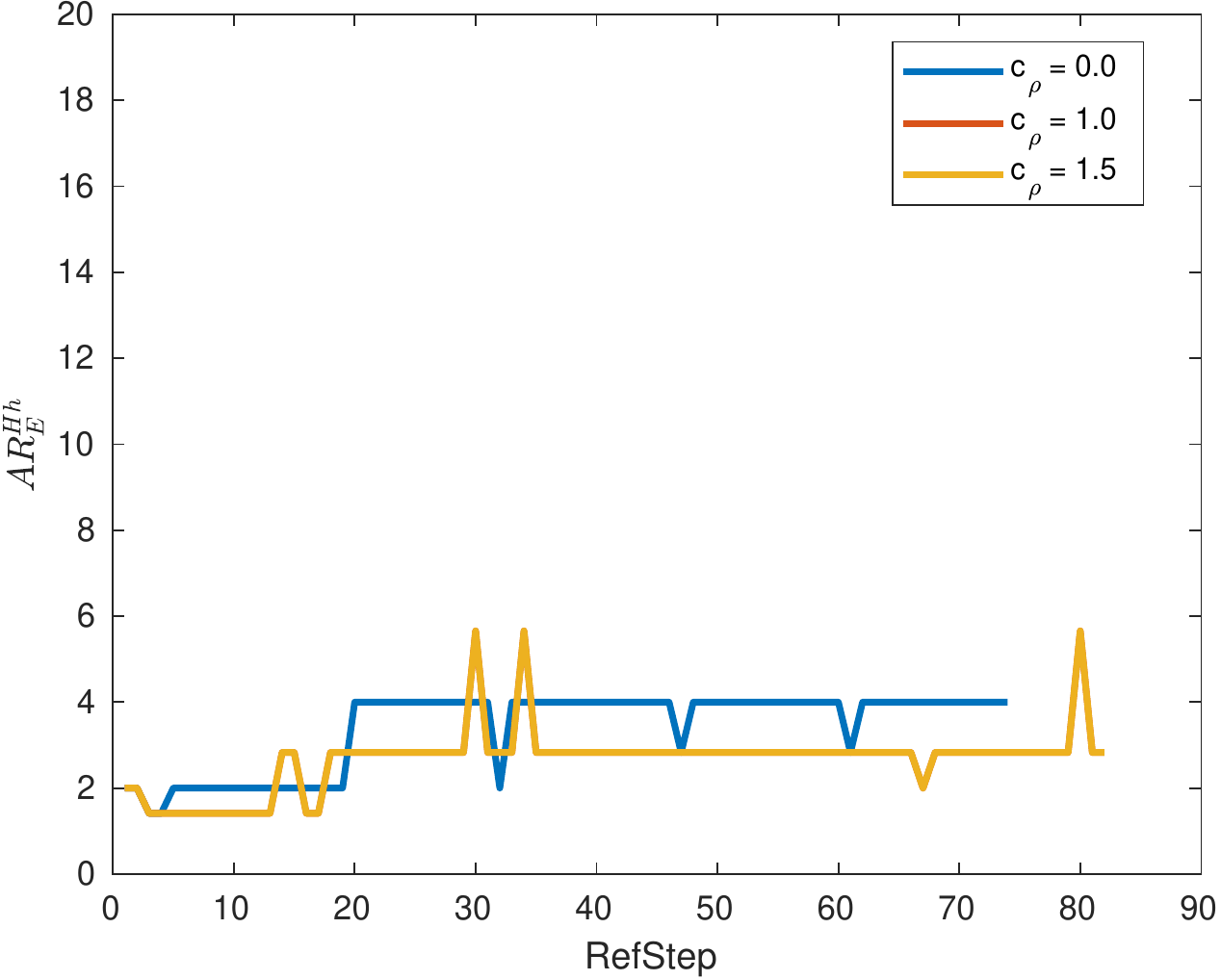}
    \caption{TRAP\_MM}
    \label{fig:TMV_MaxEta}
  \end{subfigure}
  \hfill
  \begin{subfigure}[b]{0.32\linewidth}
    \includegraphics[width=\linewidth]{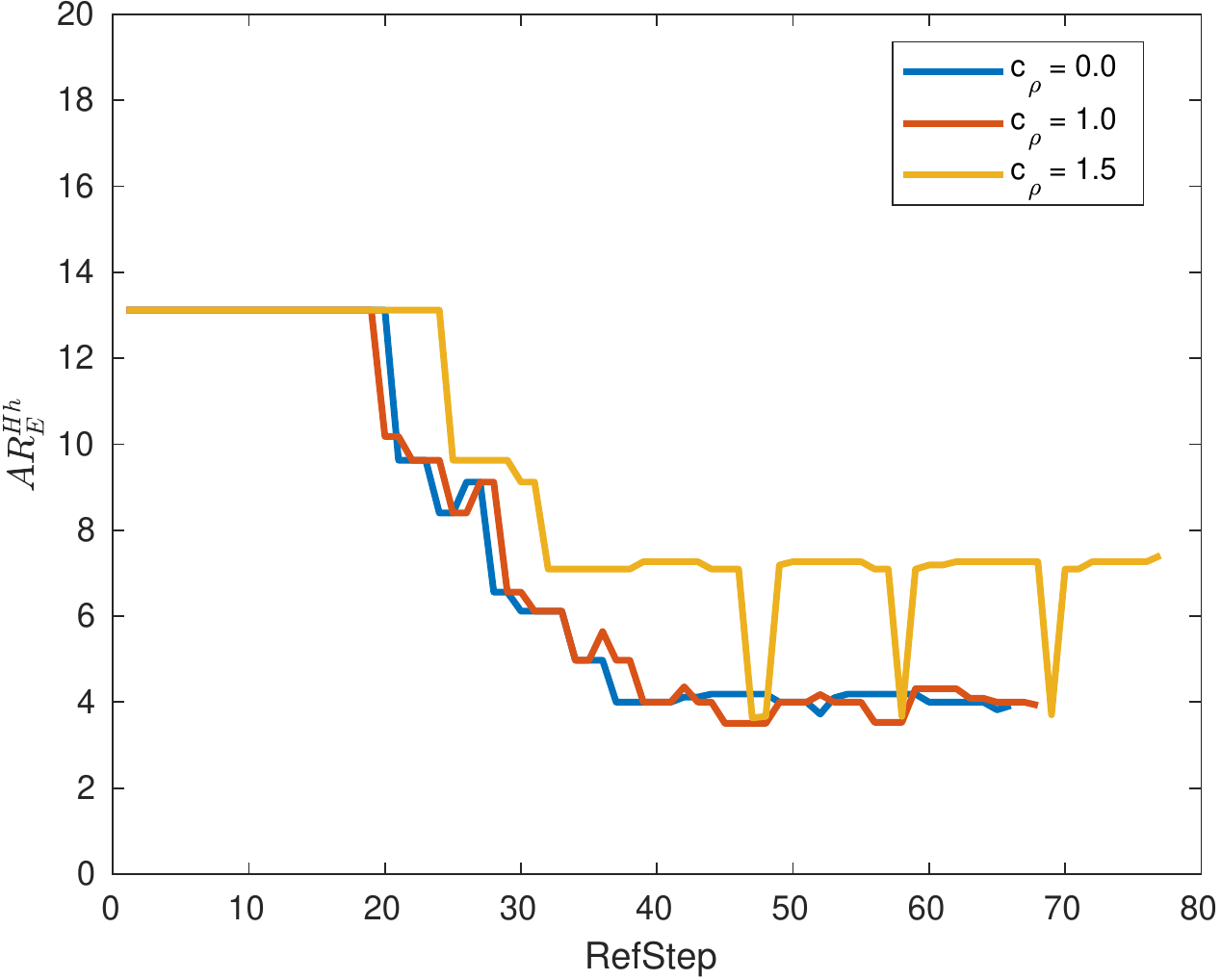}
    \caption{POLY\_MM}
    \label{fig:PMV_MaxEta}
  \end{subfigure}
  \hfill
  \begin{subfigure}[b]{0.32\linewidth}
    \includegraphics[width=\linewidth]{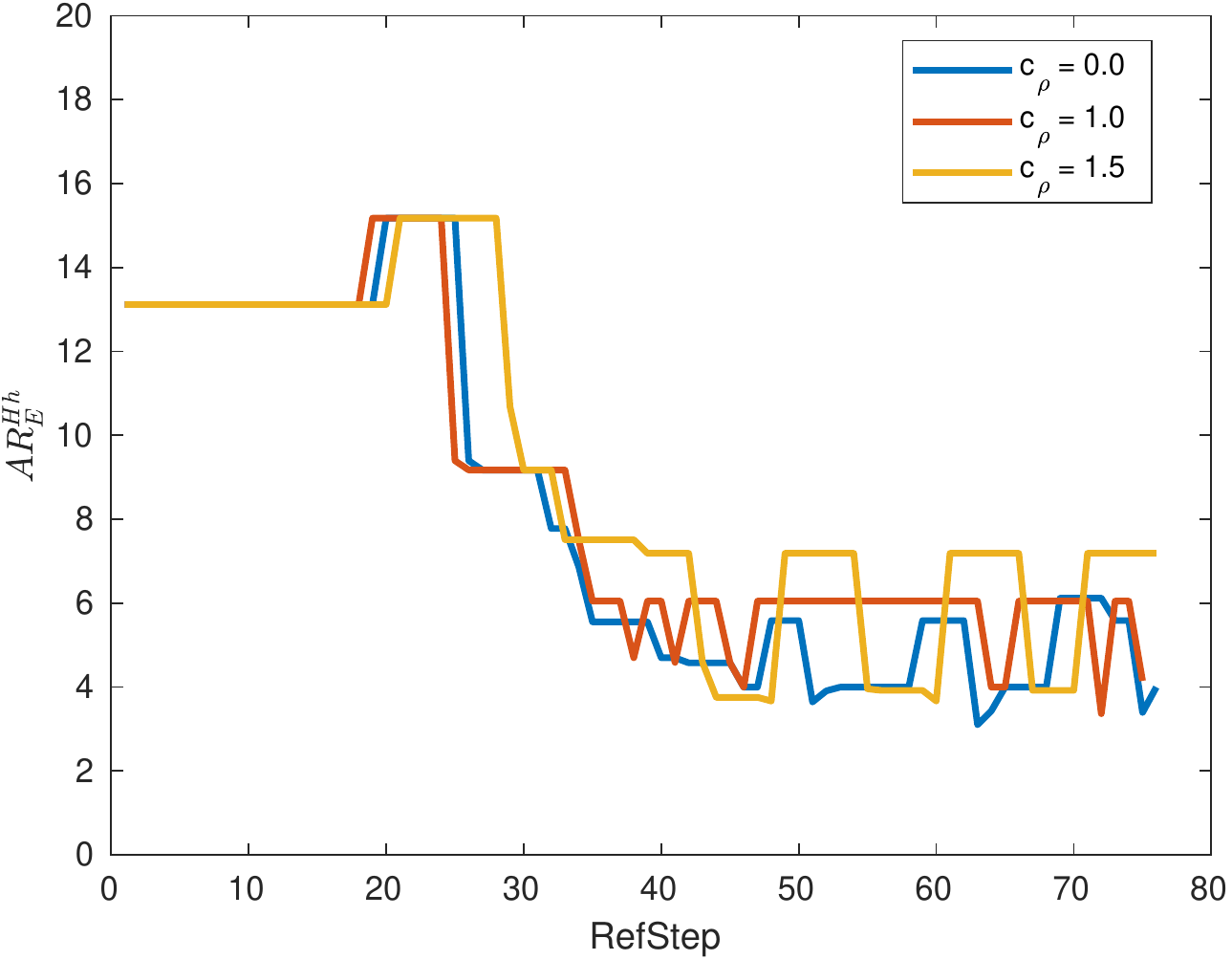}
    \caption{POLY\_LD}
    \label{fig:PLV_MaxEta}
  \end{subfigure}
  \caption{Max $AR^{Hh}_E$}
  \label{fig:MaxEta}
\end{figure}

\begin{figure}
  \centering
  \begin{subfigure}[b]{0.32\linewidth}
    \includegraphics[width=\linewidth]{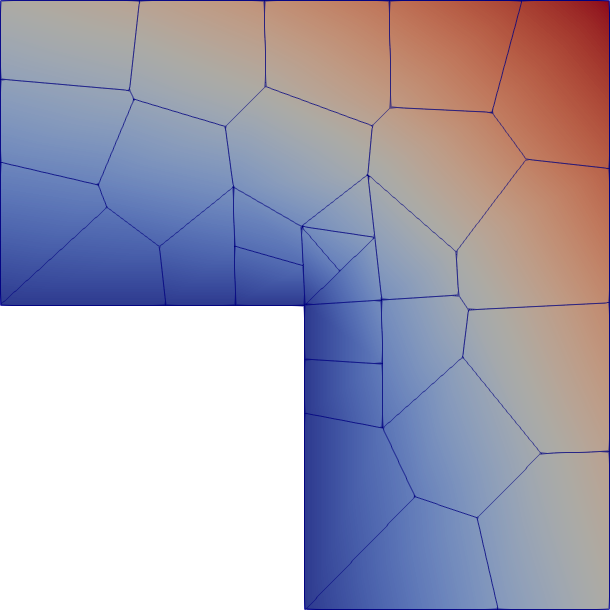}
    \caption{Step 5 $c_\rho = 0.0$}
  \end{subfigure}
  \hfill
  \begin{subfigure}[b]{0.32\linewidth}
    \includegraphics[width=\linewidth]{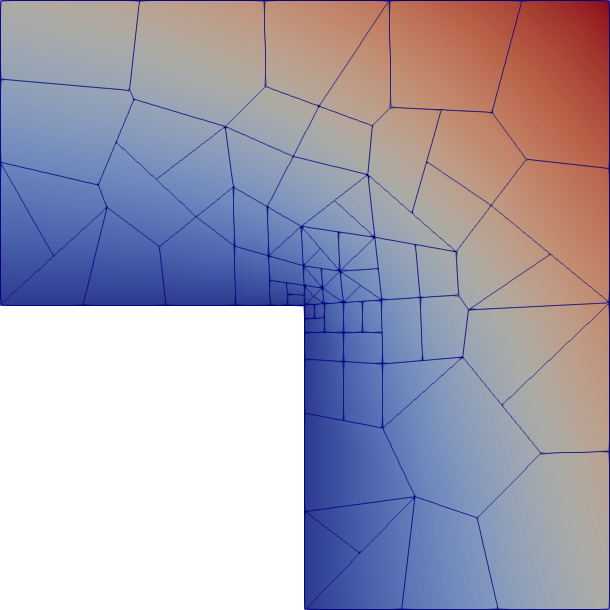}
        \caption{Step 15 $c_\rho = 0.0$}
  \end{subfigure}
  \hfill
  \centering
  \begin{subfigure}[b]{0.32\linewidth}
    \includegraphics[width=\linewidth]{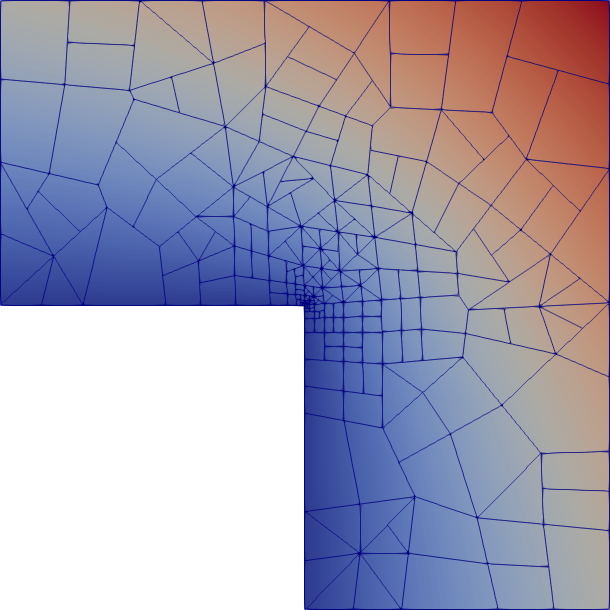}
    \caption{Step 25 $c_\rho = 0.0$}
  \end{subfigure}
  \hfill
  \begin{subfigure}[b]{0.32\linewidth}
    \includegraphics[width=\linewidth]{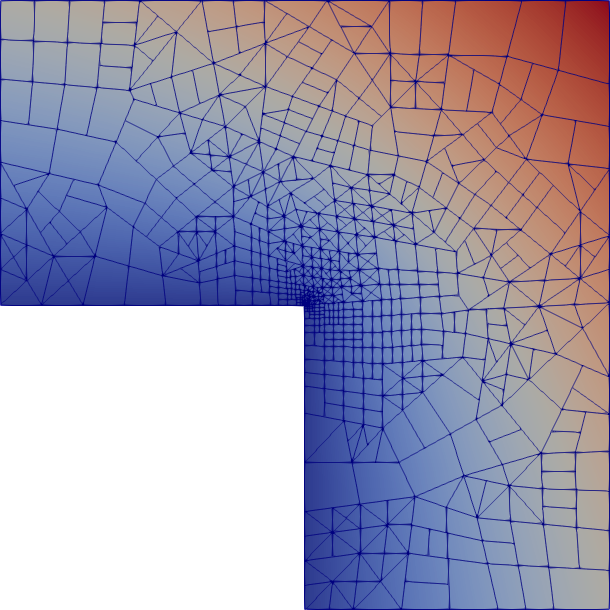}
    \caption{Step 35 $c_\rho = 0.0$}
  \end{subfigure}
  \hfill
  \begin{subfigure}[b]{0.32\linewidth}
    \includegraphics[width=\linewidth]{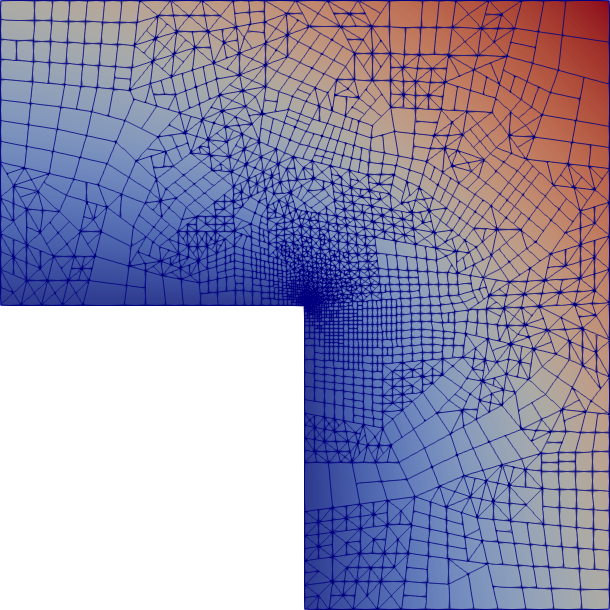}
    \caption{Step 45 $c_\rho = 0.0$}
  \end{subfigure}
  \hfill
  \begin{subfigure}[b]{0.32\linewidth}
    \includegraphics[width=\linewidth]{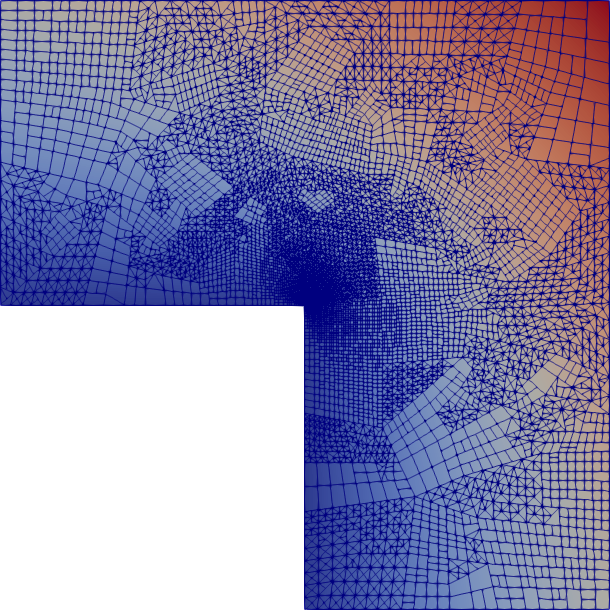}
    \caption{Step 55 $c_\rho = 0.0$}
  \end{subfigure}
  \hfill
  \begin{subfigure}[b]{0.32\linewidth}
    \includegraphics[width=\linewidth]{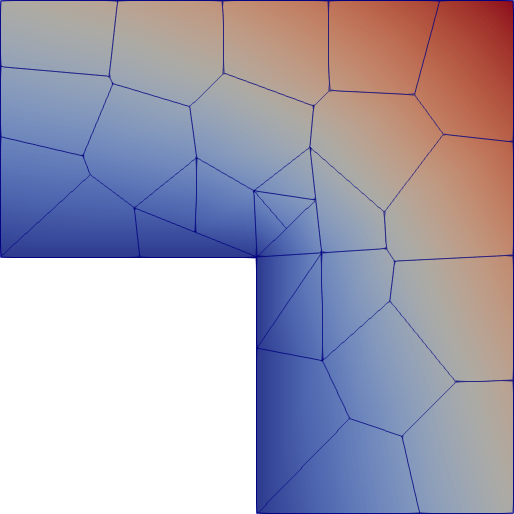}
    \caption{Step 5 $c_\rho = 1.5$}
  \end{subfigure}
  \hfill
  \begin{subfigure}[b]{0.32\linewidth}
    \includegraphics[width=\linewidth]{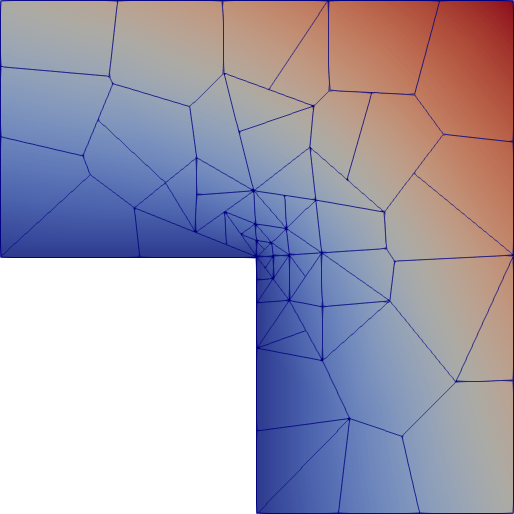}
        \caption{Step 15 $c_\rho = 1.5$}
  \end{subfigure}
  \hfill
  \centering
  \begin{subfigure}[b]{0.32\linewidth}
    \includegraphics[width=\linewidth]{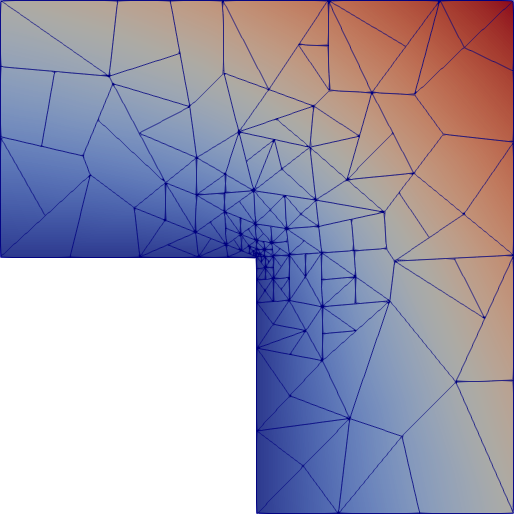}
    \caption{Step 25 $c_\rho = 1.5$}
  \end{subfigure}
  \hfill
  \begin{subfigure}[b]{0.32\linewidth}
    \includegraphics[width=\linewidth]{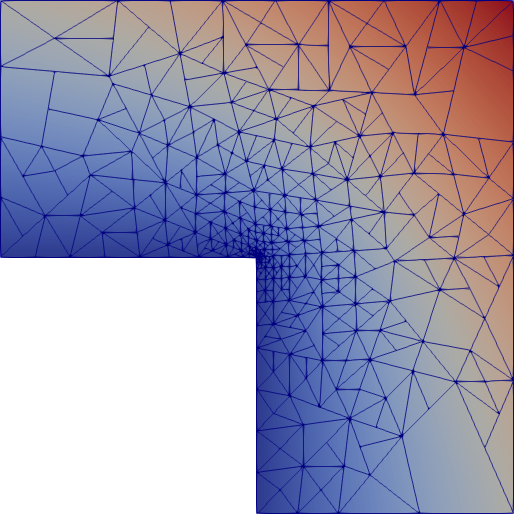}
    \caption{Step 35 $c_\rho = 1.5$}
  \end{subfigure}
  \hfill
  \begin{subfigure}[b]{0.32\linewidth}
    \includegraphics[width=\linewidth]{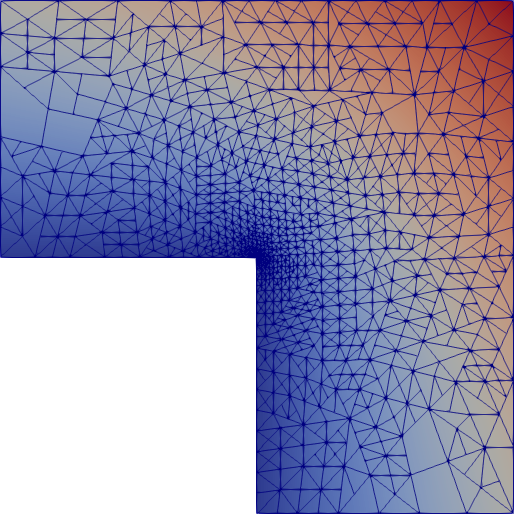}
    \caption{Step 45 $c_\rho = 1.5$}
  \end{subfigure}
  \hfill
  \begin{subfigure}[b]{0.32\linewidth}
    \includegraphics[width=\linewidth]{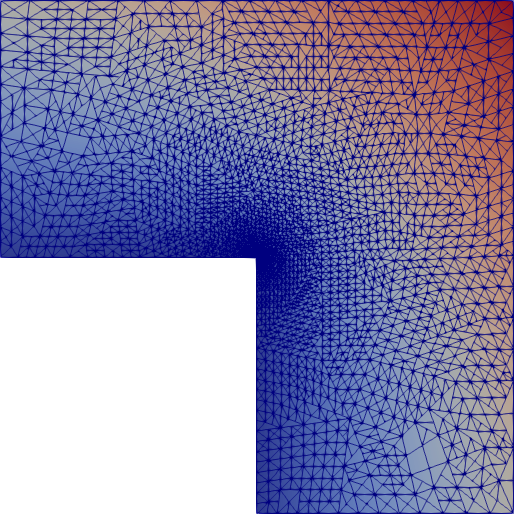}
    \caption{Step 55 $c_\rho = 1.5$}
  \end{subfigure}
 \caption{Refinement steps POLY\_MM test with $c_\rho = 0.0, 1.5$.}
  \label{fig:Polygon-AnalysisMesh}
\end{figure}

\begin{figure}
  \centering
  \begin{subfigure}[b]{0.49\linewidth}
    \includegraphics[width=\linewidth]{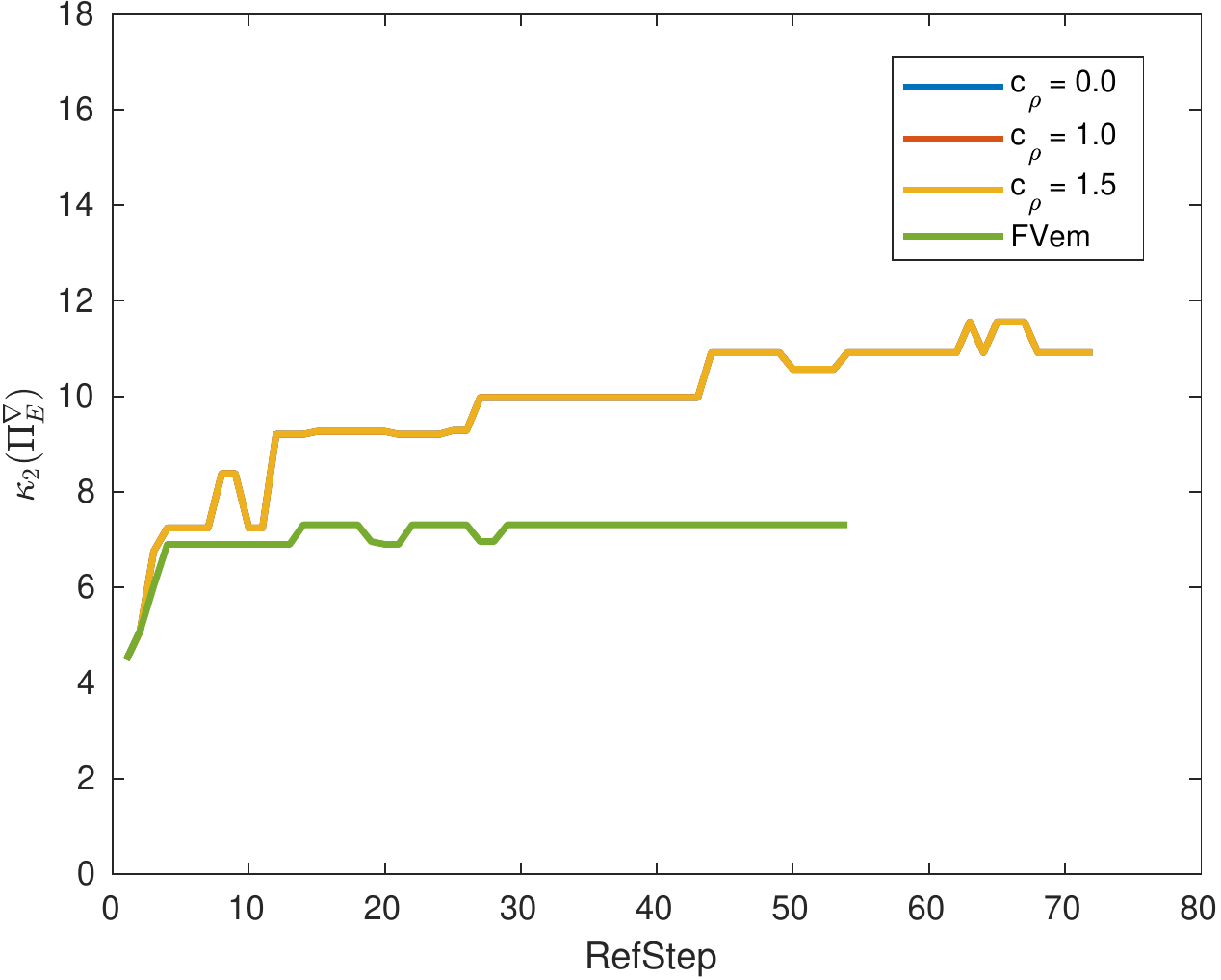}
    \caption{TRI\_MM-LE}
    \label{fig:TMVLF_MaxConditionNumberPiNabla}
  \end{subfigure}
  \hfill
    \begin{subfigure}[b]{0.49\linewidth}
    \includegraphics[width=\linewidth]{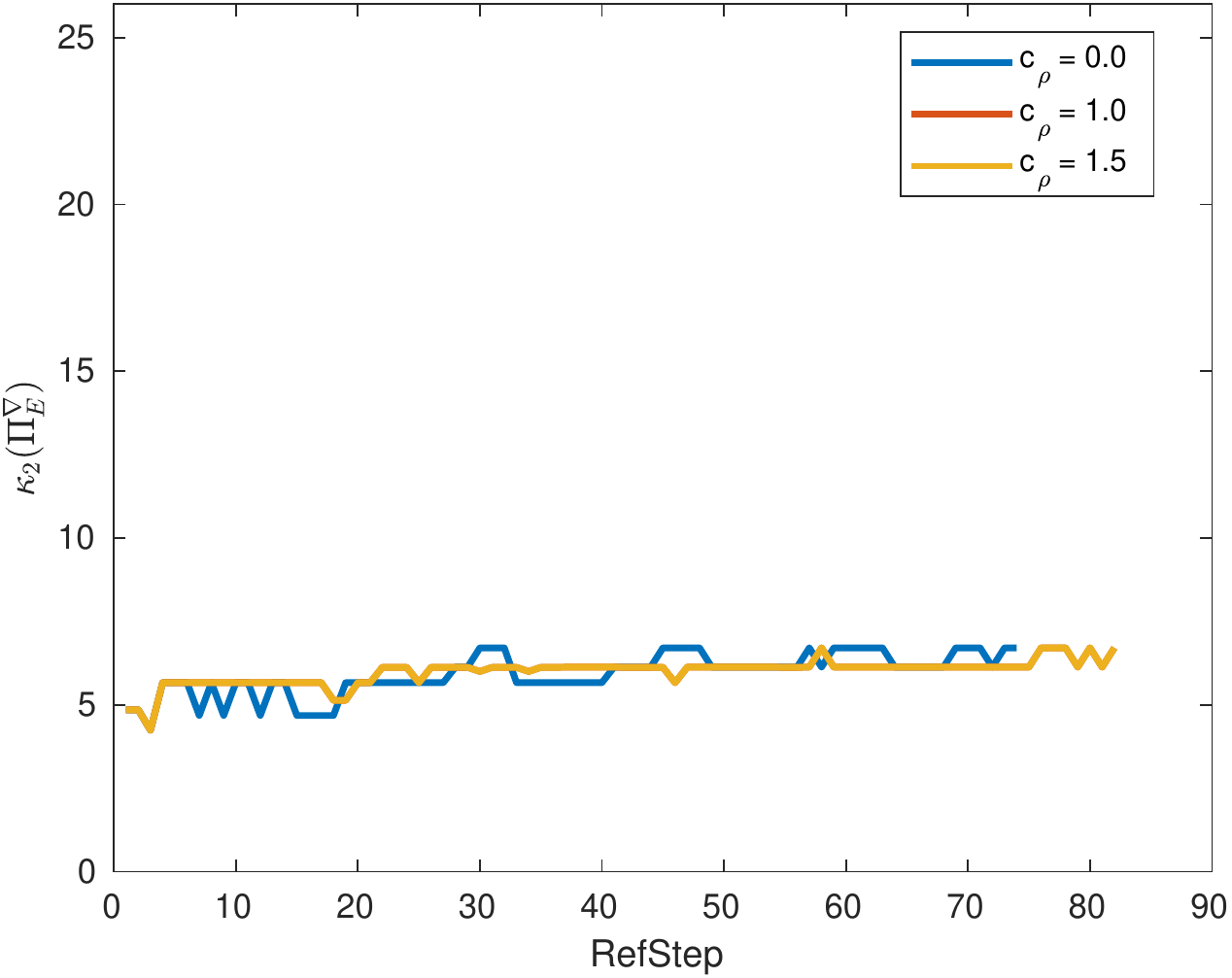}
    \caption{TRAP\_MM}
    \label{fig:TMV_MaxConditionNumberPiNabla}
  \end{subfigure}
  \hfill
  \begin{subfigure}[b]{0.49\linewidth}
    \includegraphics[width=\linewidth]{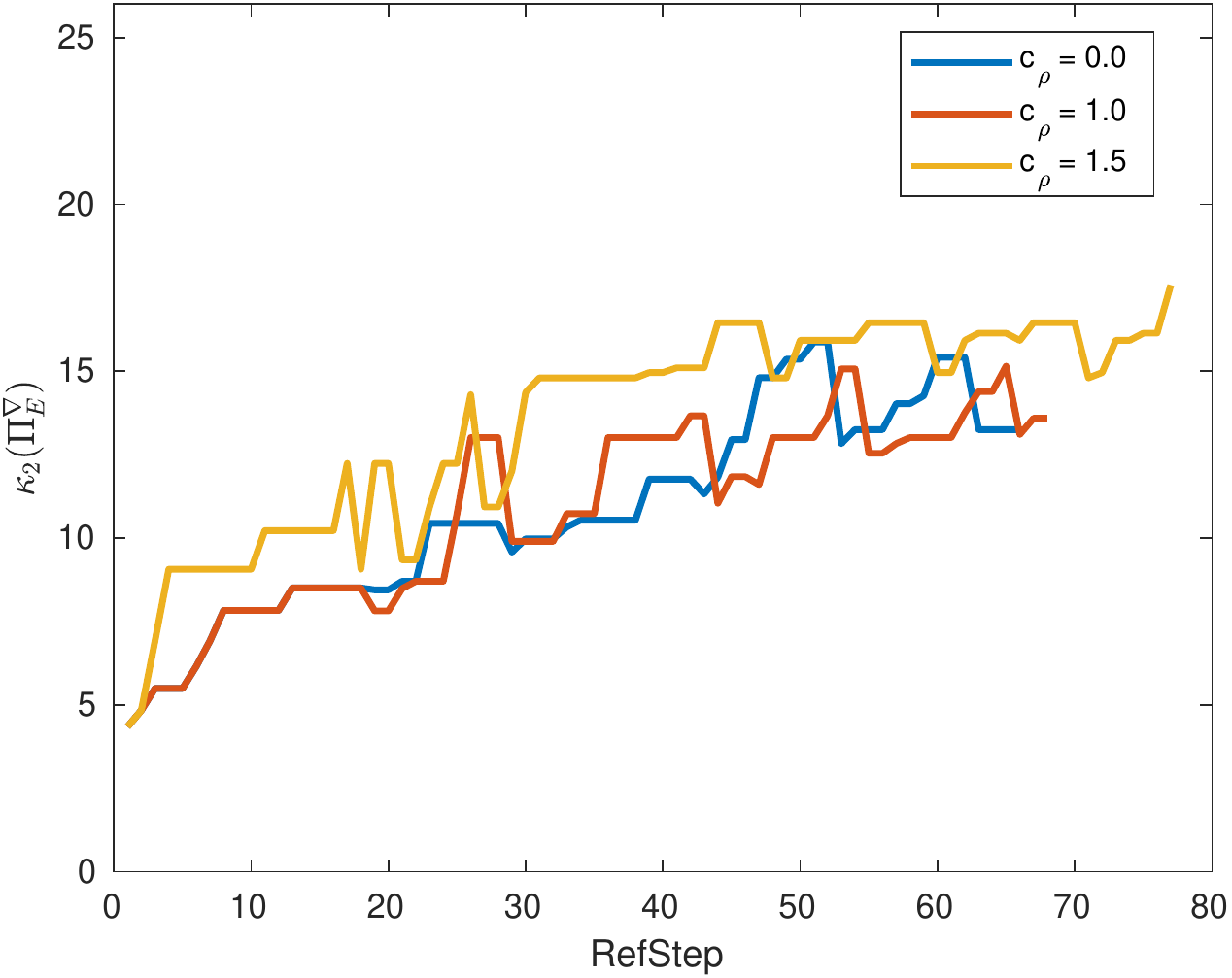}
    \caption{POLY\_MM}
    \label{fig:PMV_MaxConditionNumberPiNabla}
  \end{subfigure}
  \hfill
  \begin{subfigure}[b]{0.49\linewidth}
    \includegraphics[width=\linewidth]{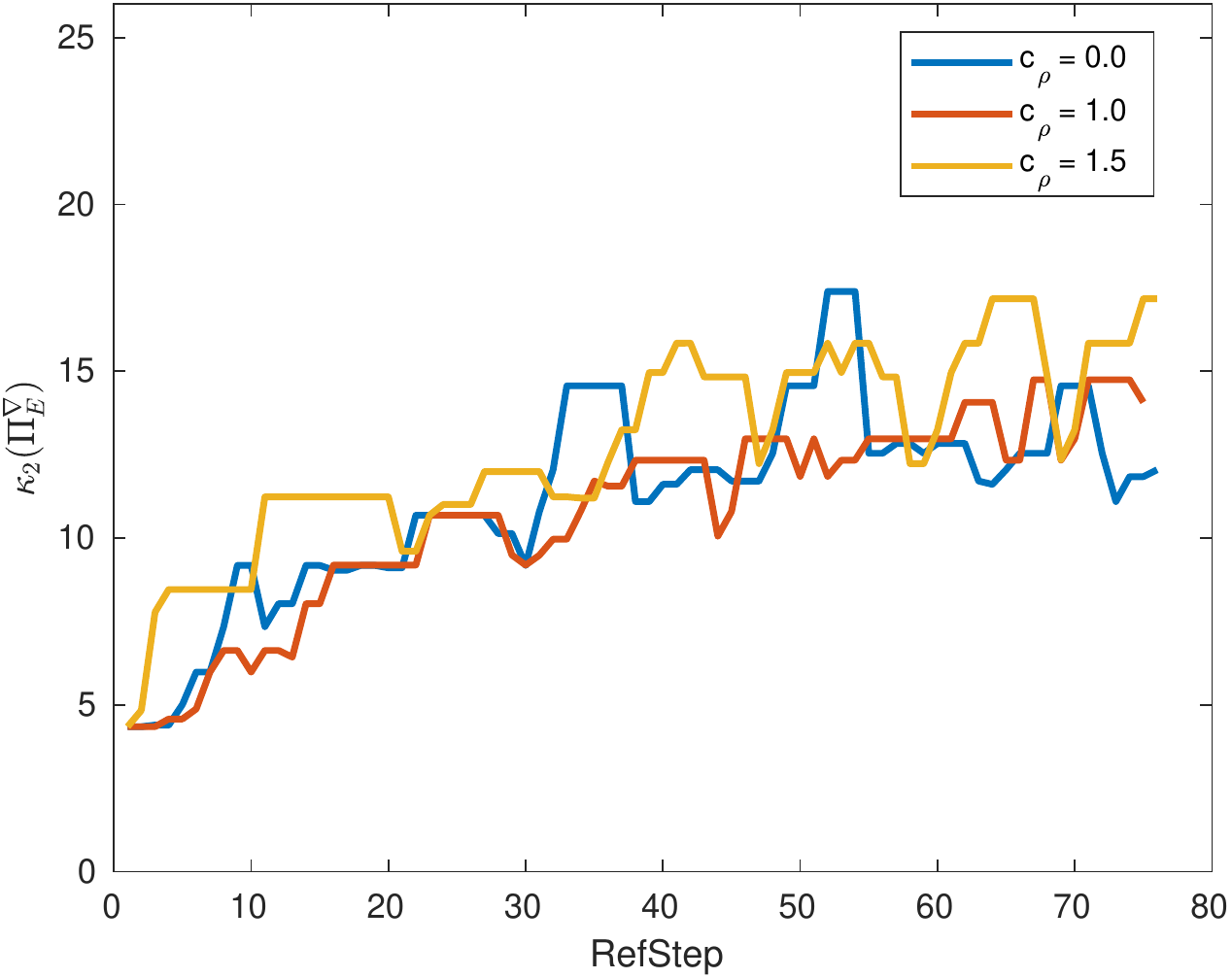}
    \caption{POLY\_LD}
    \label{fig:PLV_MaxConditionNumberPiNabla}
  \end{subfigure}
  \hfill
  \begin{subfigure}[b]{0.49\linewidth}
    \includegraphics[width=\linewidth]{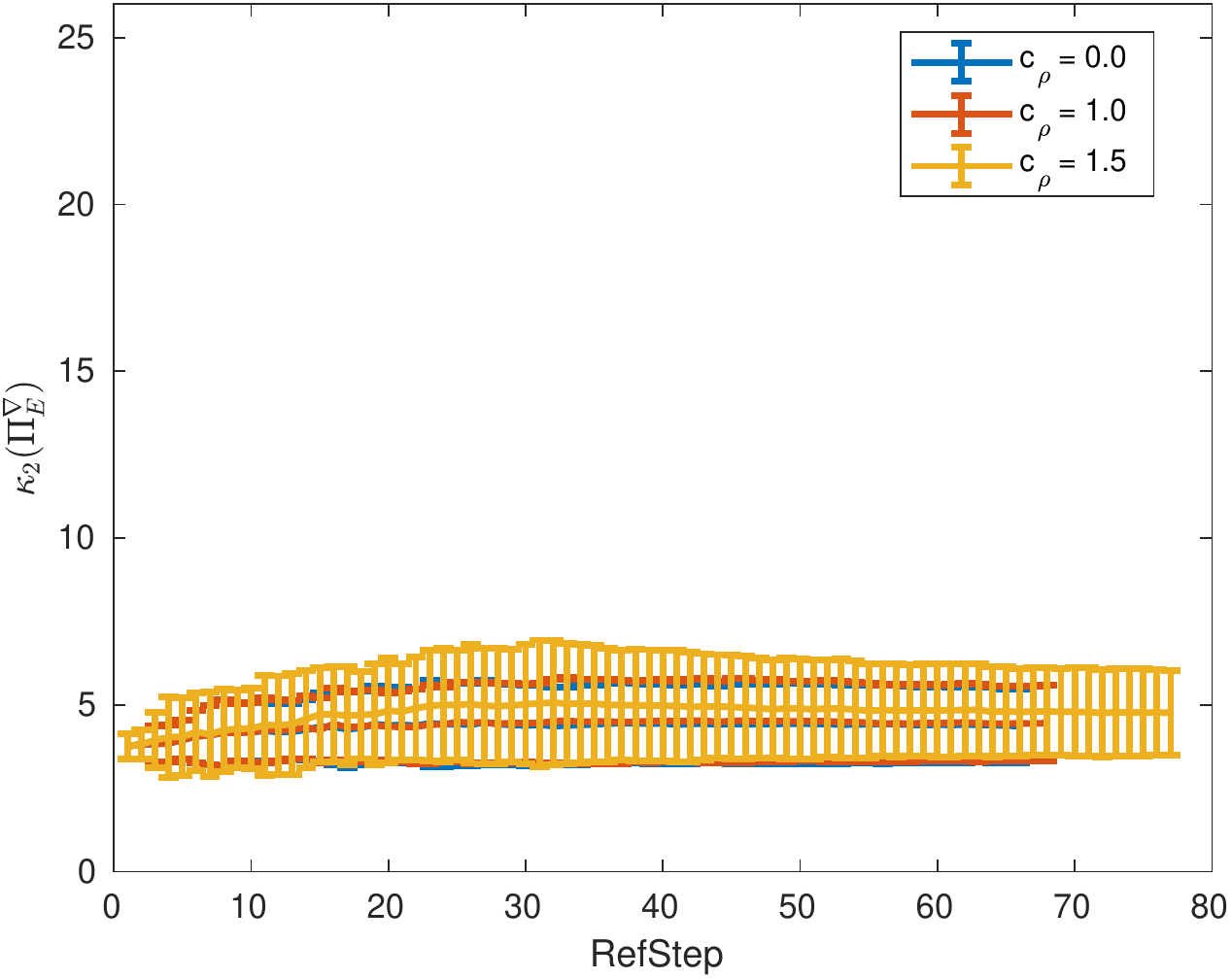}
    \caption{POLY\_MM}
    \label{fig:PMV_MeanConditionNumberPiNabla}
  \end{subfigure}
  \hfill
  \begin{subfigure}[b]{0.49\linewidth}
    \includegraphics[width=\linewidth]{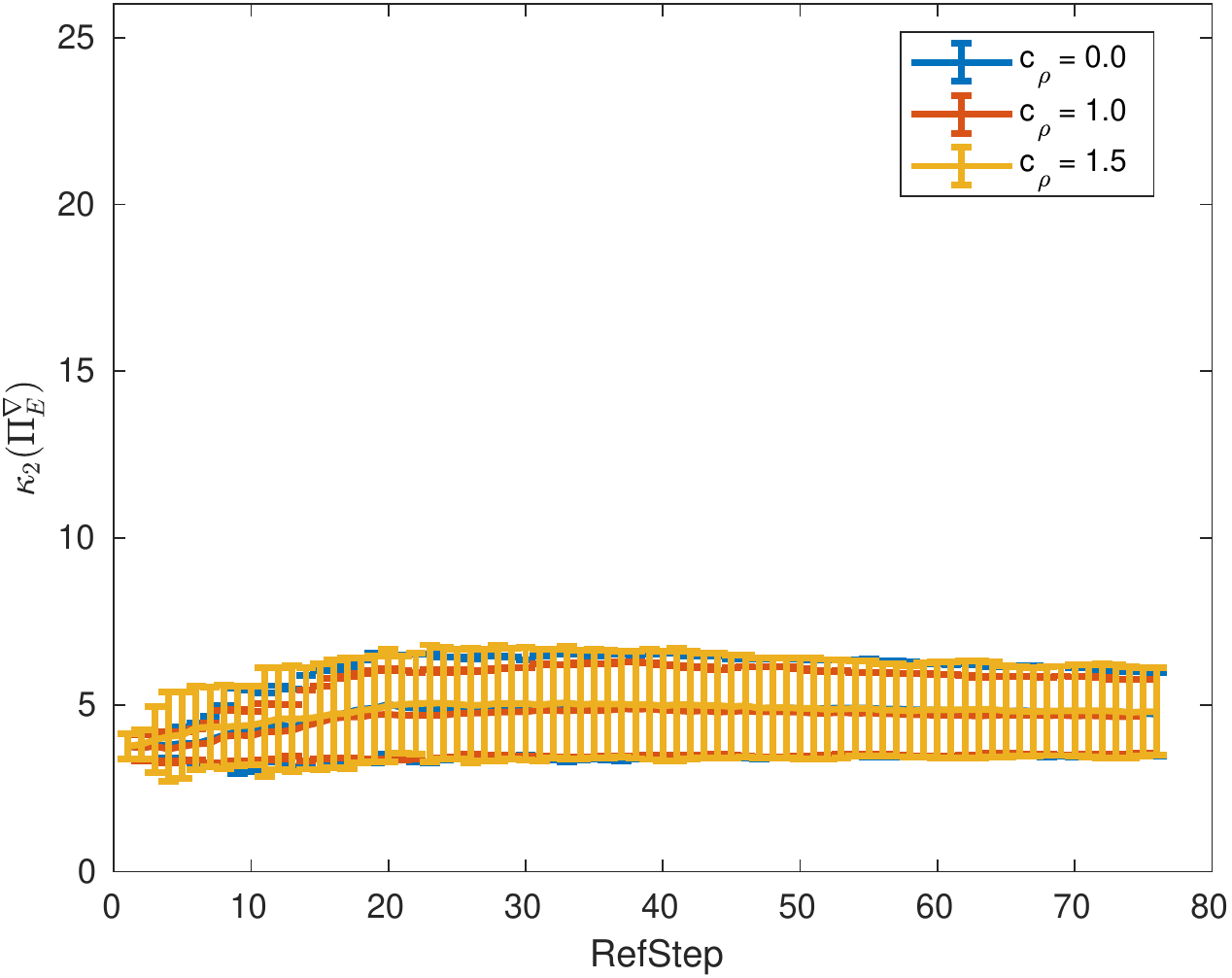}
    \caption{POLY\_LD}
    \label{fig:PLV_MeanConditionNumberPiNabla}
  \end{subfigure}
  \caption{Max $\kappa_2(\Pi^\nabla)$: (a), (b), (c), (d), Mean $\kappa_2(\Pi^\nabla)\pm$ standard deviation: (e), (f)}
  {}
  \label{fig:ConditionNumberPiNabla}
\end{figure}

In Figure~\ref{fig:MaxARR} we display the behaviour of $AR_E^{Rr}$, in subfigure \ref{fig:TMVLF_MaxARR}-\ref{fig:PLV_MaxARR} the maximum value of $AR_E^{Rr}$ for each mesh is plotted, and in Subfigures~\ref{fig:PMV_MeanARR} and \ref{fig:PLV_MeanARR} the mean values over each mesh and the standard deviation are displayed. We can see that the maximum values are quite stable and during the refinement process does not increase so much being the quality of the mesh stable for all the elements.

The mean values are slightly increasing in the first iterations and then almost constant during refinement. The almost constant standard deviation and  mean value suggest that, also if the maximum is increasing, the number of elements for which we have this higher values are very few and most of the elements of the mesh have an improved quality. In Figure~\ref{fig:MaxARH} we display the same type of information concerning $AR^{edge}_E$. The behaviour of this parameter is even better with respect to $AR_E^{Rr}$ and the mean values of this aspect ratio are decreasing and the standard deviation strongly decreasing denoting a very good effect of the refinement on the produced polygonal meshes. In Figures~\ref{fig:MaxGamma} and \ref{fig:MaxEta} we report some plot of the aspect ratios $ AR_E^{Hr},\ AR_E^{Hh}$ with a behaviour very similar to $AR_E^{Rr},\ AR^{edge}_E$, respectively.

In Figures~\ref{fig:NumQuad} and \ref{fig:NumTri} we compare the fraction of cells that during refinement become triangles or quadrilateral polygons for TRAP\_MM mesh and POLY\_MM and POLY\_LD meshes. As previously noted the TRAP mesh with $c_\rho=0.0$ tend to produce quadrilateral cells, whereas in all the other cases the mesh is progressively and automatically transformed in a good quality triangular mesh.

In Figure~\ref{fig:Polygon-AnalysisMesh} we display a sequence of meshes generated by the refinement algorithm. Figures~\ref{fig:Polygon-AnalysisMesh}a-f describe the behaviour of the algorithm with $c_\rho=0.0$, whereas Figures~\ref{fig:Polygon-AnalysisMesh}g-l provide the meshes generated with $c_\rho=1.5$. In  Figure~\ref{fig:Polygon-AnalysisMesh}f we can clearly observe a large number of cells with four or three edges and that some cells have an hanging node. In  Figure~\ref{fig:Polygon-AnalysisMesh}l most of the cells are triangles. 

In Figure~\ref{fig:ConditionNumberPiNabla}, we analyse the behaviour of the spectral condition number of the projector $\Pi^\nabla$. Although the maximum condition number is slightly increasing for polygonal meshes all the values are in the order of ten and can be considered stable during the refinement process.
Comparing the plots of $\kappa_2(\Pi^\nabla)$ with the behaviour of $AR^{Rr}_E$ and $AR^{Hh}_E$ reported in Figures~\ref{fig:MaxARR} and \ref{fig:MaxGamma}, respectively, we can see that $\kappa_2(\Pi^\nabla)$ is strongly related to the geometrical shape of the cell.
As for the other quality parameter the mean value is almost constant as well.
The condition number in the case FVem is not the identity matrix because it is the change of basis between the scaled monomial basis and the triangle Lagrangian basis.

\begin{figure}
  \centering
    \begin{subfigure}[b]{0.32\linewidth}
    \includegraphics[width=\linewidth]{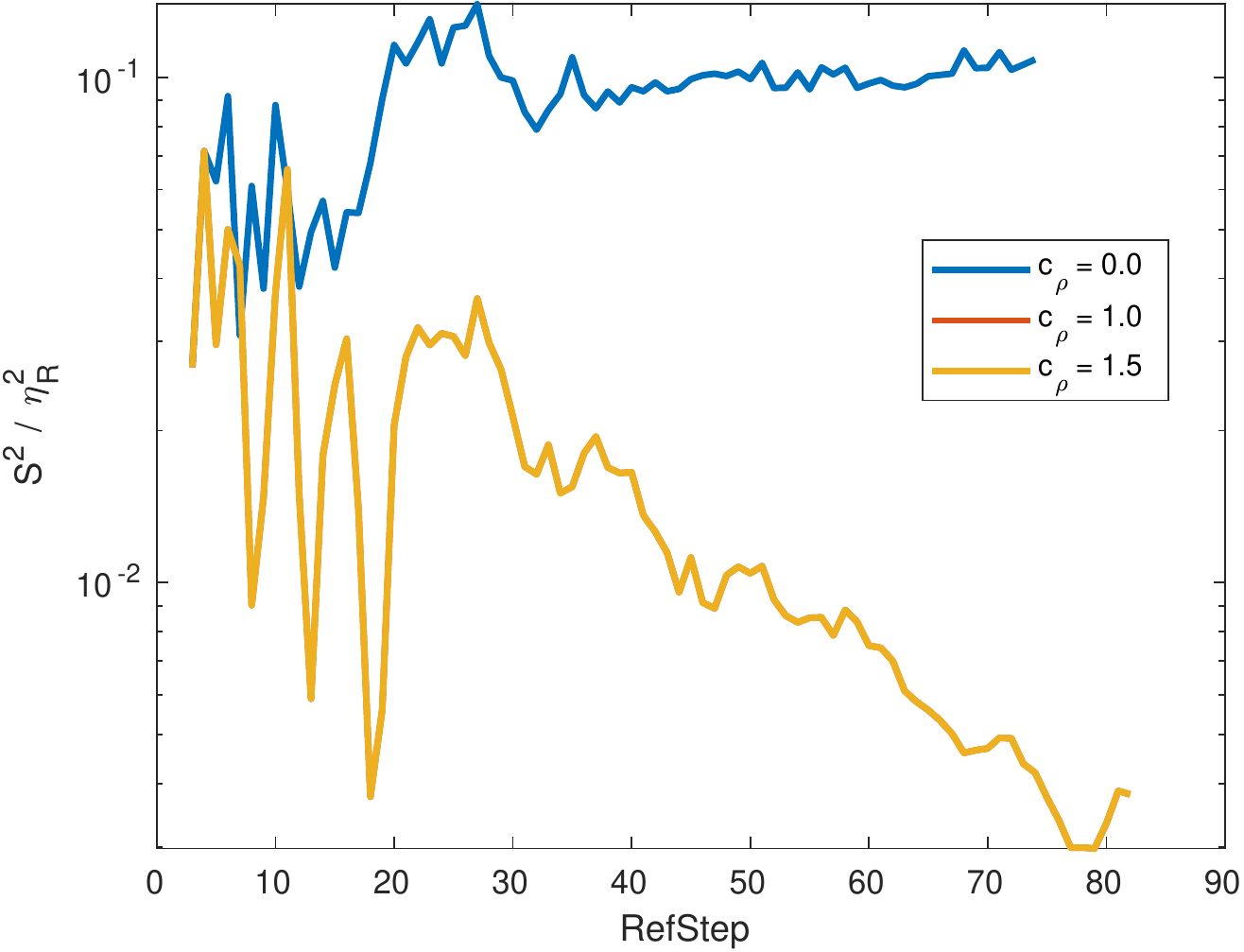}
    \caption{TRAP\_MM}
    \label{fig:TMV_StabEstRatio}
  \end{subfigure}
   \begin{subfigure}[b]{0.32\linewidth}
    \includegraphics[width=\linewidth]{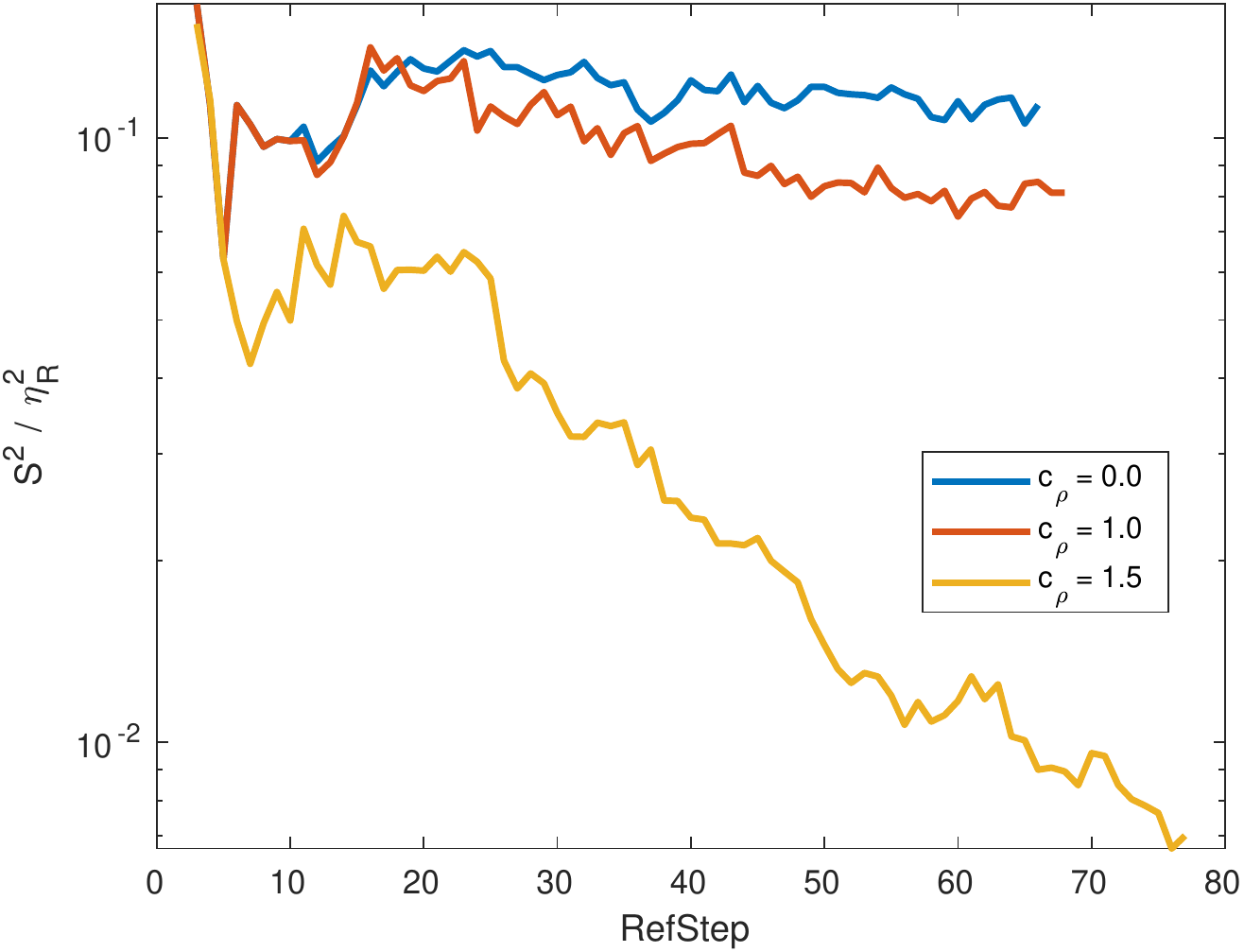}
    \caption{POLY\_MM}
    \label{fig:PMV_StabEstRatio}
  \end{subfigure}
  \begin{subfigure}[b]{0.32\linewidth}
    \includegraphics[width=\linewidth]{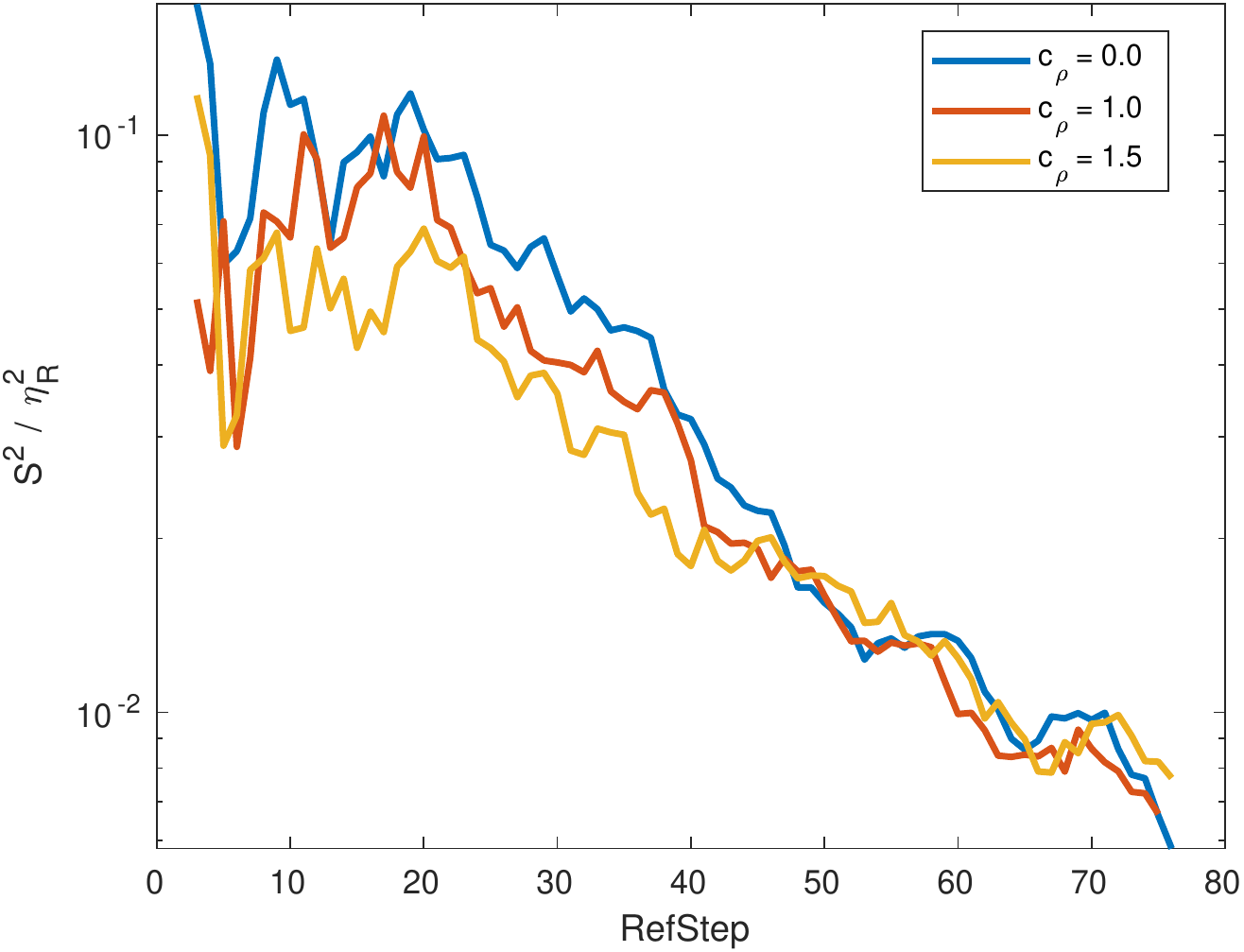}
    \caption{POLY\_LD}
    \label{fig:PLV_StabEstRatio}
  \end{subfigure}
  \caption{Ratio between $\vemstab[]{u_\delta}{u_\delta}^2$ and $\eta^2_R$}
   \label{fig:StabEstRatio}
\end{figure}

The estimator used in the numerical experiments is equivalent to the error provided some mesh properties are assumed. In \cite{BBapost} these assumptions are discussed and the stability of a suitable quasi interpolation operator is numerically investigated. In \cite{daveiga2021adaptive} is proved that the VEM stabilization term, present in the estimator in \cite{Cangiani_et_al-apost:2017}, can be bounded by the estimator \eqref{eq:defetaR} on triangular meshes with hanging nodes.
 In Figure \ref{fig:StabEstRatio} we report the behaviour of the ratio between the stabilization term computed with the solution and the squared estimator along the refinement iterations. The decay of this ratio confirm the negligible role of the stabilization terms in the error estimator within the proposed framework that converges to triangular or quadrilateral (including triangles with hanging nodes) meshes. In the case TRAP\_MM (Figure \ref{fig:TMV_StabEstRatio}) we can see a fast decay for the values of $c_\rho$ that converge to triangular meshes (see Figures \ref{fig:NumTri}  and  \ref{fig:NumQuad}), whereas we find a larger and stable ratio for $c_\rho=0.0$ converging to mixed triangular-quadrilateral meshes. The same behaviour is confirmed by  Figure \ref{fig:PMV_StabEstRatio} where the case  $c_\rho=1.5$ is rapidly converging to an almost triangular mesh and the ratio $\vemstab[]{u_\delta}{v_\delta}^2/\eta^2_R$ is decreasing, whereas the cases $c_\rho=0.0$ and $c_\rho=1.0$ do not decrease due to the large number of non triangular cells. The refinement strategy based on the longest diagonal (Figure~\ref{fig:PLV_StabEstRatio}) is converging towards a triangular mesh for all the values of $c_\rho$ and the ratio is decreasing for all the values.

%% file: num_res_Polygon.tex
\section{Uniform refinement of a regular polygon}
\label{sec:UnifPolyg}

\begin{figure}
  \centering
  \begin{subfigure}[b]{0.32\linewidth}
    \includegraphics[width=\linewidth]{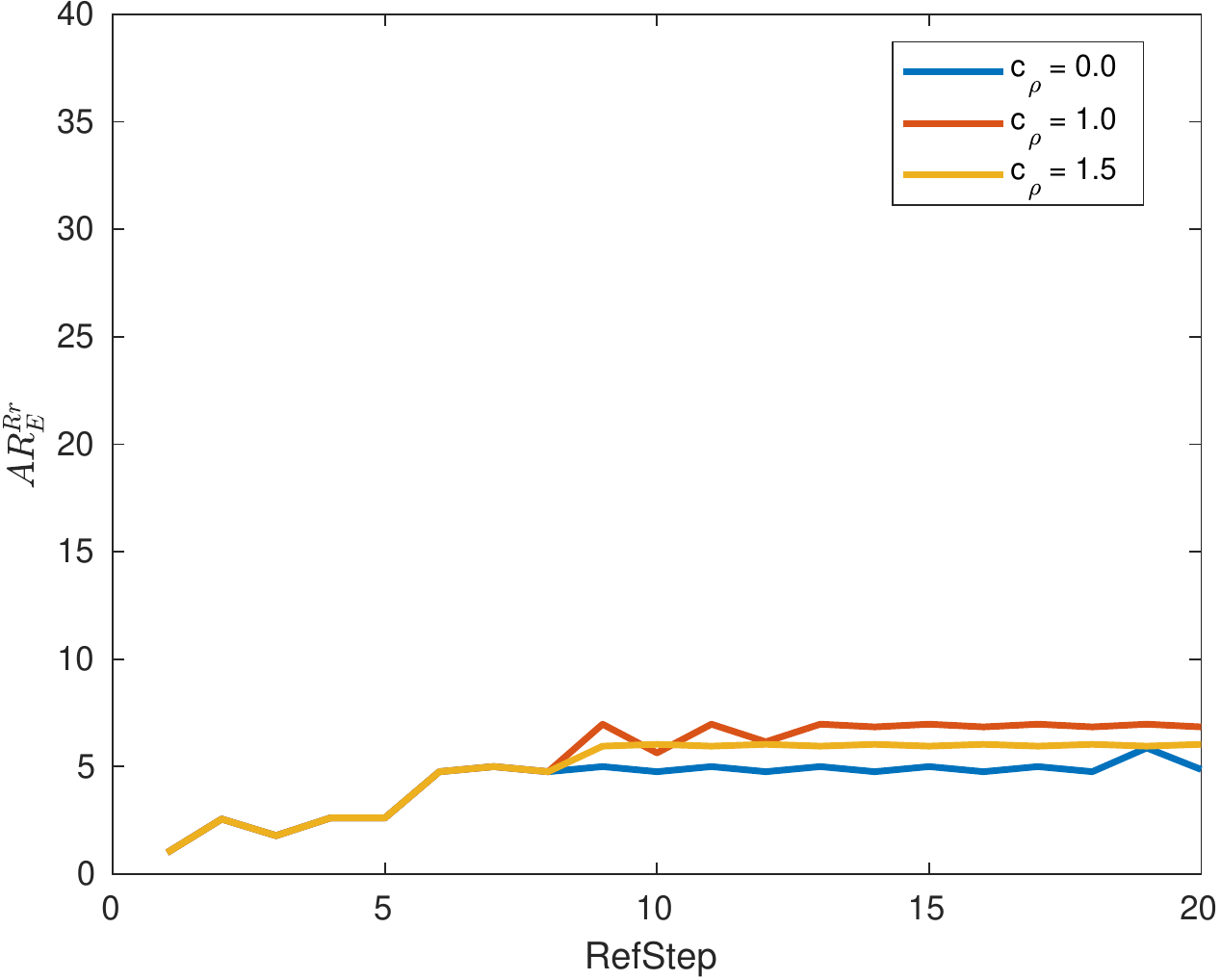}
    \caption{Max $AR^{Rr}_E$}
        \label{fig:p1_maxARR}
  \end{subfigure}
  \hfill
  \begin{subfigure}[b]{0.32\linewidth}
    \includegraphics[width=\linewidth]{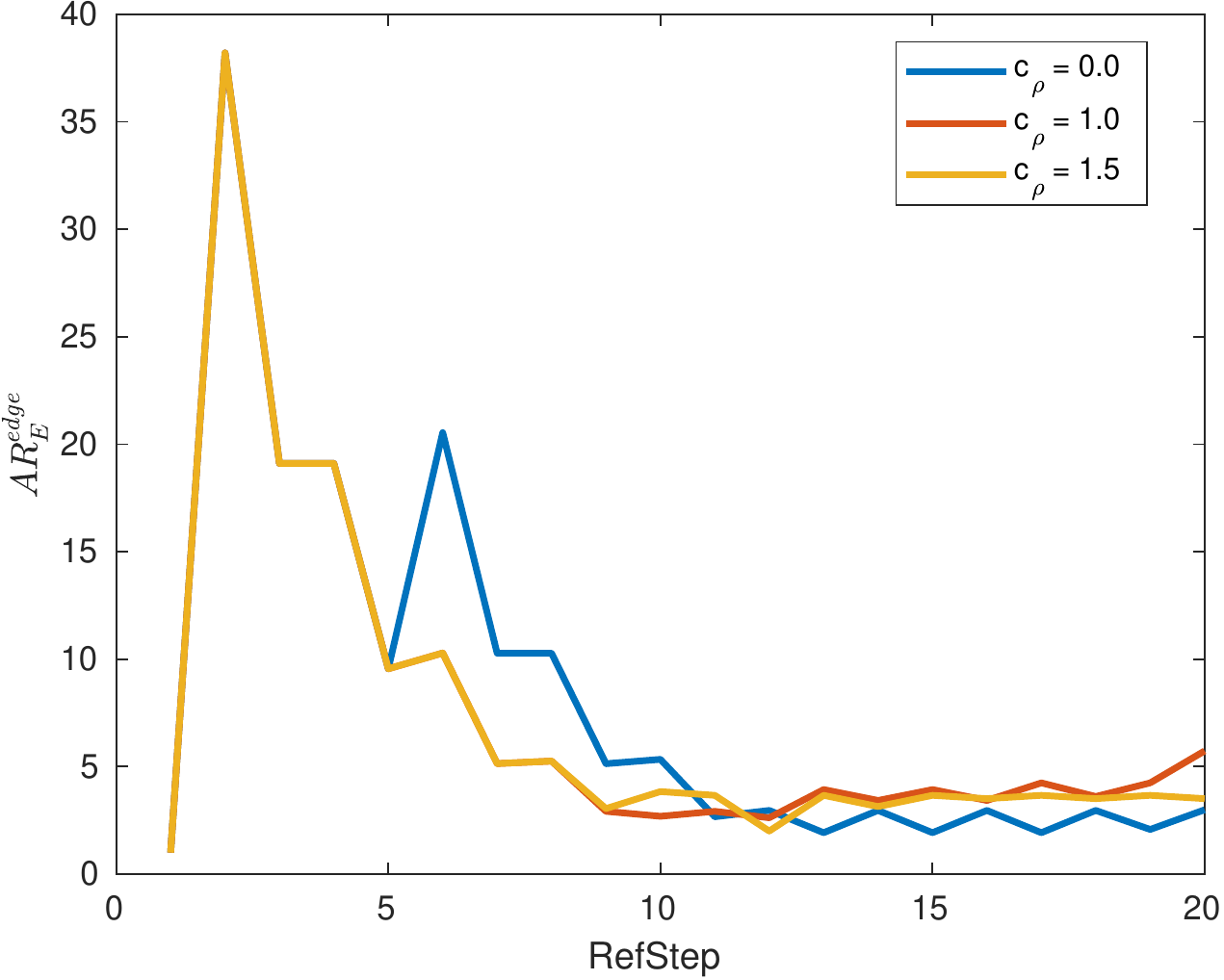}
    \caption{Max $AR^{edge}_E$}
    \label{fig:p1_maxARH}
  \end{subfigure}
  \hfill
  \centering
  \begin{subfigure}[b]{0.32\linewidth}
    \includegraphics[width=\linewidth]{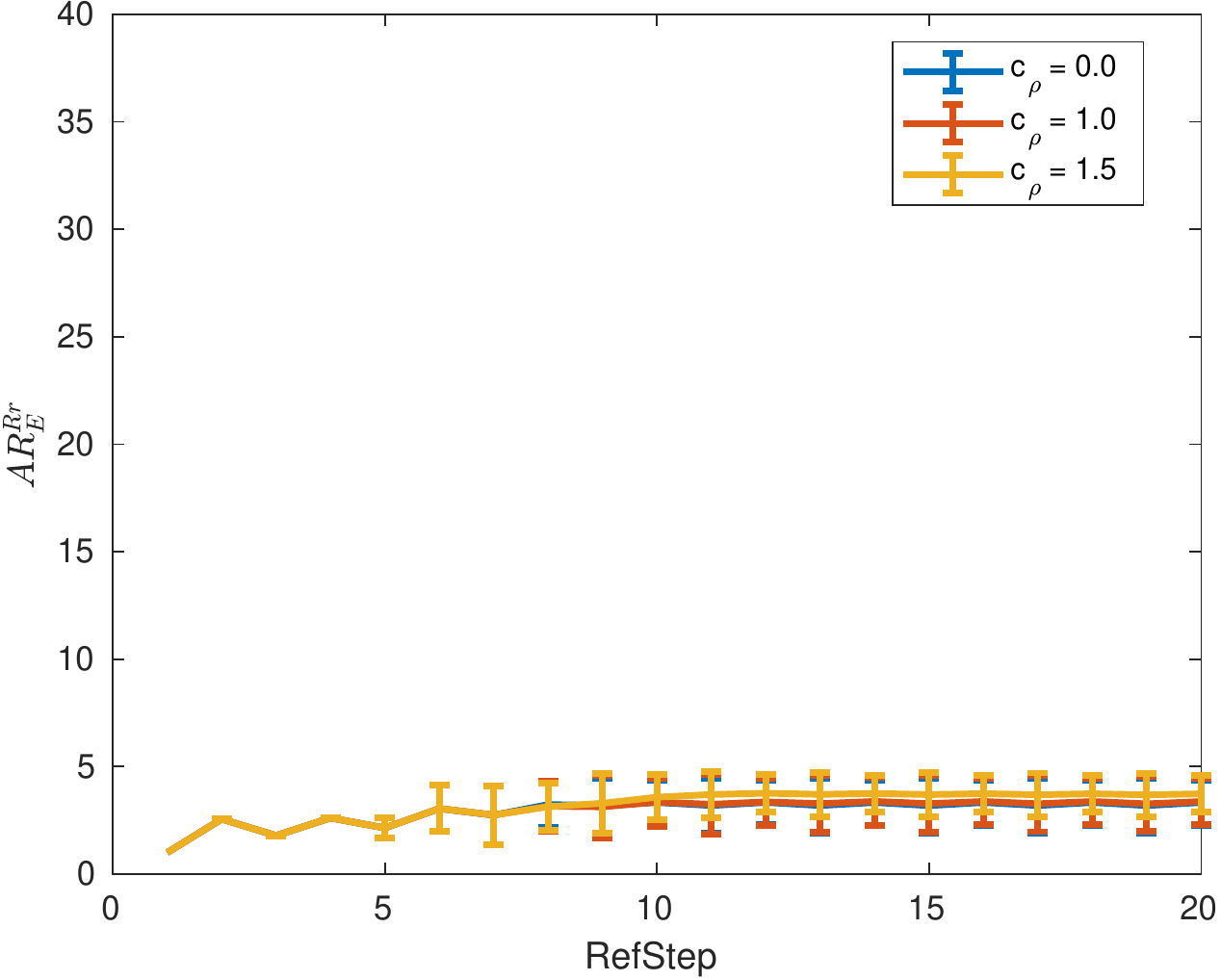}
    \caption{Mean $AR^{Rr}_E$}
    \label{fig:p1_meanARR}
  \end{subfigure}
  \hfill
  \begin{subfigure}[b]{0.32\linewidth}
    \includegraphics[width=\linewidth]{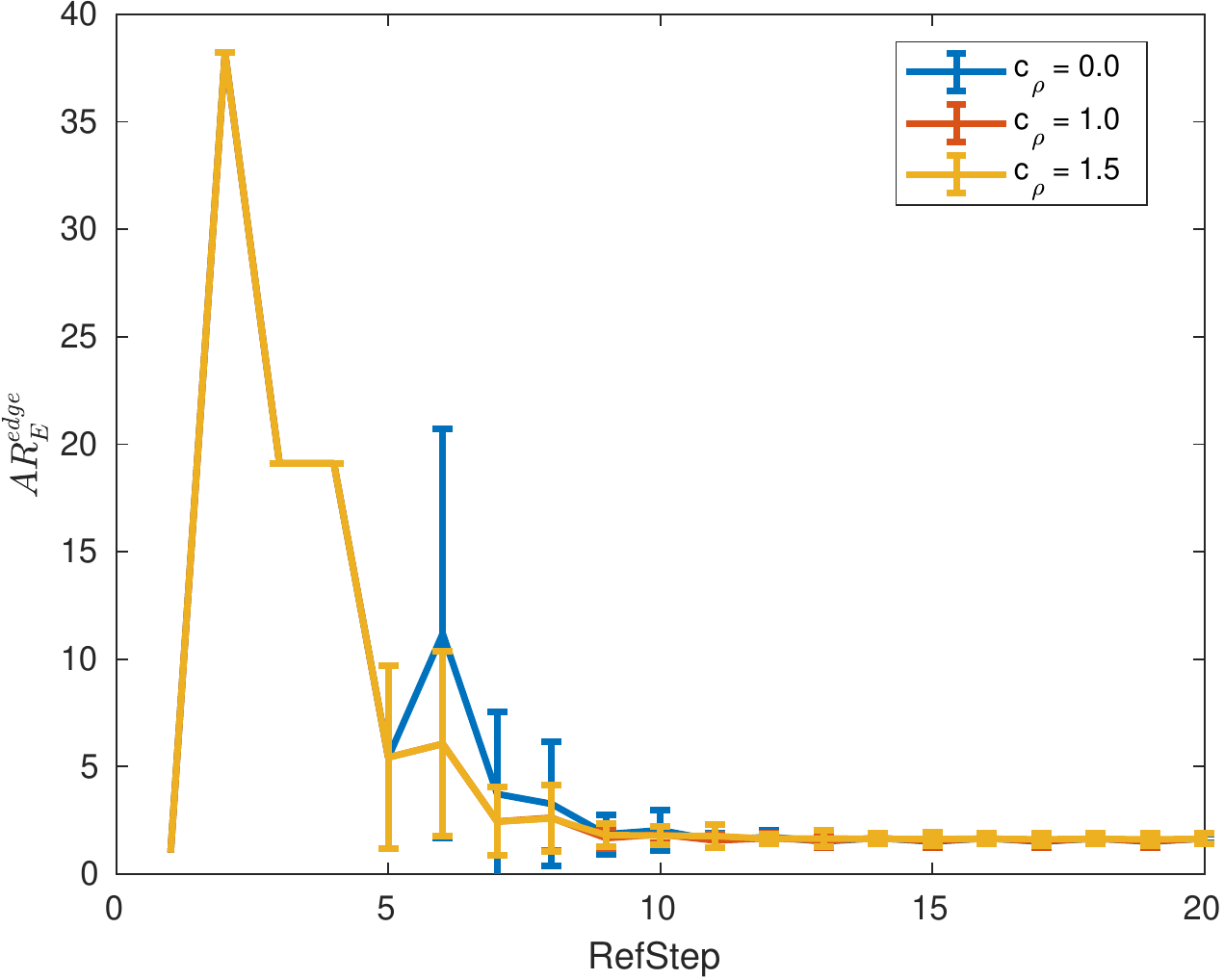}
    \caption{Mean $AR^{edge}_E$}
    \label{fig:p1_meanARH}
  \end{subfigure}
  \hfill
  \begin{subfigure}[b]{0.32\linewidth}
    \includegraphics[width=\linewidth]{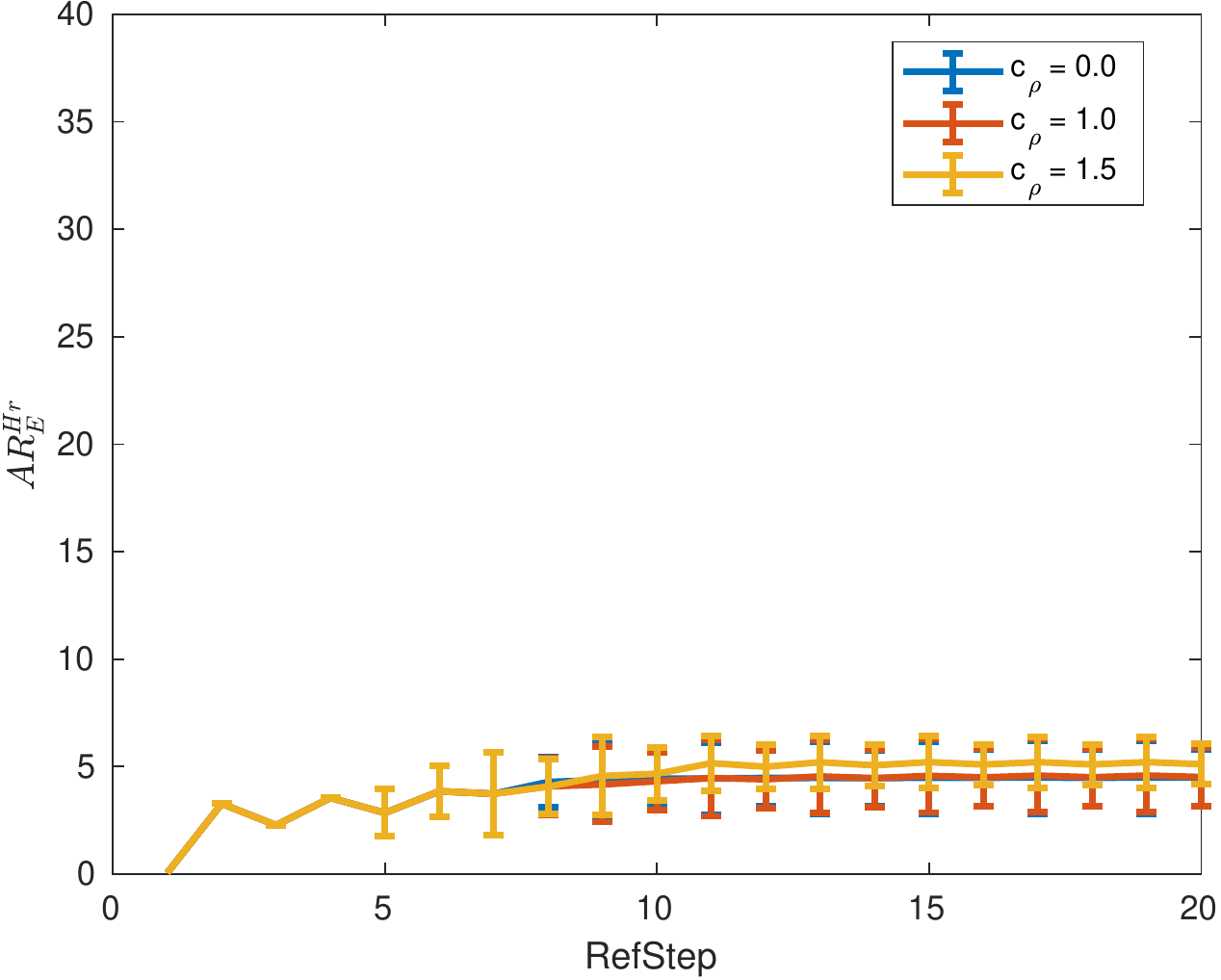}
    \caption{Mean $AR^{Hr}_E$}
    \label{fig:p1_MeanGamma}
  \end{subfigure}
  \hfill
  \begin{subfigure}[b]{0.32\linewidth}
    \includegraphics[width=\linewidth]{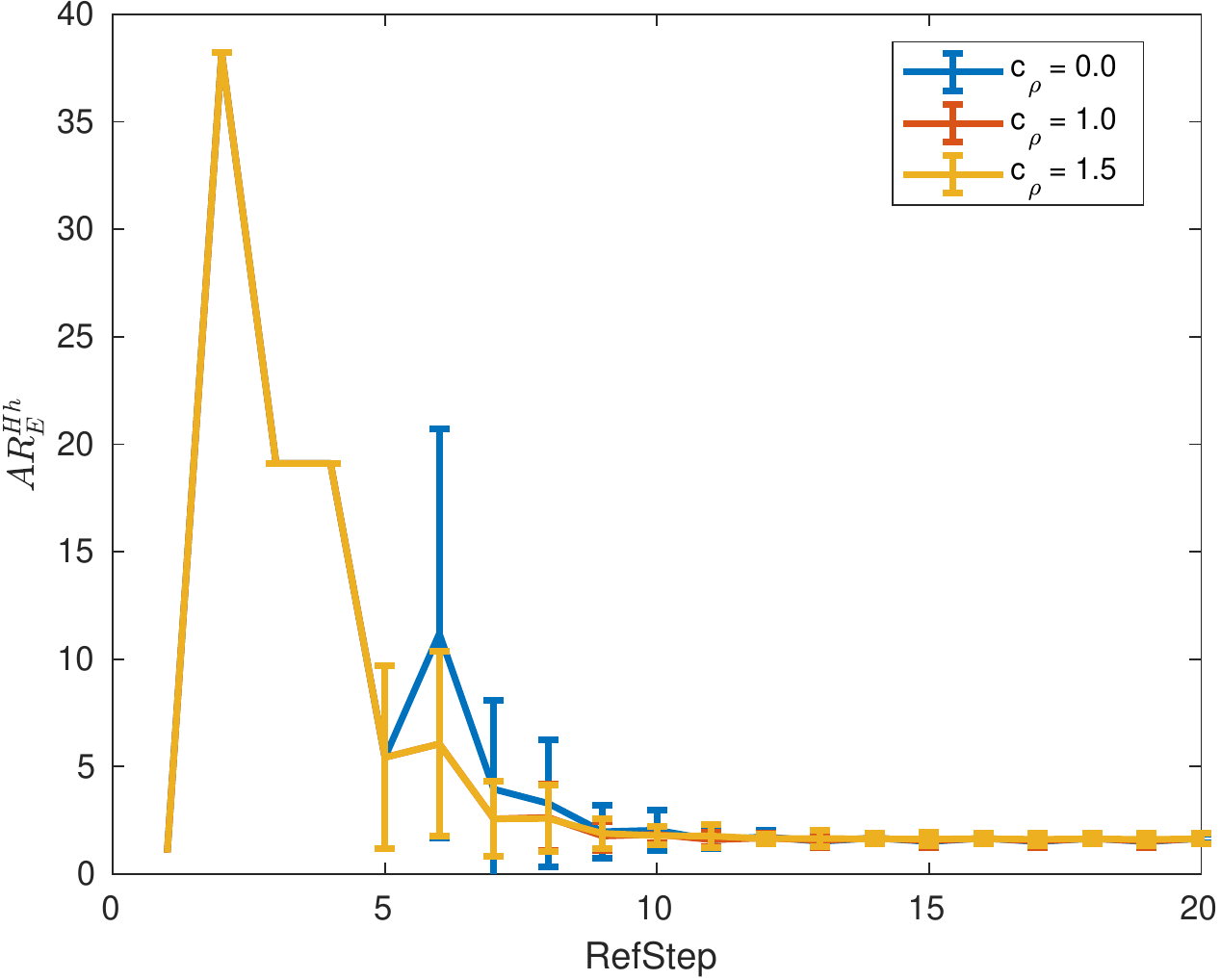}
    \caption{Mean $AR^{Hh}_E$}
    \label{fig:p1_MeanEta}
  \end{subfigure}
  \caption{P120, uniform refinement: $AR^{Rr}_E, AR^{edge}_E, AR^{Hr}_E,AR^{Hh}_E$}
  \label{fig:p1_Polygon-AnalysisMesh}
\end{figure}

\begin{figure}[t]
  \begin{subfigure}[b]{0.4\linewidth}
    \includegraphics[width=\linewidth]{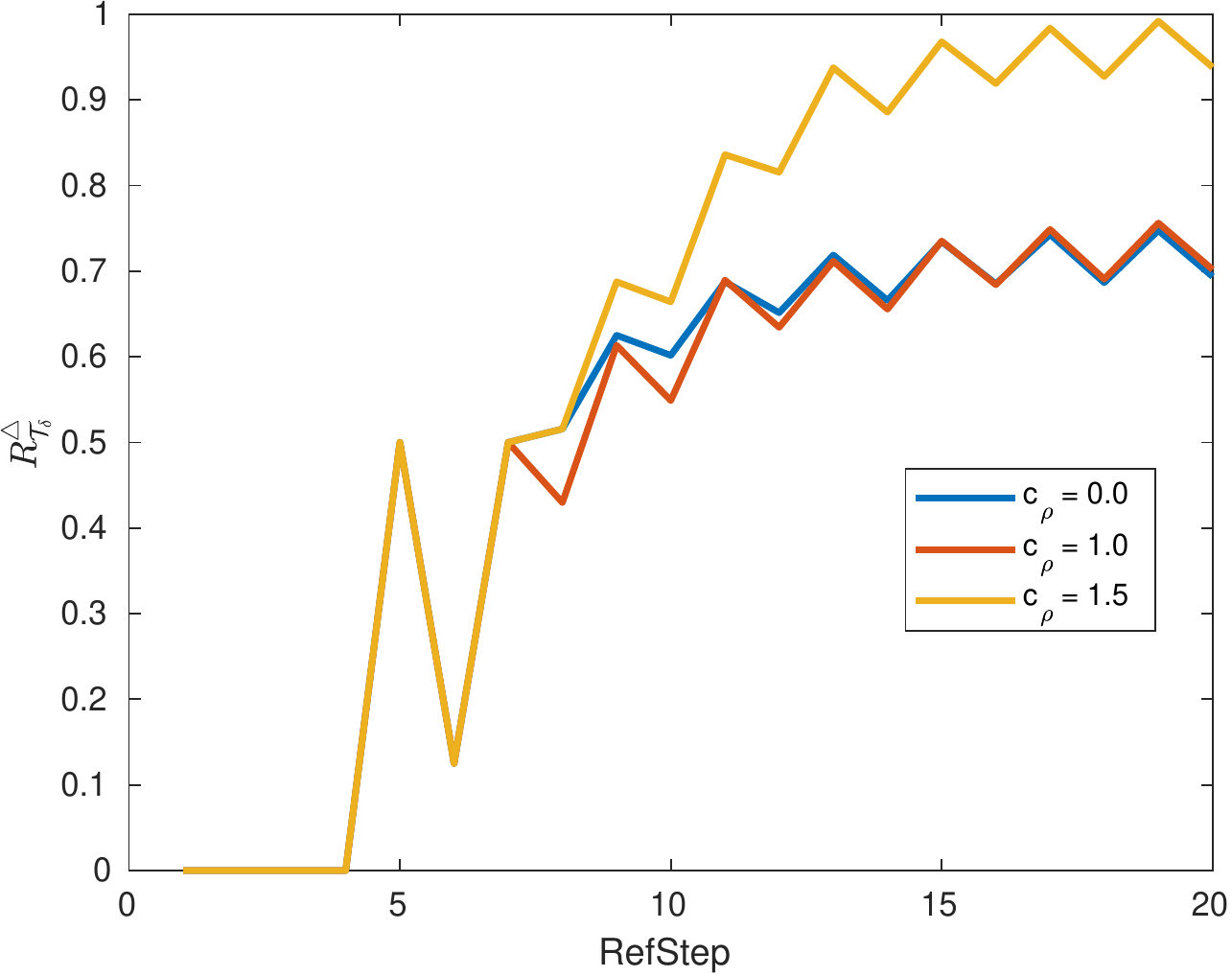}
    \label{fig:p1_Triangles}
  \end{subfigure}
  \hfill
    \begin{subfigure}[b]{0.4\linewidth}
    \includegraphics[width=\linewidth]{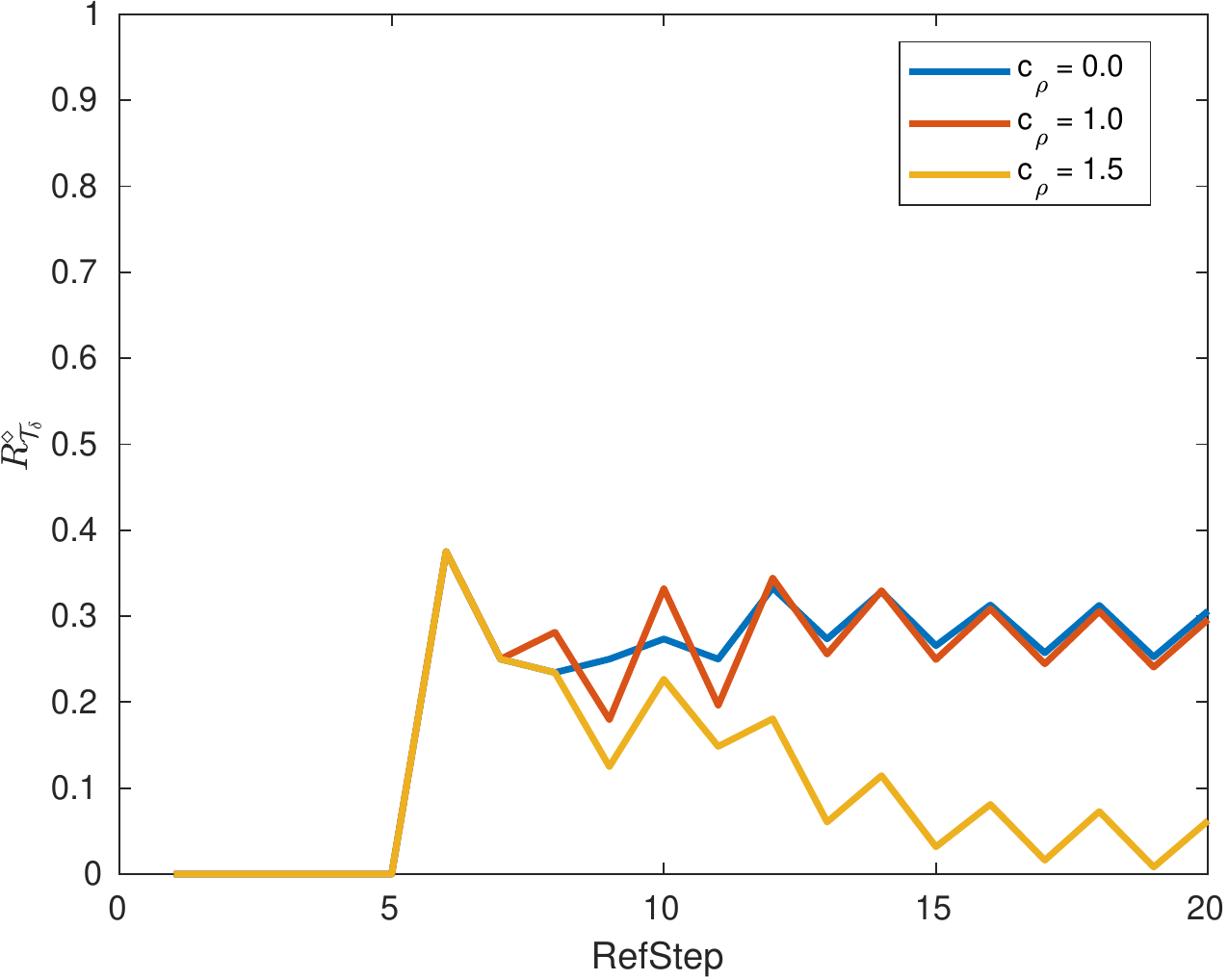}
        \label{fig:p1_Squares}
  \end{subfigure}
 \caption{P120, uniform refinement:  $R^\triangle_{\Th}$ (left) and $R^\diamond_{\Th}$ (right).}
 \label{fig:p1_Polygon-AnalysisMesh2}
\end{figure}

In this Section we analyse the behaviour of the meshes obtained by a uniform refinement of a regular polygon with 120 edges, initially not partitioned (one cell), with the strategy of the Algorithms~\ref{alg:MaximumMomentum} and \ref{alg:Refinement}.
In the first refinement iterations a strong degeneration of the elements quality is unavoidable starting from a regular polygon that is cut in two cells by a diagonal, nevertheless the algorithm is capable to fix the element quality and converge towards a good quality mesh composed by triangles and quadrilateral elements depending on the parameter $c_\rho$.

In Figure~\ref{fig:p1_Polygon-AnalysisMesh} we collect the behaviour of the parameters we have chosen for the characterization of the mesh quality.

Figures~\ref{fig:p1_maxARR},  \ref{fig:p1_meanARR} suggest that the parameter $AR_E^{Rr}$ is less sensitive to the geometric degeneration of the elements in the first step with respect to $AR_E^{edge}$, Figures~\ref{fig:p1_maxARH},  \ref{fig:p1_meanARH}. The strong jump in $AR_E^{edge}$ in the first iteration is due to the presence of a long edge (the diameter of the polygon) close to one of the 120 edges of the polygon. The behaviour of the parameters $AR_E^{Hr}$ is similar to the one of the parameter $AR_E^{Rr}$, and $AR_E^{Hh}$ is similar to  $AR_E^{edge}$.

In Figure~\ref{fig:p1_Polygon-AnalysisMesh2} (left) we report the fraction of cells that have three vertices and in Figure~\ref{fig:p1_Polygon-AnalysisMesh2} (right) the fraction of cells with four vertices. As noticed in the previous tests the refinement strategy characterized by the parameter $c_\rho=1.5$ clearly tends to naturally converge to a good quality triangular mesh and in all the cases, after few iterations, the mesh tends to be composed by elements with three or four vertices and all the cells can be definitely defined good quality cells.


%% file: conclusions.tex
\section{Conclusions}
\label{sec:conclusions}
In this paper we propose a refinement algorithm for polygonal meshes with convex cells.
We have tested its behaviour with an adaptive mesh refinement algorithm for virtual element methods (VEM), but it can be applied to any other polygonal method such as Discontinuous Galerkin methods, Hybrid High Order methods or Mimetic Finite Differences.
The algorithm requires the selection of a suitable splitting direction for each element, the proposed strategies are based on the maximum moment direction or the longest diagonal direction that perform in a comparable way.

With respect to the algorithm proposed in \cite{BBD} the main improvement is a strategy to avoid multiple splitting of an edge in the same refinement step that are often responsible of a degeneration of the quality of the elements. Moreover, the algorithm always avoids the creation of relative small edges during the refinement.

A set of geometrical parameters, easily computable, is tested in order to characterize the quality of the elements. Among the tested parameters we have found that the parameters $AR_E^{Rr}$ and $AR_E^{edge}$ can provide a measure of the quality of an element in all the situations.
An additional important property in the VEM discretization is related to an upper bound of the ratio between the number of vertices of the mesh over the number of cells ($E^{\#P}_{\#{\Th}}$), in particular, the proposed algorithm allows the increasing by one of the number of the vertices of the polygons produced by the refinement in a very particular unusual situation, in all the other situations the number of vertices of the produced polygons is less or equal to the number of vertices of the refined polygon and in the most common situations this number decreases if the number of vertices is larger than four. Depending on the parameter $c_\rho$ we have verified that, starting from general polygons, the mesh tends toward a good quality triangular, quadrilateral or mixed triangular-quadrilateral mesh. These meshes are characterized by a better efficiency (smaller $E^{\#P}_{\#{\Th}}$) with respect to other meshes of convex elements.

A general advantage of a conforming polygonal method that allows aligned edges is that it does not require a conformity recovery after an adaptive refinement iteration. This property trivially bounds the number of cells produced in each refinement iteration with the number of cells marked for refinement that is a property useful for optimality analyses of an adaptive mesh refinement approach.


The advantage of the proposed refinement method relies on the fact that the generation of an initial polygonal mesh can be usually replaced by a description of the domain geometry in terms of very coarse convex cells only considering the geometric constraints. This given coarse-mesh/geometry-description can be refined by an adaptive method producing a good quality mesh and a reliable solution fully driven by the problem.
